\input amstex\documentstyle{amsppt}  
\pagewidth{12.5cm}\pageheight{19cm}\magnification\magstep1
\topmatter
\title Comments on my papers\endtitle
\author G. Lusztig\endauthor
\address{Department of Mathematics, MIT, Cambridge MA 02139}\endaddress
\thanks{Supported in part by NSF grant DMS-1566618}\endthanks
\abstract{Changes have been made to the comments on \cite{24},\cite{37},\cite{92}, \cite{95}.
Comments to \cite{7},\cite{77},\cite{227},\cite{228},
\cite{229},\cite{231},\cite{232},\cite{233},\cite{234},\cite{235},\cite{236},\cite{237},\cite{238},\cite{239},
\cite{240},\cite{241},\cite{242},\cite{243},\cite{244},\cite{245},\cite{246},\cite{247},\cite{248},\cite{249},
\cite{250},\cite{251},\cite{252},\cite{253},\cite{254},\cite{255} have been added.}\endabstract 
\endtopmatter   
\document

\define\Irr{\text{\rm Irr}}

\define\mpb{\medpagebreak}

\define\sqc{\sqcup}

\define\lb{\linebreak}

\define\op{\oplus}
   
\define\part{\partial}

\define\iy{\infty}
\define\m{\mapsto}
\define\do{\dots}

\define\bsl{\backslash}

\define\lra{\leftrightarrow}

\define\sub{\subset}    

\define\T{\times}
\define\ti{\tilde}
\define\nl{\newline}
\redefine\i{^{-1}}

\define\un{\underline}

\define\ot{\otimes}

\define\ind{\text{\rm ind}}
    
\define\res{\text{\rm res}}

\define\tr{\text{\rm tr}}

\define\card{\text{\rm card}}

\define\a{\alpha}

\redefine\c{\chi}
\define\g{\gamma}

\define\e{\epsilon}

\define\r{\rho}

\redefine\t{\tau}

\define\k{\kappa}
\redefine\l{\lambda}

\redefine\G{\Gamma}

\define\Om{\Omega}

\redefine\L{\Lambda}

\define\kk{\bold k}

\define\BB{\bold B}
\define\CC{\bold C}

\define\HH{\bold H}

\define\NN{\bold N}

\define\QQ{\bold Q}

\define\ZZ{\bold Z}

\define\cb{\Cal B}

\define\ce{\Cal E}

\define\cg{\Cal G}
\define\ch{\Cal H}
\define\ci{\Cal I}

\define\cl{\Cal L}

\define\co{\Cal O}

\define\cs{\Cal S}

\define\fg{\frak g}

\define\fl{\frak l}

\define\fs{\frak s}

\define\fC{\frak C}

\define\tc{\ti c}

\define\sha{\sharp}

This document contain comments on some of my papers. I hope to add to it more
comments in future versions.

\head{7} A property of certain non-degenerate holomorphic vector fields, 1969\endhead
Let $X$ be a compact complex manifold and let $\xi$ be a holomorphic vector field
on $X$ with zero set $X_0$. For any $x\in X_0$, $\xi$ defines an
endomorphism $\t_x$ of the tangent space of $X$ at $x$. Let $[\t_x]$ be the multiset
of eigenvalues of $\t_x$. One says that $\xi$ is nondegenerate if for any $x\in X_0$,
$[\t_x]$ does not contain $0$; this implies that $X_0$ is finite. In this paper it is
conjectured that the multiset $\cup_{x\in X_0}[\t_x]$ is a union of sets of the form $\{z,-z\}$.
The conjecture is proved in this paper in the case where $\dim X=2$ using the
Atiyah-Bott fixed point formula and the non-vanishing of Bernoulli numbers
with even index. The general case was proved later in \cite{9}.

\head \cite{8} (with J.Milnor and F.P.Peterson) Semicharacteristics and
cobordism, 1969\endhead
I did the work on this paper during a two months stay in Oxford (fall of 1968). 
During my first meeting with Atiyah, he and Singer explained to me the following 
question on the (Kervaire) semicharacteristic. 
A compact smooth oriented manifold $M$ of dimension $4n+1$ has a 
semicharacteristic $c(M,p)=\sum_{i\in[0,2n]}\dim H^i(M,k)\mod 2$ with respect 
to a field $k$ of characteristic $p\ge0$. At the time it was known that
the obstruction to the existence of two independent vector fields on $M$ is
equal to $c(M,2)$ if $M$ is spin [E. Thomas, Bull. Amer. Math. Soc. 1969] and 
to $c(M,0)$, without assumption (Atiyah-Singer). The question was whether 
there is a relation between $c(M,0)$ and $c(M,2)$ which would make clear that the 
Atiyah-Singer result implies the Thomas result. (A few months earlier they 
asked F.Peterson at MIT the same question.) Then I and (independently) 
F.Peterson and J.Milnor found (different) 
formulas for $c(M,0)-c(M,2)$; one of those formulas expressed
$c(M,0)-c(M,2)$ as the Stiefel-Whitney number $w_2w_{4n-1}$ which 
clearly vanishes
for spin manifolds. The initial proofs of both formulas used C.T.C.Wall's results on the structure 
of the oriented cobordism group (that is the formulas were checked on the generators of 
that group) but in the final version Wall's results are not used.
The result of this paper is used 
in [Atiyah,Singer, Ann.Math. 1971] and is generalized in \cite{10}. 

\head \cite{9} Remarks on the holomorphic Lefschetz formula, 1969\endhead
I did the work on this paper during a two months stay in Oxford (fall of 1968).
The main result of the paper is the following rigidity result: the holomorphic
Lefschetz number of a circle action on a compact complex manifold is constant.
(This is of interest only in the case where the manifold is non-Kaehler, when
the induced action on the Dolbeault cohomology may be nonconstant.) The 
argument that I used in this paper has played a role in the proof of the
vanishing of the $\hat A$-genus of a $4k$-dimensional
 spin manifold with nontrivial circle action given in [Atiyah,Hirzebruch, in 
Essays on topology and related topics, 1970], where there is a reference to my 
result (but not to the paper itself). There is a similar reference in: 
[Bott and Taubes, Jour. Amer. Math. Soc. 2(1989)]. 

\head \cite{10} (with J.Dupont) On manifolds satisfying  $w_1^2=0$, 1971\endhead
This paper was written in 1970 during my stay (1969-71) at IAS.
It contains a generalization of the result on semicharacteristic in \cite{8}
 where the orientability assumption $w_1=0$ is weakened to $w_1^2=0$ 
($w_1$ is the first Stiefel-Whitney class). This was used in the paper 
[Davis and Milgram, Trans. Amer. Math. Soc., 1989].  
The appendix of this paper is a study of the symmetric power $SP^nX$ where 
$X$ is a compact unorientable smooth $2$-manifold whose first rational Betti 
number is $g$. We show that $SP^nX$ is a $(2n-g)$-dimensional bundle over a 
$g$-dimensional realtorus with fibre $RP^{2n-g}$, the real projective space 
of dimension $2n-g$. In particular if
$X$ is a projective plane ($g=0$) then $SP^nX=RP^{2n}$; if $X$ is a 
Klein bottle ($g=1$) then $SP^nX$ is a $RP^{2n-1}$-bundle over the circle. 
($SP^nY$ for $Y$ a compact Riemann surface
was studied earlier in [Macdonald, Topology 1962].) We also show that 
$S^{\iy}X$ is a product of a $g$-dimensional torus with 
$RP^{\iy}$. Our result $SP^n(RP^2)=RP^{2n}$ is reinterpreted
in [Arnold, Topological content of the Maxwell theorem on multiple 
representations of spherical functions, Topological methods in nonlinear analysis 7(1996)].
This paper contains also a study of a cobordism ring $G_*$ based on 
closed manifolds together with an element $\Gamma$ in $H^1(M,\ZZ/4\ZZ)$ 
which reduces to the first Stiefel Whitney class by reduction $\mod2$. 
Our explicit determination of this ring relies on earlier work of
[C.T.C.Wall, Ann. Math. 1960]; a special role in our description is played by 
the manifolds $SP^nX$ where $X$ is a Klein bottle and $n$ is a power of $2$. 
(There are two natural choices for $\Gamma$ but they represent the same element
in $G_*$.) 

\head \cite{11} Novikov's higher signature and families of elliptic operators, 1972\endhead
This paper was written in 1970 during my stay (1969-71) at IAS. I used it as 
my Ph.D. thesis at Princeton University (may 1971).
The main contribution of this paper is to introduce the analytic approach
(based on the index theorem) to attack the Novikov's conjecture on higher
signature. That conjecture states that, if one multiplies the
Hirzebruch $L$-class of a compact oriented manifold $M$ with a cohomology class
which comes from the cohomology of the classifying space of the fundamental
group of $M$ and then one integrates the result over $M$, one obtains a homotopy
invariant of $M$.
In this paper I introduce, for $M$ of even dimension with fundamental group $Z^n$,
a family of elliptic operators on $M$. These operators are obtained by twisting 
the Atiyah-Singer signature operator by a variable flat vector bundle
on $M$ coming from a unitary one dimensional representation of the fundamental
group. While the index of each of these operators is the same as that of the
untwisted operator it turns out the family of operators has an interesting
index in the K-theory of the parameter space of the space of flat bundles
considered and I showed that this index is on the one hand
a homotopy invariant and on the other hand from it one can recover
the whole Novikov higher signature thus proving Novikov's conjecture in this case.

Another contribution of this paper is to formulate a version of the
Hirzebruch signature theorem in which cohomology is taken with coefficients in 
a flat vector bundle with a flat hermitian form which is not necessarily
positive definite. In this case Hirzebruch's original proof (with constant
coefficients) does not work but the Atiyah-Singer theorem can be used
instead. In the paper I show that from this "twisted" signature theorem one
can derive various examples where Novikov's conjecture holds for
certain nonabelian fundamental groups.
The analytic approach of this paper has been extended by Mischenko and
Kasparov to the case where the fundamental group is a discrete subgroup of Lie
groups and then by Connes, Moscovici, Gromov, Higson and others to even more
general fundamental groups. See [Ferry,Ranicki,Rosenberg: Novikov signatures,
index theory and rigidity, London Math.Soc.Lect.Notes, 1995] for a review of
these developments. 
The twisted signature theorem of this paper is used in: [Gromov and Lawson, 
Ann. Math. 1980], [Atiyah, Math. Annalen 1987], [Gromov, in 
"Functional analysis on the eve of the 21st century, II", Progress in Math. 
132,Birkhauser 1996]. 

\head \cite{14}. Introduction to elliptic operators, 1974\endhead
This (mainly expository) paper is based on a lecture that I gave at a Trieste summer school in 1972.
It contains the definition of elliptic operators and their index. The only part which is perhaps non-standard 
is the definition of analytical index as a homomorphism $ind:K(BT^*M,ST^*M)@>>>\ZZ$ where $T^*M$ 
is the cotangent bundle of a compact manifold $M$, $BT^*M$ is its unit disk bundle, $ST^*M$ is its unit sphere bundle and 
$K()$ is $K$-theory. The usual (Atiyah-Singer) definition of $ind$ is via the theory of pseudo-differential operators. But in 
this paper I show that if we are willing to increase $\ZZ$ to $\ZZ[1/2]$, one can define $\ind:K(BT^*M,ST^*M)@>>>\ZZ[1/2]$
in a more elementary way, using only differential operators.

\head \cite{17}. On the discrete series representations of the classical groups over a finite field, 1974\endhead
This paper represents my talk at the ICM in August 1974. I was originally invited to give a talk
in the Algebraic Topology section but I requested to change to the Lie Groups section.
In 1973 the only cases where cuspidal characters of a reductive group over a finite field were
constructed were $GL_n(F_q)$ (by J.A.Green), $Sp_4(F_q)$ (by B.Srinivasan) and $G_2(F_q)$ (by B.Chang, R.Ree).
In 1973, after my study [13] of the "Brauer lifting" of the standard n-dimensional representation of $GL_n(F_q)$,
I tried to find the constituents of the Brauer lifting $X$ of the standard representation (of dimension $N$) of a 
symplectic orthogonal or unitary group $G(F_q)$. The result that I found is that $X=X_1+X_2+...+X_N$
where $X_i$ is $\pm$ an irreducible representation, $X_1$ is up to sign a cuspidal representation 
(new at the time) of dimension $|G(F_q)|$ divided by the order of a "Coxeter torus"
and by the order of a maximal unipotent subgroup; moreover $X_i$ for $i>1$ were noncuspidal and could be
explicitly described as components of certain induced representations from analogous cuspidal representations
of Coxeter type of smaller classical groups or $GL_n$ by determining explicitly the relevant Hecke algebras.
Thus this method gives a way to approach at least the "Coxeter series" of cuspidal representations of a classical
group. This is what is explained in the first part of this paper. (The proofs of the results in the first part
were never published since they were superseded by later developments.)

In the second part of the paper I described my joint work with Deligne (done
during the spring 1974 at IHES) in which
$l$-adic cohomology is used to construct representations of $G(F_q)$ where $G$ is a connected reductive
group over $F_q$. This method was first used by Tate and Thompson [Tate, Algebraic cycles and poles...,1965]
who observed that the obvious action of the unitary group $U_3(F_q)$ on the projective
Fermat curve $x^{q+1}+y^{q+1}+z^{q+1}=0$ over $\bar F_q$ induces an action on $H^1$ which is the (interesting)
irreducible representation of degree $q^2-q$ of $U_3(F_q)$. Around 1973, Drinfeld observed that the 
cuspidal representations of $SL_2(F_q)$ can be realized in the cohomology with compact support of the 
Dickson curve $xy^q-x^qy=1$ over $\bar F_q$ by taking eigenspaces of the action of $T=\{t\in\bar F_q^*;t^{q+1}=1\}$ 
which acts on the curve by homothety. (I learned about this fact from T. A. Springer in 1973.) (Note that
Dickson's curve can be viewed as the (open) part of the Fermat curve where $z\ne0$.
This open part is stable under $SL_2(F_q)\T T$ viewed as a subgroup of $U_3(F_q)$.)
The main result of this section was the introduction for any element $w$ in the Weyl group of $G$, of two
new algebraic varieties: the variety $X_w$ of Borel subgroups $B$ of $G$ such that $B$ and its image under the
Frobenius map $F$ are in relative position $w$ and the finite principal covering $\tilde X_w$ of $X_w$ whose group
is the group $T_w$ of rational points of an $F$-stable torus of type $w$ of $G$. (Note 
that the Drinfeld curve is a special case of $\tilde X_w$ in the case $G=SL_2$ 
and the Tate-Thompson curve is the compactification of an $X_w$ in the case
$G=GL_3$ with a nonsplit $F_q$-structure.)
These varieties admit natural action of $G(F_q)$ and the principal covering above gives rise to
$G(F_q)$-equivariant local systems on $X_w$, one for each character $\theta$ of $T_w$. By passing to cohomology with 
compact support with coefficients in such a local system one obtains representations of $G(F_q)$. By taking 
alternating sums one obtains certain virtual representations $R(w,\theta)$ of $G(F_q)$ indexed by the various $\theta$.
At the time when this paper was written (summer 1974) we conjectured that $R(w,\theta)$ for $\theta$ generic
are up to sign irreducible characters which provide a solution of the Macdonald conjecture. 

\head \cite{18} Sur la conjecture de Macdonald, 1975\endhead
By my joint work with Deligne (spring 1974) described in \cite{17}, a conjectural solution to the
Macdonald conjecture for the irreducible representations of $G(F_q)$ was known in terms of the 
the virtual representations $R(w,\theta)$ defined by the subvarieties $X_w$ of the flag manifold, their finite
coverings $\tilde X_w$ and their cohomology with compact support. But it was not clear how to prove the 
irreducibility of $R(w,\theta)$ for $\theta$ generic or how to compute the degree of $R(w,\theta)$.
In the fall 1974 (at Warwick) I completed the proof of the fact that the virtual representations 
$R(w,\theta)$ defined in \cite{17} are indeed a solution to the Macdonald conjecture. In this paper (written
in late 1974) I sketched this proof; it is based on the following principle:

($*$) Assume that $H$ is a finite group acting on an algebraic variety $Y$ in such a way that the space of orbits
$Y/H$ is again an algebraic variety $Y'$; assume further that there is a partition of $Y$ into finitely many
locally closed $H$-stable pieces $Y_i$ and on each $Y_i$ the action of $H$ extends to an action of a connected
algebraic group $H_i$. Then $Y',Y$ have the same Euler characteristic.

The main observation of this paper is that ($*$) is applicable in the following two cases:

(A) $Y=\tilde X_w/P(F_q)$, $Y'=X_w/P(F_q)$ where $P$ is a parabolic subgroup defined over $F_q$;

(B) $Y=(\tilde X_w\T X_{w'})/G(F_q)$, $Y'=(X_w\T X_{w'})/G(F_q)$,
where $w,w'$ are Weyl group elements.
\nl
Now ($*$) in case (A) implies easily the expected formula for the degree of $R(w,\theta)$ and ($*$) in case (B) implies 
easily an explicit formula for the inner product of an $R(w,\theta)$ with an $R(w',\theta')$ from which the desired 
irreducibility result follows.
The proof of the degree formula and that of the inner product formula given later in \cite{22} are
quite different: they use a disjointness theorem.

\head \cite{19} Divisibility of projective modules of finite Chevalley groups by the Steinberg module, 1976 \endhead
This paper was written during my stay at IHES in the spring of 1974. The motivation for this paper
was to find evidence for the Macdonald conjecture for a connected reductive group $G$ defined over $F_q$.
(However, by the time this paper appeared, Macdonald's conjecture was proved.)
If $T$ is a maximal torus of $G$ defined over $F_q$ and $\theta$ is a character of $T(F_q)$ then the induced
representation $R=Ind_{T(F_q)}^{G(F_q)}(\theta)$ is defined. If $\theta$ is in general position and if
$R(T,\theta)$ is the irreducible representation of $G(F_q)$ provided by Macdonald's conjecture (assumed to hold)
then $R$ is isomorphic to $R(T,\theta)\otimes S$ where $S$ is the Steinberg representation. Therefore if we can
prove apriori that $R$ is isomorphic to $S$ tensor some virtual representation, then this would be evidence
for the Macdonald conjecture. Now $R$ can be viewed as a representation coming from a
projective $G(F_q)$-module over the ring of integers in a suitable $p$-adic field. Hence it would be enough to
show that any such projective $G(F_q)$-module is "divisible" by $S$. This is what is shown in this paper.

\head \cite{20} A note on nilpotent matrices of fixed rank, 1976\endhead
This paper was written in early fall 1974. At that time the series of
representations of a reductive group attached to a maximal torus over $F_q$
were already constructed (in the joint work with Deligne, see \cite{17}) but their
irreducibility was not yet proved (it was proved shortly afterwards in \cite{18}). 
But in the case of the even nonsplit orthogonal group over $F_q$ and the series 
corresponding to the Coxeter torus the character was explicitly computable
(in this case I could compute the Green functions, by some computations
which later became part of \cite{23}) hence in this case irreducibility could be
proved directly by using orthogonality of the explicit Green functions; to do
this it was necessary to know the number of unipotent elements $u$
such that $u-1$ has fixed rank. (It turned out out that the Green function
on a unipotent element $u$ depended only on the rank of $u-1$.)
This number is computed in this paper. The result of this paper gives a new 
proof (for classical groups over $F_q$) for 
Steinberg's theorem on the total number of unipotent elements. (For
$GL_n$ that theorem is due to Fine and Herstein and independently to
Ph.Hall.)

\head \cite{22} (with P.Deligne) Representations of reductive groups over finite fields, 1976\endhead
This paper (written during the first half of 1975)
contains a detailed study of the varieties $X_w,\tilde X_w$ associated in \cite{17} (by me and Deligne)
to an element $w$ in the Weyl group of a connected reductive group $G$ defined over a finite field $F_q$ and of
the associated virtual representations $R(w,\theta)$ of $G(F_q)$. Here $\theta$ is a character of the finite "torus"
$T_w$ of type $w$. In Section 3 a Lefschetz type fixed point formula for a transformation of finite order of an
algebraic variety is given. (This formula was already used implicitly in \cite{17}.) This is used in Section 4
to prove a formula for the character of $R(w,\theta)$ assuming that the Green functions of $G$ and smaller groups
are known. In Section 6 the disjointness theorem is proved: the virtual representations $R(w,\theta)$,
$R(w',\theta')$ are disjoint unless $\theta,\theta'$ are conjugate after extension of the ground field (and
composition with the trace); it is also shown that the equivalence classes of various $\theta$ as above can be 
viewed as semisimple conjugacy classes defined over $F_q$ in the "dual" group $G^*$ (at least if $G$ has
connected centre).  The proof of the disjointness uses the possibility of extending the action of a finite group 
on $(\tilde X_w\times\tilde X_{w'})/G(F_q)$ to the action of a higher dimensional group, not on the whole variety, 
but separately on each piece of a partition of the variety in pieces stable under the finite group (compare with
($*$) in the comments to \cite{18}).
The disjointness theorem has several applications (which were proved in a different way in \cite{18}):  the 
degree formula for $R(w,\theta)$; the inner product formula for $R(w,\theta),R(w',\theta')$; the orthogonality 
formula for Green functions.
In Corollary 7.7 it is shown that any irreducible representation of $G(F_q)$ appears in some 
$R(w,\theta)$; a completely different proof of this result based on the theory of perverse sheaves was later
given in \cite{178}. Combining this result with the disjointness theorem one obtains a canonical map from the set of
irreducible representations of $G(F_q)$ (up to isomorphism) to the set of 
semisimple conjugacy classes in $G^*(F_q)$ (at least if $G$ has connected centre); this is an initial but 
crucial step in the classification of irreducible representations of
$G(F_q)$. In Section 9 it is shown that $X_w$ is affine assuming that $q$ is greater than the Coxeter number.
Recently it has been shown that $X_w$ is affine without restriction on $q$ assuming that $w$ has minimal length in its
twisted conjugacy class [Orlik and Rapoport, J. Algebra, 2008] and [He, J. Algebra, 2008], see also 
[Bonnaf\'e and Rouquier, J. Algebra, 2008]. In 9.16 it is shown that any Green function 
evaluated at a regular unipotent element is equal to $1$. In Section 10 (assuming that the centre of G is connected)
it is shown how to parametrize explicitly the irreducible components of the Gelfand-Graev representations and that
these components are explicit linear combinations of $R(w,\theta)$. In Section 11 the results of the paper are
extended to Ree and Suzuki groups. In this case the inner product formula cannot be handled by the methods of
this paper in the case $q=\sqrt{2}$ or $q=\sqrt{3}$, but it can be handled by the proof given later in \cite{30}.

\head \cite{23} On the Green polynomials of classical grouyps, 1976\endhead
This paper was written in the summer of 1975, after the completion of \cite{22}. In \cite{22} a general study
of the variety $X_w$ (all Borel subgroups in relative position $w$ with their transform under Frobenius in
a reductive group $G$ over $F_q$) was made. In the present paper I tried to study in detail the first 
nontrivial class of examples of the variety $X_w$, namely the case where $G$ is a classical group and $w$ is
a Coxeter element of minimal length. In this case I obtained explicit formulas for (a) the number of rational points 
of $X_w$ over any extension of the ground field  and for (b) the Green functions
(alternating sums of traces of a unipotent element of $G(F_q)$ on the cohomology with compact support of $X_w$. 
In the case of symplectic groups, (b) solved a conjecture of B.Srinivasan. 
This paper was a preparation for my next project \cite{25} in which I studied $X_w$ for $w$ a (twisted) Coxeter element of 
minimal length for a general $G$.

\head \cite{24} On the finiteness of the number of unipotent classes, 1976\endhead
The work on this paper was started during a visit to IHES in December 1974 and was completed during another
visit to IHES in December 1975. Let $G$ be a connected reductive group defined over $F_q$, let $P$ be a parabolic
subgroup of $G$ (not necesarily defined over $F_q$) and let $L$ be a Levi subgroup of $P$ defined over $F_q$.
In this paper I define for any virtual representation $\r$ of $L(F_q)$ a virtual representation $R(L,P,\r)$ of $G(F_q)$.
In the case where $P$ is defined over $F_q$ this is the usual induced representation from $P(F_q)$ to $G(F_q)$ where
$\r$ is viewed as a virtual representation of $P(F_q)$. In the case where $L$ is a maximal torus defined over $F_q$ and
$\r$ is a character of $T(F_q)$, this reduces to the construction in \cite{22}.
 One of the main results of this paper 
is an inner product formula for two virtual representations $R(L,P,\r),R(L',P',\r')$ under a genericity assumption.
One consequence of this is that, under a genericity assumption, $R(L,P,\r)$ is irreducible (up to sign) for 
$\r$ irreducible. This irreducibilty result plays a key role in my later papers \cite{29},\cite{57}; it allows one to 
construct new irreducible representations starting with known irreducible representations of $L(F_q)$. Another
consequence is that the number of unipotent conjugacy classes in $G$ is finite, answering a conjecture of Steinberg
at the ICM in 1966. Previously this finiteness result was known only for classical groups or for exceptional 
groups in good characteristic. 
(In the paper I attribute the finiteness result in characterictic zero to Dynkin and
Kostant. But in fact this was known earlier: it follows by combining the 1942 paper of
Morozov with the 1944 paper of Malcev.)
In 1980 Mizuno gave another proof of the finiteness result for 
exceptional groups in bad 
characteristic based on extensive computations. I believe that the proof given in this paper is still the only 
proof of finiteness which does not use classification.
After this paper was written, Deligne stated a refinement of the inner product formula for $R(L,P,\r),R(L',P',\r')$
(without genericity assumptions) as an analogue of Mackey's theorem (generalizing the case where $L,L'$ are maximal
tori, known from \cite{22}) and proved it assuming that q is large. Deligne's proof remains unpublished. Around 1985
I proved a character sheaf version of this formula (see \cite{65}); in view of the results of [89] this implies the 
refined inner product formula for representations (again for large $q$). A version of Deligne's proof appeared in 
[Bonnaf\'e, J.Alg. 1998]. Recently [Bonnaf\'e and Michel, J. Alg., 2011] gave a proof of this formula with a
very mild assumption on $q$, using computer calculation. The general case is still not proved. 
The (refined) inner product formula would imply that $R(L,P,\r)$ is independent of $P$. 

\head \cite{25} Coxeter orbits and eigenspaces of Frobenius, 1976\endhead
The work on this paper was done in late 1975. This paper continues the project started in \cite{23} to
study in detail the variety $X_w$ of \cite{22} in the case where $w$ is a (twisted) Coxeter of minimal length in the
Weyl group of a connected reductive almost simple group $G$ defined over $F_q$. (The case of Suzuki and Ree groups
is also treated in the paper.) One of the main results 
of this paper is a construction of several new unipotent cuspidal representations of $G(F_q)$ in the case where $G$ is
exceptional. Let $d$ be the smallest integer $\ge1$ such that $X_w$ is stable under $F^d$, the $d$-th power of
the Frobenius map. In this paper I give

(a) an explicit computation of the eigenvalues of $F^d$ on $H^*_c(X_w)$ and an explicit formula for
the dimensions of its eigenspaces;

(b) a proof that $F^d$ acts semisimply on $H^*_c(X_w)$ and its eigenspaces are irreducible, mutually nonisomorphic $G(F_q)$-modules.
\nl
In the case where $G$ is of type $E_7$ we have $d=1$ and two of the eigenvalues of $F$ are of the form $\sqrt{-q^7}$, providing
the first examples in a split case where $F:H^*_c(X_w)@>>>H^*_c(X_w)$ can have eigenvalues with absolute value not an integer power of $q$. 
In the case where $G$ is a Suzuki or Ree group of type $B_2$ or $G_2$ with $q=\sqrt{Q}$ ($Q$ is an odd power 
of $m$, $m=2$ or $3$), the variety $X_w$ for $w$ a simple reflection is an affine curve and its compactification $\bar X_w$ is a
smooth projective curve defined over $F_Q$ such that $\bar X_w-X_w$ consists of $Q^m+1$ points. On the other hand by theorem 3.3(i) 
of this paper, $X_w$ has no $F_Q$-rational points. It follows that the number of $F_Q$-rational points of $\bar X_w$ ie equal to
$Q^m+1$. Note also that the genus of $\bar X_w$ is determined explicitly from (a) or \cite{22}.
In a letter to me dated May 11, 1983, J.-P.Serre made the following remarks.

(1) $\bar X_w$ has the maximum number of $F_Q$-rational points compatible with its genus.

(2) If a smooth curve over $F_Q$ has $Q^m+1$ rational points and has the same genus as $\bar X_w$ then it has the
same zeta function as $\bar X_w$ (which is determined from (a)).
\nl
Due to property (1), these curves have been used to produce Goppa (error correcting) 
codes. See [N. Hurt: Many rational points; coding theory and algebraic geometry, Kluwer, 2003].

\head\cite{29} Irreducible representations of finite classical groups, 1977\endhead
The work on this paper was started at Warwick during the summer of 1976 and completed
during my visit to MIT in the fall 1976. This paper contains the classification and
degrees of the irreducible complex representations of classical groups (with connected
centre) other than $GL_n$, over a finite field. It relies on: the use of the "cohomological
induction" \cite{22},\cite{23}; the use of the dimension formulas for the irreducible representations
of Hecke algebras of type $B$ with two parameters [Hoefsmit, UBC Ph.D. Thesis, 1974] (of which 
I learned from B.Chang during my visit to Vancouver at ICM-1974). This paper establishes 
what was later called "the Jordan decomposition" for the representations of classical 
groups (with connected centre). It also establishes the parametrization of unipotent 
representation for these groups in terms of some new combinatorial objects,
the "symbols" and the classification of unipotent cuspidal representations of classical
groups (for example $Sp_{2n}(F_q)$ has such a representation if and only if $n=k^2+k$ which 
is then unique). Also it is shown that the endomorphism algebra of the representation 
induced from a unipotent cuspidal (or more generally isolated cuspidal) representation 
to a larger classical group is an Iwahori-Hecke algebra (anticipating a later result of 
[Howlett and Lehrer, Inv. Math. 1980]) and giving also precise information on the values of 
the parameters of that Iwahori-Hecke algebra. For this we need to count in terms of 
generating functions the number of conjugacy classes in a classical group with connected 
centre. This together with an inductive hypothesis and the methods outlined above give a
way to predict the number of isolated cuspidal representations. The degrees of these 
isolated cuspidal representations can be guessed using the technique of symbols by "interpolation" from 
the degrees of noncuspidal representations. To prove that these guessed are correct we 
need to calculate the sum of squares of the (guessed) degrees of unipotent representations 
which is perhaps the most interesting part of this paper. To do this I find explicit 
formulas (for each irreducible representation $E$ of the Weyl group) of the polynomial $d_E(q)$ 
whose coefficients record the multiplicities of $E$ in the various cohomology spaces of the 
flag manifold. (I do this first for $GL_n$ and then reduce the case of classical groups to 
that of $GL_n$.) Then I show that the (guessed) degree polynomials can be expressed as linear 
combinations of the $d_E(q)$ with constant coefficients of the form plus or minus $1/2^s$. This 
anticipates the notion of family of representations of the Weyl group and the role of the 
nonabelian Fourier transform \cite{34} (which in this case happens to be abelian.) Here the use 
of the technique of symbols (introduced in this paper) is crucial. It is remarkable that 
suitable variations of the 
notion of symbol (used here in connection with unipotent representations) were later shown 
to be exactly what one needs to describe explicitly the Springer correspondence (including 
the generalized one) for classical groups \cite{59},\cite{61} and for classical Lie algebras in 
characteristic 2 [T. Xue, 2009]. 

\head \cite{30} Representations of finite Chevalley groups, 1978\endhead
This paper represents lectures that I gave in August 1977 at Madison, Wisconsin.
Let $G$ be a connected reductive group defined over $F_q$. 
Among other things, in this paper I give two refinements (see (a),(b) below)
of the inner product formula for the virtual representations $R(w,\theta)$ (see 
\cite{22}) of $G(F_q)$:

(a) a proof (in 2.3) which applies equally well to 
the Ree and Suzuki groups with $q=\sqrt{2}$ or $q=\sqrt{3}$ which were not covered by earlier 
proofs in \cite{18},\cite{22};

(b) a proof (see 3.8) of the fact that, if $w,w'$ are in the Weyl group $W$ and $X_w,X_{w'}$ are the varieties 
of \cite{22}, then $|((X_w \T X_{w'})/G(F_q))(F_{q^s})|$ is equal to the trace of a linear transformation 
$h\m T_whT_{w'{}^{-1}}$ of the Hecke algebra with parameter $q^s$. (Assume that $G$ is split over $F_q$.)
\nl
Note that (b) is a refinement of the inner product formula for $R(w,\theta)$ (in
the case where $\theta=1$) since for that formula one needs the Euler 
characteristic of $(X_w\T X_{w'})/G(F_q)$ which is the limit of 
$|((X_w\T X_{w'})/G(F_q))(F_{q^s})|$ as $s$ goes to $0$. I also show that to any unipotent representation of 
$G(F_q)$ one can attach an eigenvalue of Frobenius well defined up to an integer power of $q$. This result
was later used in [Digne and Michel, C.R. Acad. Sci. Paris, 1980] and by [Asai, Osaka J.Math., 1983]. In 3.34 
these eigenvalue of Frobenius are described in all cases arising in type $\ne E_8$.
On page 26 (see (d)) it is shown that $X_w$ is irreducible if and only if for any simple reflection $s$, some element
in the Frobenius orbit of $s$ appears in a reduced expression of $w$; this result has been rediscovered in 
[Bonnaf\'e and Rouquier, C.R. Acad. Sci. Paris, 2006]. (Another proof is given in [Goertz, Repres.Th., 2009]).
In 3.13 and 3.16 an explicit formula for the sum of squares of unipotent representations of $G(F_q)$ is given. This
is used in 3.24 to classify the unipotent representations in the case where $G$ is split of type $E_6$ or $E_7$ (it 
turns out that all cuspidal unipotent representations arise from the analysis in \cite{25}). In the case where $G$ is 
nonsplit $E_6$, triality $D_4$ or $F_4$, a classification of unipotent representations is again given assuming that $q$ is large; in
these cases there are cuspidal unipotent representations which do not arise from the analysis in \cite{25}.
On page 24 (see (b)) (assuming that Frobenius acts on $W$ as conjugation by the longest element $w_0$) I define for 
each $w\in W$ a bijective morphism $t_w:X_{w_0}@>>>X_{w_0}$ as follows: $t_w(B)=B'$ where $B\in X_{w_0}$ and $B'$ is 
defined by $pos(B,B')=w,pos(B',Frob(B))=w^{-1}w_0$ (so that $B'\in X_{w_0}$); I show that the $t_w$ define a 
homomorphism of the braid group in the group of permutations of $X_{w_0}$. Moreover I show that after passage to
cohomology one obtains a representation of the Hecke algebra of $W$ with parameter $-q$ on $H^*_c(X_{w_0})$. Several 
years later (around 1982) I used a similar idea in the case where $G$ is of type $D_4$ so that $W$ has simple
reflections $s_0,s_1,s_2,s_3$ with $s_1,s_2,s_3$ commuting and $w=s_1s_0s_2s_0s_3s_0\in W$. (This is unpublished but 
there is a reference to it in [Brou\'e and Malle, Ast\'erisque 1993, 5A] 
and [Brou\'e and Michel, Progr.in Math.141, 1997, page 114].) 
Let $Z(w)$ be the centralizer of $w$ in $W$. Namely, I defined three permutations $A,B,C$ of $X_w$ into itself 
(similar to $t_w$ above) such that $A,B,C$ commute with Frobenius and $ABC=BCA=CAB=$ Frobenius. The maps $A,B,C$ 
are associated to three generators $a,b,c$ of $Z(w)$ which satisfy $abc=bca=cab=w$. This suggests that the "braid group" 
corresponding to $Z(w)$ (a complex reflection group) should have the relation $ABC=BCA=CAB$. Indeed, this later appeared 
as a special case of the relations of such "braid groups" given in [Brou\'e and Malle, Ast\'erisque 1993]. The idea in this example was further 
pursued in [Digne and Michel, Nagoya Math.J.,2006] and in \cite{203}.

\head \cite{31} (with W.M.Beynon) Some numerical results on the characters of exceptional Weyl groups, 1978\endhead
For any irreducible representation $E$ of a Weyl group, the fake degree $d_E(q)$ of $E$ is defined in \cite{30}
as the polynomial in $q$ which records the multiplicities of $E$ in the various cohomology spaces of the
flag manifold. After the polynomials $d_E(q)$ were explicitly computed in \cite{29} in the case of classical
groups, it was natural to try to compute them for simple exceptional groups. This is what is done in 
the present paper, using a computer and the known character tables of $W$ (but we found and corrected 
some errors in the character table for type $E_8$). These computations were later used in \cite{34}. One 
observation of this paper is that the polynomials $d_E(q)$ are palindromic apart from a small number of 
exceptions in type $E_7$ (for $E$ of degree $512$) and $E_8$ (for $E$ of degree $4096$). My collaborator, W.M.Beynon,
was a computer expert at Warwick; I was introduced to him ny R.W.Carter.

\head \cite{33} On the reflection representation of a finite Chevalley group, 1979\endhead
The work on this paper was done in the spring of 1977; the results were presented at 
an LMS Symposium on Representations of Lie Groups in Oxford (July 1977).
I will explain the main result of this paper using concepts which were
developed several years after the paper was written (theory of character
sheaves). Let $G$ be a connected reductive group over an algebraic closure
of a finite field $F_q$ with a fixed $F_q$-split rational structure and Frobenius map 
$F:G\to G$. For each $w$ in the Weyl group one can consider (following \cite{22}) the
variety $X_w$ of Borel subgroups $B$ of $G$ such that $B,FB$ are in position $w$. Then
$G(F_q)$ acts naturally on the $l$-adic cohomology $H^i_c(X_w)$. Replacing $F$ by conjugation
by an element $g\in G$ one can consider the variety $Y_{w,g}$
of Borel subgroups $B$ of $G$ such that $B,gBg\i$ are in position $w$. The union over $g$ in $G$
of these varieties maps naturaly to $G$ and we can take the direct image $K_w$ with
compact support of the sheaf $\bar Q_l$ under this map. Then ${}^pH^iK_w$ are perverse
sheaves on $G$. Now for any irreducible representation $E$ of the Weyl group we denote by
$E_q$ the corresponding irreducible representation of $G(F_q)$ which appears in $H^0_c(X_1)$
(functions on the flag manifold of $G(F_q)$) and we denote by $E_1$ the simple perverse sheaf on G 
corresponding to $E$ which appears in $K_1$ (a perverse sheaf on $G$ up to shift, with $W$-action). The main 
result of this paper is that for any $w$ we have
$$\sum_i(-1)^i(E_q:H^i_c(X_w)=(-1)^{\dim G}\sum_i(-1)^i(E_1:{}^pH^iK_w)$$
where $(:)$ denotes multiplicity.
Here the left hand side can be interpreted as the value of the character of $E_q$ on a regular
semisimple element in a maximal torus of type $w$. This result does not compute these character 
values but it shows that these values are universal 
invariants which make sense also over complex numbers. 
Now the multiplicity $(E_1:{}^pH^iK_w)$ does not change when $E_1,{}^pH^iK_w$ are 
restricted to the variety of regular semisimple elements of $G$. But after this restriction
$E_1,{}^pH^iK_w$ become local systems and $(E_1:{}^pH^iK_w)$ is equal to the corresponding multiplicity
of local systems which can be considered independently of the theory
of perverse sheaves. It is in this form that the result above appears in the present paper
where the varieties $Y_{w,g}$ are introduced only for $g$ regular semisimple (in which case they
are shown to be smooth of dimension equal to the length of $w$). Thus this paper can be viewed
as a precursor of the theory of character sheaves which was developped in \cite{63-65,68,69}. 
As an application I determine
explicitly the value of the character of the "reflection representation" of $G_(F_q)$ 
(constructed earlier by Kilmoyer) on a regular semisimple elements of type $w$
assuming that $G$ is of type $A,D$, or $E$. Namely, it is shown that this value is equal to
the trace of $w$ on the reflection representation of $W$. At the time when the
paper was written this result was new for types $D,E$.

\head \cite{34} Unipotent representations of a finite Chevalley group of type $E_8$, 1979\endhead
This paper was written in the spring of 1978, soon after my arrival to MIT
(January 1978). This paper introduces a new type of Fourier transform (the 
"non-abelian Fourier transform"). It is a unitary involution of the vector 
space of functions on a set $M(G)$ associated to a finite group $G$; here $M(G)$ is 
the set of all pairs $(x,r)$ where $x$ is an element of $G$ (up to conjugacy) and $r$ 
is an irreducible representation of the centralizer of $x$ (up to isomorphism). 
About ten years later:

-I found \cite{77} an interpretation of the "non-abelian Fourier transform" as the 
"character table" of the equivariant complexified $K$-theory convolution algebra 
$K_G(G)$ (where $G$ acts on itself by conjugation): this (commutative) algebra has 
a natural basis indexed by $M(G)$ and its (one dimensional) representations are
also indexed by $M(G)$ hence its character table is defined;

-physicists [Dijkgraaf, Vafa, E.Verlinde, H.Verlinde, Comm. Math. Phys. 1989],
[Dijkgraaf, Pasquier, Roche, Nuclear Phys. 1990] rediscovered this Fourier 
transform (possibly with a twist by a 3-cocycle);

-Drinfeld explained it in terms of his "double" of the group algebra of G (he told me
about this around 1990). 
\nl
In this paper, the "non-abelian Fourier transform" is used to complete the 
classification and computation of degrees of the unipotent representations of 
finite Chevalley groups (started in \cite{29,30}). Note that for the analogous 
problem for classical groups, the standard (abelian) Fourier transform is 
sufficient.

It is remarkable that the "non-abelian Fourier transform" enters in an 
essential way in subsequent works in representation theory: the multiplicity 
formulas in the virtual representations $R(T,\theta)$ of \cite{22}, see \cite{57}; the 
analogous multiplicity formulas for character sheaves \cite{63-65,68,69}; the 
relation of character sheaves to irreducible characters, see \cite{71,102} and [Shoji, Adv.Math.1995].
This paper (see Section 8) also introduces the concept of "special representation"
of a Weyl group (which was further developed in \cite{36}) and that of "family" of 
unipotent representations (which contains as a particular case the notion of family 
of irreducible representations of a Weyl group). The concept of special representation
of a Weyl group was suggested by the calculations in \cite{29}.
It is nowadays used extensively in the representation theory of
reductive groups over real or complex numbers. 

\head \cite{35} (with N.Spaltenstein) Induced unipotent classes, 1979\endhead
Let $G$ be a connected reductive group over an algebraically closed field,
let $L$ be a Levi subgroup of a parabolic subgroup $P$ of $G$ and let $C$ be a unipotent
class of $G$. In this paper we associate to $L,C$ a unipotent class $C'$ of $G$ (said to
be induced by $C$); it is the unique unipotent class of $G$ whose intersection with $CU_P$ is dense in $CU_P$
(here $U_P$ is the unipotent radical of $P$). In this paper we show that $C'$ does
not depend on the choice of $P$ and that $C'\cap CU_P$ is a single $P$-conjugacy class.
When $C=\{1\}$, then $C'$ is the Richardson class defined by $L$. We give two
proofs for the independence on $P$; one of these depends on some results on representations
of a reductive group over a finite field and on the Lang-Weil estimates; the other
is more elementary but uses some case by case arguments. In this paper we also
introduce the idea of truncated induction for representations of Weyl groups
generalizing a construction of Macdonald. We show that the Springer representation
attached to an induced unipotent class is obtained from the Springer representation
of the original unipotent class by truncated induction. This has been used in
subsequent works (such as \cite{36, 48}) to compute the Springer correspondence in
certain cases arising from exceptional groups.

\head \cite{36} A class of irreducible representations of a Weyl group, 1979\endhead
In this paper (written in the summer of 1978) I give an alternative definition 
of the class $S_W$ of irreducible representations of a Weyl group $W$ of
a complex adjoint group G (introduced in \cite{34} and later called "special representations").
We have two commuting involutions $A,B$ of $Irr(W)$: $A$ is tensoring by the sign representations
and $B$ is the $q=1$ specialization of an involution of the set of irreducible representations
of the Hecke algebra given by the action of the Galois group which takes $\sqrt{q}$ to $-\sqrt{q}$.
(Note that $B$ is the identity for $W$ of classical type and it is almost the identity in general.)
Let $T=AB=BA$. In this paper it is shown that 

(i) $S_W$ is preserved by the truncated induction \cite{35} from a parabolic subgroup, and 

(ii) $S_W$ is preserved by $T$;
\nl
moreover, $S_W$ is characterized by (i)-(ii) and the fact that it contains the unit representation. 
This result has the following consequence for the two-sided cells (introduced later in \cite{37}) of $W$. Let
$w_0$ be the longest element of $W$ and let $c$ be a two-sided cell of $W$ (we assume $W\ne\{1\}$). Then $c':=cw_0=w_0c$ 
is again a two sided cell and either $c$ or $c'$ meets a proper parabolic subgroup of $W$. (This kind of
result allowed me (in the later work \cite{57}) to analyze unipotent representations inductively 
by using "truncated induction" from a proper parabolic subgroup and "duality" (which interchanges $c,c'$).)
The class $S_W$ is explicitly computed in each case (using the formalism of symbols \cite{29} for classical types), the 
results of \cite{31} on "fake degrees" and the results of \cite{35} on Springer representations. 
In this paper (Sec.9) I formulate the idea of "special unipotent class" of G (although I did not use the word 
"special"): these are unipotent classes in 1-1 correspondence with the special representations of W (under the
Springer correspondence). Since the set of special representations of $W$ admits a natural involution (given by $T$ 
above) it follows that the set of special unipotent classes admits a natural involution. (For example the class 
$\{1\}$ is interchanged with the regular unipotent class.) Later, [Spaltenstein, LNM 946, III, Springer Verlag 1982],
motivated by this paper (as mentioned in [loc.cit., p.210]) gave a definition of a subset of the set of unipotent
classes of $G$ and an order reversing involution of this subset; this definition is based on properties of the 
partial order of the set of unipotent classes and is somewhat unsatisfactory [loc.cit., p.210] for exceptional 
types. One can show that the subset defined in [loc.cit.] is the same as the set of special unipotent classes but
this was stated in [loc.cit.] not as a fact but as an analogy.
In this paper I also define a class $\bar S_W$ of irreducible representations of $W$ which contains $S_W$, but unlike 
$S_W$, depends on the underlying root system. The representations in $\bar S_W$ are obtained by truncated induction 
\cite{35} from special representations of subgroups of W which are Weyl groups of Borel-de Siebenthal subgroups
(=centralizers of semisimple elements) of the dual group $G^*$ of $G$.
I believe that the most interesting and unexpected contribution of this paper is the statement that $\bar S_W$ is 
in bijection with the set of unipotent classes in $G$ via Springer's correspondence when $G$ is of classical type and 
conjecturally in general; for exceptional groups this was verified in \cite{48}, see also 13.3 of \cite{57}. (The details 
of the proof for classical groups appeared only 25 years later in \cite{188}.) This has the following consequence: 
there is a natural map from the set of special unipotent classes of a Borel-de Siebenthal subgroup of $G^*$ (or its
dual) to the set of unipotent classes in $G$; moreover, all unipotent classes in $G$ appear in this way. This map has
been later interpreted in terms of representation theory in [Barbasch and Vogan, Primitive ideals and orbital 
integrals..., Math. Ann. 1982] (for complex groups) and in \cite{100} (for groups over $F_q$ and character sheaves).

\head \cite{37} (with D.Kazhdan) Representations of Coxeter groups and Hecke algebras, 1979\endhead
The work on this paper was done in late 1978 and early 1979. 
My motivation for this work came from the desire to construct explicitly representations
of the Hecke algebra $H$ with parameter $q$  and standard basis $\{T_w;w\in W\}$ attached to a Weyl group $W$. A basis $B$ of a vector space with
$W$-action is said to be "good" if for any simple reflection $s$ of $W$ and any $b\in B$
we have either $sb=-b$ or $sb=b+\sum_{b'\in B;b'\ne b,sb'=-b'}a_{b,b',s}b'$ where $a_{b,b',s}$ are integers.
Similarly, a basis $B$ of a vector space with $H$-action is said to be "good" 
if for any simple reflection $s$ of $W$ and any $b\in B$ we have either $T_sb=-b$ or 
$T_sb=qb+\sqrt{q}\sum_{b'\in B;b'\ne b,T_sb'=-b'}a_{b,b',s}b'$ where $a_{b,b',s}$ are integers.
In late 1977 I showed that if $u$ is a unipotent element in a semisimple group $G$ over $C$ and 
$\BB_u$ is the variety of Borel subgroups containing $u$, then in the Springer representation of $W$ on 
$H_{top}(\BB_u)$ the basis given by the irreducible components of $\BB_u$ is good with $a_{b,b',s}\ge0$.
This appeared in a letter I sent to Springer (March 1978). The same idea appeared in 
[Hotta, J. Fac. Sci. Univ. Tokyo, 1982] where the letter above is cited. In the case where $u$ is subregular 
(with $G$ of type $ADE$), $H_{top}(\BB_u)$ could be identified with the reflection representation of $W$ with 
the basis formed by simple roots. In Kilmoyer's MIT thesis (which became a part of 
[Curtis, Iwahori, Kilmoyer, Publ. Math. IHES, 1971]) an explicit $q$-deformation of the reflection 
representation of $W$ to a representation of $H$ with a good basis is found. This suggested that the bases 
of $H_{top}(\BB_u)$ (as above) may admit a q-analog which are good bases for an $H$-action. This is so for
for $G=SL_4,SL_5$. 
One of the main results of this paper is the definition of a new basis $\{c_w;w\in W\}$ of $H$ (with $q$ an 
indeterminate) which is good for the left and right action of $H$ on $H$.
I will try to explain the definition of the elements $c_w$ in a 
way somewhat different from the paper. Let $w_0$ be the longest element of $W$. 
We have  $T_sT_{w_0}=qT_{sw_0}+(q-1)T_{w_0}$. Here the right hand side 
has some coefficient  $q-1$ but if you replace $T_s$ by $T_s+1$ we obtain 
$(T_s+1)T_{w_0}=qT_{sw_0}+qT_{w_0}$ and now the right hand side has only coefficients $q$. 

Assume that $W$ is of type $A_3$ with generators $1,2,3$. We have
$$\align&T_{2132}T_{w_0}=(q^4-3q^3+4q^2-3q+1)T_{312312}\\&
+(q^4-3q^3+3q^2-q)T_{13213}+(q^4-3q^3+3q^2-q)T_{32312}\\&
+(q^4-3q^3+3q^2-q)T_{12312}+(q^4-2q^3+q^2)T_{3212}+(q^4-2q^3+q^2)T_{1232}\\&
+(q^4-2q^3+q^2)T_{2312}+(q^4-2q^3+q^2)T_{3213}+(q^4-2q^3+q^2)T_{1213}\\&
+(q^4-q^3)T_{213}+(q^4-q^3)T_{123}+(q^4-q^3)T_{321}+(q^4-q^3)T_{312}+q^4T_{13}\endalign$$
and again the right hand side has several coefficients involving powers of $q$ far from $q^4$. 
As in the case of $T_s$ we can hope that by adding to $T_{2132}$ a linear combination of the $T_y$ (with $y$ strictly less 
than $2132$) with coefficients sums of  powers of $q$ very close to $1$, the resulting sum times $T_{w_0}$ is a linear 
combination of $T_{y'}$ with coefficients sums of powers of $q$ very close to $q^4$. There is a unique way to that:
$$\align&(T_{2312}+T_{231}+T_{232}+T_{312}+T_{212}+T_{23}+T_{32}+T_{12}+T_{21}\\&+
T_{13}+T_1+T_3+(q+1)T_2+(q+1))T_{w_0}\\&=
(q^4+q^3)T_{312312}+(q^4+q^3)T_{13213}+q^4T_{32312}+q^4T_{12312}\\&
+q^4T_{3212}+q^4T_{1232}+q^4T_{2312}+q^4T_{3213}+q^4T_{1213}+q^4T_{213}\\&+
q^4T_{123}+q^4T_{321}+q^4T_{312}+q^4T_{13}.\endalign$$
We then take
$$\align&c_{2312}^*=T_{2312}+T_{231}+T_{232}+T_{312}+T_{212}+T_{23}\\&+T_{32}+T_{12}+
T_{21}+T_{13}+T_1+T_3+(q+1)T_2+(q+1).\endalign$$
This procedure works in general and leads to a basis $\{c_w^*;w\in W\}$ of $H$.
More explicitly, $c^*_w=\sum_{y\le w}P_{y,w}(q)T_y$ is characterized
by 
$$c^*_wT_{w_0}=\sum_{y\le w}q^{l(w)}P'_{y,w}(q^{-1})T_{w_0y}$$ 
where $P_{y,w},P'_{y,w}$ are polynomials in $q$ of degree at most $(l(w)-l(y)-1)/2$ if $y\ne w$ and $P_{w,w}=P'_{w,w}=1$.)
Let $(c_w)$ be the basis obtained from $c_w^*$ by the involution
$T_s\to-qT_s^{-1}$ of $H$. 

In this paper it is shown that the  basis $(c_w)$ is good for both the left and right $H$-module structure on $H$
(and in fact the scalar $a_{b,b',s}$ is independent of $s$ whenever it is nonzero. 
Note also that the definition of $c_w^*$ is applicable
to any Coxeter group by replacing the operation of multiplication by $T_{w_0}$ by the bar operation which replaces
$T_x$ by $T_{x^{-1}}^{-1}$ and $q$ by $q^{-1}$ (which is what appears in the paper). 

Also in the paper left cells, right cells and two sided cells are introduced for any Coxeter group and the
left cells in type $A$ are determined explicitly. The inversion formula 3.1 shows that the inverse of the triangular matrix 
$(P_{y,w})$ is the triangular matrix which has again the entries $P_{y,w}$ in another indexing and with some sign changes. This
inversion formula was later generalized in [Vogan, Duke Math.J. 1982] to the case of 
symmetric spaces, in which case a passage to the Langlands dual of $G$ is necessary. 
In this paper it is observed that the nontriviality of $P_{y,w}$ is very closely related to
the failure of local Poincar\'e duality on a Schubert variety. 
The fact that the equivalence relation on the set of irreducible representations of W given 
by the two-sided cells (of this paper) seemed to coincide with the equivalence relation defined by
the families (introduced earlier in \cite{34} in conection with the representation theory of finite
reductive groups), suggested that the $P_{y,w}$ may have a representation theory significance. In this paper, a 
conjecture is stated to the effect that $P_{y,w}(1)$ should be equal to the multiplicity $[L_y:M_w]$ of a simple 
highest weight module $L_y$ in a Verma module $M_w$ over the Lie algebra of $G$. Some evidence for the conjecture (in addition
to the one mentioned above) came from the fact that the matrix $[L_y:M_w]$ was known in the literature for rank $\le 3$ (Jantzen) and the 
$P_{y,w}$ could be explicitly computed in rank $\le 3$ and they matched the $[L_y:M_w]$. Another evidence came from 
[Joseph, W-module structure on the primitive spectrum...,1979] which showed among other things that the basis of the 
regular representation of $W$ given by $\{\sum_{y\in W}sign(y)sign(w)[L_y:M_w]y;w\in W\}$ is good.
A step towards proving the conjecture was made in \cite{39} where
$P_{y,w}(1)$ is interpreted as the Euler characteristic of
a certain local intersection cohomology space. The remaining statement was established in 
[Beilinson and Bernstein, C.R. Acad. Sci. Paris, 1981] and [Brylinski and Kashiwara, Invent. Math. 1981].
Later I formulated an extension of the conjecture above to express the character of ireducible
 highest weight modules with positive central charge of any affine Lie algebra which involved the
values at 1 of the entries of the matrix inverse to $(P_{y,w})$; I communicated this conjecture (and also
the similar conjecture for any Kac-Moody Lie algebra) to V.Kac and the conjecture appeared in
[Deodhar, Gabber and Kac, 1982]. (A proof was given in [Kashiwara and Tanisaki, Grothendieck Festschrift II, 1990].)
Even later \cite{88}, I formulated an extension of the conjecture above to express the 
character of ireducible highest weight modules with negative central charge of any affine Lie algebra 
which involved the values at $1$ of the $P_{y,w}$ themselves. (A proof was given in [Kashiwara and Tanisaki, Duke Math.J., 1995].)

A generalization of the notion of left/right/two-sided cells of this paper
to the case of complex reflection groups has been proposed in
[Bonnaf\'e and Rouquier, Cellules de Calogero-Moser, arxiv:1302.2720];
earlier, [Gordon and Martino, Math. Res. Lett. 16(2009)] proposed a 
generalization to complex reflection groups of the notion of family of 
irreducible representations of a Weyl group introduced in \cite{34}.

\head \cite{38} (with D.Kazhdan) A topological approach to Springer's representations, 1980\endhead
This paper was written in 1979 but the work on it was done in late summer of 1978 (except Sec.7).
In 1976, Springer defined an action of the Weyl group $W$ on the cohomology $H^*(\BB_u)$ of the variety $\BB_u$
of Borel subgroups containing a unipotent element $u$ of a reductive algebraic group over $C$, using
methods in characteristic $p>0$. Moreover in a letter to me (1977) Springer defined an action
of $W\T W$ on the cohomology $H^*_c(Z)$ of the Steinberg variety $Z$ of triples $(u,B,B')$ where $u$
is a variable unipotent element and $B,B'$ are Borel subgroups containing $u$; in the same letter
he conjectured that the representation in $H_c^{top}(Z)$ is the bi-regular representation.
In this paper an elementary construction of the Springer representation of $W$ on $H^{top}(\BB_u)$
and of the Springer representation of $W\T W$ on $H_c^{top}(Z)$ is given and the conjecture
of Springer mentioned above is proved. The construction in this paper is based on an explicit
homotopy equivalence $s_i$ from $\BB_u$ to $\BB_u$ for any simple reflection $s_i$ in $W$. We were expecting 
(but unable to prove) that the maps $s_i$ give a representation of $W$ in the group of homotopy equivalences 
modulo homotopy of $\BB_u$; we could only prove this after passage to cohomology and only in top degree. The 
stronger statement has been established later in [Rossmann, J. Funct. Analysis, 1991].
This paper's use of the Steinberg variety $Z$ of triples reappeared in \cite{72} in connection with the study
of representations of an affine Hecke algebra. 

\head \cite{39} (with D.Kazhdan) Schubert varieties and Poincar\'e duality, 1980\endhead
The work on this paper was done in early 1979.
The appendix to \cite{37} showed that the nontriviality of the polynomials $P_{y,w}$ of \cite{37} (for ordinary Weyl 
groups) was very closely related to the failure of local Poincare duality on a 
Schubert variety. It looked like the computations made in that appendix were actually computations of
local intersection cohomology in case of an isolated singularity or in the case where one meets the
singular locus for the first time. (R.Bott has suggested to Kazhdan that the results in that appendix could be
related to intersection cohomology. On the other hand I have attended a lecture of MacPherson on intersection
cohomology at Warwick in 1977 which dealt with the failure of Poincar\'e duality, as did the appendix to \cite{37}, 
and I was wondering about the connection between the two.) However, a preprint of Goresky, MacPherson gave a 
different value for the local intersection
cohomology than what we found in our case. Therefore Kazhdan and I arranged to meet MacPherson (at Brown University) to 
clarify this point. It turned out that the Goresky-MacPherson preprint had a misprint and in fact it should
have matched our computation. After this Kazhdan and I tried to identify all of $P_{y,w}$ with the local 
intersection cohomology of a Schubert variety and we succeeded in doing so (using results of Deligne). This is 
what is done in this paper. This can be viewed as a step in the proof of the conjecture on multiplicities in
Verma modules in \cite{37}. 
The idea to consider the affine Schubert variety (attached to an element in the affine Weyl group) 
as an algebraic variety also appears (perhaps for the first time) in this paper.
This paper also gives a proof of the positivity of coefficients of $P_{y,w}$ (for Weyl groups and
affine Weyl groups). In 2012 an elementary proof of the positivity valid for any Coxeter group was given by
B.Elias and G.Williamson. 

\head\cite{40} Some problems in the representation theory of finite Chevalley groups, 1980\endhead
This paper is based on a talk given in July 1979 at the Santa Cruz Conference on Finite Groups.
It states several problems. Problem I states as a conjecture the multiplicity formula for unipotent
representations in the virtual representations $R(T,1)$ of \cite{22}. This was solved in \cite{42,45,46,57},
(the last three papers make use of the results in \cite{39}). Problem II is about assigning a unipotent 
support to an irreducible representation. This was solved in large characteristic in \cite{100} and later in general 
in [Geck, Malle, Trans. Amer. Math. Soc., 2000]. Problem III relates the families \cite{34} of irreducible representations 
of the Weyl group with the two-sided cells \cite{37}; it has been solved in [Barbasch and Vogan, Math. Ann. 1982 and J.Alg.
1983]. Problem IV is a conjecture on the characters of irreducible
modular representations of a semisimple group in characteristic $p>0$ in terms of the polynomials of \cite{37}
attached to the affine Weyl group of the Langlands dual group. This was solved for $p$ larger than a fixed unknown number
by the combination of [Andersen, Jantzen and Soergel, Ast\'erisque, 1994], \cite{108,109,115,116},
[Kashiwara and Tanisaki, Duke Math.J, 1995, 1996], \cite{117}. An explicit, rather large, bound for $p$ was found in 
[Fiebig, J. reine angew. math. 2012]. The bound found by Fiebig cannot be much improved, see [Williamson, arxiv:1309.5055], as shown by [Williamson, arxiv:1502.04914] (parly in collaboration with X.He).
Problem V states that the unipotent classes of a semisimple group in large characteristic 
are in bijection with the two-sided cells (see \cite{37}) of the affine Weyl group of the Langlands dual group. 
This was solved in \cite{86}. 

\head \cite{41} Hecke algebras and Jantzen's generic decomposition patterns, 1980\endhead
In this paper I introduce and study a certain module over an affine Hecke algebra, which I now call the periodic 
module. For simplicity I define it here in type $A_1$. Let $E$ be an affine euclidean space of dimension $1$ with a 
given set $P$ of affine hyperplanes (points) which is a single orbit of some nontrivial translation of $E$. Then the 
group $\Om$ of affine transformations of E generated by the reflections with respect to the various $H$ in $P$ is 
an infinite dihedral group. The connected components of $E-P$ are called alcoves; they form a set $X$ on which $\Om$ 
acts simply transitively. Let $S$ be the set of orbits of $\Om$ on $P$. It consists of two elements. If $s\in S$ then 
$s$ defines an involution $A\to sA$ of $X$ where $sA$ is the alcove $\ne A$ such that $A$ and $sA$ contain in their closure a 
point in the orbit $s$. The maps $A\to sA$ generate a group of permutations of $X$ which is a Coxeter group $(W,S)$ (an 
affine Weyl group of type $A_1$ acting on the left on $X$). We assume that for any two alcoves $A,A'$ whose closures 
contain exactly one common point (in $P$) we have a rule which says which of the two alcove is to left (or to the 
right) of the other in a manner consistent with translations. Let $v$ be an indeterminate. Let $\HH$ be the Hecke 
algebra attached to $W,S$ and let $M$ be the free $\ZZ[v,v^{-1}]$-module with basis $X$. There is a unique $\HH$-module 
structure on $M$ such that for $s\in S,A\in X$ we have $T_sA=sA$ if $sA$ is to the right of $A$ and $T_sA=v^2sA+(v^2-1)A$ if 
$sA$ is to the left of $A$. For each $H\in P$ let $e_H\in M$ be the sum of the two alcoves in $X$ whose closures contain $H$.
Let $M^0$ be the $\HH$-submodule of $M$ generated by the elements $e_H$.
Now in the paper the higher dimensional analogue of the situation above is studied. The analogue of $X$ and the
$\HH$-modules $M,M^0$ are introduced. A bar involution of $M^0$ is introduced; it is semilinear with respect to the bar
involution \cite{37} on $\HH$. A canonical basis of the $\ZZ[v,v^{-1}]$-module $M^0$ is constructed using the bar operator on 
$M^0$ by a method similar to that of \cite{37} (but the construction is more intricate). This canonical basis is indexed
by the alcoves in X. The polynomials which give the coefficient of an alcove $B$ in the basis element corresponding
to an alcove $A$ are periodic with respect to a simultaneous translation of $A$ and $B$. They can be related to the 
polynomials attached in \cite{37} to W; this relation proves a periodicity property for these last polynomials the 
proof of which was the main motivation for this paper. (An analogous periodicity property for the multiplicities 
in the Weyl modules of a simple algebraic group in characteristic $p$ was first pointed out by [Jantzen, J. Algebra,
1977] and the periodicity result of this paper provided support for the conjecture in \cite{40} on these 
multiplicities). Shortly after writing this paper I found the folowing geometric interpretation of the results of
this paper. Let $G$ be a simply connected almost simple group over $\CC$. Let $K=\CC[[t]]$. Let $U$ be the unipotent radical 
of a Borel subgroup of $G$ and let $I$ be an Iwahori subgroup of $G(K)$. Then the set of double cosets $U(K)\bsl G(K)/I$ is 
(noncanonically) the affine Weyl group and (canonically) the set $X$ of alcoves as above. (A closely related 
statement is contained in [Bruhat and Tits, Groupes r\'eductifs sur un corps local, Publ. IHES, 1972, Prop.(4.4.3)(1).)
This led me to the statement that the periodic polynomials of this paper can be interpreted as local intersection
cohomologies of the (semiinfinite) $U(K)$-orbits on $G(K)/I$. This statement appears without proof in \cite{59}; a proof 
appears in [Finkelberg and Mirkovic, Semiinfinite flags I; Feigin, Finkelberg, Kuznetsov and Mirkovic, Semiinfinite flags 
II, Transl. of Amer. Math. Soc., 1999]. 

\head \cite{42} On the unipotent characters of the exceptional groups over finite fields, 1980\endhead
In this paper I determine the multiplicities of the unipotent representations
of an exceptional group over a finite field $F_q$ in the virtual representations
$R(w,1)$ of \cite{22}, assuming that $q$ is large. These multiplicity formulas were
conjectured in \cite{40}. The proof given in the paper uses the formulas (known at
the time) for the dimensions of unipotent representations and unlike the
later proof \cite{54} (where the restriction on $q$ was removed) it does not use
intersection cohomology methods. The method of this paper does not seem to
be strong enough in the case of classical groups which was treated later in \cite{45,46} 
using intersection cohomology methods. 

\head \cite{43} On a theorem of Benson and Curtis, 1981\endhead
Let $H$ be the Hecke algebra over $\QQ$ associated to the Weyl group $W$ and to the
parameter $q$, a power of a prime number. In 1964, Iwahori conjectured that $H$ is isomorphic
to the group algebra $\QQ[W]$. Tits showed (in an exercise in Bourbaki) that this
conjecture holds if $\QQ$ is replaced by its algebraic closure. In 1972 Benson
and Curtis showed that Iwahori's conjecture was true as originally stated but
Springer found a gap in the proof (for type $E_7$). (The Benson-Curtis proof 
was correct for types other than $E_7,E_8$.)  Springer in fact showed that the character of
a $512$-dimensional irreducible representation of $H$ (of type $E_7$) 
definitely involves a square root of q. 
In this paper (written in 1980) I construct an algebra isomorphism of
$\QQ[\sqrt q]\ot H$ with $\QQ[\sqrt q][W]$.
The key new observation is as follows. Consider the vector space
spanned by the elements of a fixed two-sided cell of $W$. There is a left action on this
vector space for the Hecke algebra $H$ with parameter $q$ in which the basis elements are identified with
the elements of the new basis \cite{37} of $H$; there is also a right action on this vector space for
the Hecke algebra $H'$ with another parameter $q'$ in which the basis elements are identified with
the elements of the new basis \cite{37} of $H'$. The two actions obviously commute with each other if $q=q'$
but surprisingly they also commute with each other when $q,q'$ are independent.
The proof is based on some properties of primitive ideals in an enveloping algebra. 
The isomorphism I construct is explicit unlike those in 
earlier approaches. The use of the theory of primitive ideals can nowadays
be eliminated and replaced by the use of the "a-function" introduced in \cite{60}.
This paper also gives $W$-graphs (in the sense of \cite{37}) for the left cell representations of $H$ in the 
noncrystallographic case $H_3$ and an example analogous to the $512$-dimensional representation (for $E_7$) is pointed 
out in type $H_3$.

\head\cite{44} Green polynumials and singularities of unipotent classes, 1981\endhead
In this paper I find a relation between:

(1) the local intersection homology groups of the closure of a unipotent class in $GL_n$;

(2) the local intersection homology of an affine Schubert variety in an affine grassmannian of type $A$;

(3) the character value at a unipotent element of an irreducible unipotent representation of $GL_n(F_q)$.
\nl
The connection (1)-(3) is a precursor of the theory of character sheaves which was developped in \cite{63-65,68,69}. The connection 
(2)-(3) implies that the groups in (2) can be described in terms of multiplicities of weights for 
the finite dimensional representations of $GL_n(\CC)$ which was an inspiration for the paper \cite{53} (a generalization from $GL_n$ to 
a general reductive $G$). This paper also formulated the idea (new at the time) that 
the Springer resolution is a small map and uses this idea to give a new 
definition of the Springer representations of a Weyl group in terms of
intersection cohomology (unlike previous definitions this was valid in 
arbitrary characteristic). This shows in particular that the direct image 
of the constant sheaf under the Springer resolution is a perverse sheaf 
up to shift.
Conjecture 2 of this paper was subsequently proved by [Borho and MacPherson, C.R. Acad. Sci. Paris 1981].
The method introduced in this paper to construct Springer's representations
has been used in later papers:

(a) to construct the "generalized Springer correspondence" \cite{59};

(b) to construct analogues of the Springer representation over parameter spaces which yield representations of graded affine Hecke algebras \cite{81};

(c) to construct a version of Springer representations for affine Weyl groups \cite{125}; 

(d) to construct a Weyl group action on the cohomology of certain quiver varieties \cite{149}.

\head \cite{45,46} Unipotent characters of symplectic and odd orthogonal groups over a finite field, 1981;
Unipotent characters of the even orthogonal groups over a finite field, 1982\endhead
The first of these papers was conceived during a visit to the Australian National University, Canberra (January, 1981);
the second one was written later in 1981. Let $G(F_q)$ be a group as in the title. In \cite{30, Conj.4.3} I conjectured the precise pattern which
gives the multiplicities of the various unipotent representations in the virtual representations
$R(w,1)$ of \cite{22} or equivalently in the linear combinations $R_E$ of the $R(w,1)$ with coefficients given by an
irreducible character $E$ of the Weyl group; namely the pattern should be the same as the pattern \cite{29}
describing the dimensions of unipotent representations as linear combinations of fake degrees.
This conjecture is established in this paper. The main new technique in the proof 
is the use of the local intersection cohomology of the closures of the varieties $X_w$ of \cite{22}  which I
show that is the same as the local intersection cohomology of a Schubert variety and hence \cite{39}
is computable in terms of Hecke algebras. Another new technique used in the paper is the systematic
use of the leading coefficients of character values of the Hecke algebra. These techniques were
later generalized to any reductive group (see \cite{57}).
 
\head \cite{47} (with P.Deligne) Duality for representations of a reductive group over a finite field, 1982\endhead
In 1977 I found a definition of an operation $D$ in the complex representation 
group of a reductive group over $F_q$ which to any representation $E$ associates
$\sum_P(sgn_P)\ind_P aa^*res_P(E)$ where $P$ runs over the parabolic subgroups over 
$F_q$ containing a fixed Borel subgroup over $F_q$, $\ind_P$ is induction from 
$P(F_q)$ to $G(F_q)$, $a$ is lifting from $P(F_q)/U_P(F_q)$ to $P(F_q)$,
$a^*$ is the adjoint of $a$ and $\res_P$ is restriction to $P(F_q)$; $sgn_P$ is a sign. It was known at the time 
(Curtis) that $D$ takes the unit representation to the Steinberg representation. If $E$ is cuspidal
then $a^*res_P(E)$ is zero if $P\ne G$ hence $DE=\pm E$. But my main motivating example
was one which I encountered in \cite{13} where $G=GL_n(F_q)$, the complex numbers are
replaced by $F_q$ and $E$ is the natural representation of $G$ on $F_q^n$.
In that case $DE$ can be defined as above and can be viewed as a reduction
$\mod p$ of a cuspidal complex representation of $G$ of dimension $(q-1)(q^2-1)...(q^{n-1}-1)$; 
this was the main observation on which the work \cite{13} was based. I conjectured that over complex numbers, $D$
takes any irreducible $E$ to an irreducible representation (up to sign) and that $D^2=1$. 
In 1977 I communicated this conjecture to D.Alvis and C.W.Curtis (at the Corvallis Conference) and (separately) 
to N.Kawanaka. My conjecture was proved (at the level of characters) by [Alvis, Bull. AMS, 1979],
[Curtis, J. Algebra, 1980] 
and independently by [Kawanaka, Invent. Math., 1982]. In the present paper a version of $D$ at
the level of representations (rather than characters) is given. As an application another proof of the conjecture
is given. The operation $D$ played a key role in my later work \cite{57} where it was used to analyze unipotent 
representations inductively (in conjunction with "truncated induction" from a proper parabolic subgroup). An 
analogous operation plays a key role in the classification of character sheaves. In 
[A.M.Aubert, Trans.Amer. Math. Soc., 1995] a study of a p-adic analogue of the operation D defined in this paper is 
made.

\head \cite{48} (with D.Alvis) On Springer's representations for simple groups of  type $E_n$ $(n=6,7,8)$, 1982\endhead
Let $G$ be as in the title (over $\CC$). In this paper we compute the Springer representation of the Weyl group $W$ of $G$
corresponding to any unipotent class and the local system $\CC$ on it. There are three tools that are used in the 
proof: (a) the compatibility of truncated induction with the Springer correspondence \cite{35}; (b) the conjecture (2)
in \cite{44} which was just proved by Borho and MacPherson; (c) an 
induction formula for the total Springer representation
for a unipotent element contained in a proper Levi subgroup. Moreover, using the induce/restrict tables of Alvis 
we showed that the class of irreducible representations of $W$ thus obtained coincides with the class $\bar S_W$ 
introduced in \cite{36}, thereby completing the proof of the conjecture at the end of \cite{36} (which at the time of \cite{36}
was already known for classical types and $G_2$). In the appendix (by Spaltenstein) the rest of the Springer 
correspondence (involving irreducible local systems $\ne\CC$) is determined.

\head \cite{49} (with D. Alvis) The representations and generic degrees of the
Hecke algebra of type $H_4$, 1982\endhead
In this paper the irreducible representatios of a Hecke
algebra of type $H_4$ are explicitly constructed in terms of $W$-graphs.
Moreover, the generic degrees of these representations are explicitly 
computed. Remarkably, these turn out to be polynomials in $q$ rather than
rational functions. This fact suggested  to me that a theory of unipotent representations for $H_4$
should exist, and led to my paper \cite{110}.

\head\cite{50} A class of irreducible representations of a Weyl group II, 1982\endhead
This paper was written in early 1981. Let $W$ be a Weyl group and let $\Irr W$ be
the set of irreducible representations of $W$ (up to isomorphism).
In \cite{34} a partition of $\Irr W$ into subsets called families was described.
The definition was such that the degrees of unipotent representations
of a finite Chevalley groups were linear combinations of fake degrees of
objects of $\Irr W$ in a fixed family. In the present
paper an elementary definition of families is given. More precisely a 
collection of possibly reducible representations (called cells and in later
papers, constructible representations) is defined by induction. Namely it is required that by applying a certain kind of truncated induction to a cell of
a proper parabolic subgroup one obtains a cell of $W$; moreover by tensoring
a cell by the sign representation of W one obtains again a cell. The cells are
obtained by applying a succession of such operation starting with the trivial
one dimensional representation of $W$. In this paper the cells of any W are
explicitly determined. It is shown that any $E$ in $\Irr W$ appears in some cell;
every cell contains a unique special representation (in the sense of \cite{36})
which in fact has multiplicity one; and two cells have a common irreducible component if and only if they contain the same special representation. Therefore
we can define an equivalence relation on $\Irr W$ as follows: $E,E'$ in $\Irr W$ are
equivalent if there exist cells $c,c'$ such that $E$ appears in $c$, $E'$ appears in $c'$
and $c,c'$ have the same special component. The equivalence classes are called
families. In the paper it is conjectured that the cells of $W$ are exactly the
representations of $W$ that are carried by the left cells of $W$ (in the sense of \cite{37}). This conjecture was proved in \cite{70}.
Using the results of this paper one can give a new definition of the involution
of the set of special representations of $W$ (see the comments to \cite{36}) which 
bypasses the consideration of a Galois group action: namely the involution 
maps a special representation $E$ in a family $f$ to the unique special 
representation in the family ($f$ tensored by sign).

\head\cite{51} (with D.Vogan) Singularities of closures of $K$-orbits on a flag manifold, 1983\endhead
The work on this paper was done in late 1980. 
Its main object of study was the local intersection cohomology (l.i.c.) of the closure of a
$K$-orbit on the flag manifold of $G$ where $K$ is the identity component of the fixed point set of an involution of a complex reductive
group $G$. At the time it was known from the work of Beilinson and Bernstein that
this l.i.c. is closely related to the computation of multiplicities in standard
module of the various irreducible representations of a real reductive group
attached to the involution in the same way as the l.i.c. of Schubert varieties
was known to be closely related to multiplicities in Verma modules. 
The problem of determining the l.i.c. in the present case was a generalization of
the problem of determining the l.i.c. of Schubert varieties solved in \cite{39}. But
the method of \cite{39} did not work in the present case, partly due to the presence
of non-trivial equivariant local systems (of order two) on the $K$-orbits. Unlike
in \cite{39} in this paper the connection with the representation theory of real groups
is used in the computation; also the purity theorem of Gabber (which was not available
at the time of \cite{39}) plays a key role in the proof. The main result of this paper
is that the l.i.c. are described in terms of some new polynomials $P_{\gamma,\delta}$, where
each of $\gamma$  and $\delta$ is a $K$-orbit together with a $K$-equivariant irreducible local 
system on it, which are explicitly computable and which generalize
the polynomials $P_{y,w}$ of \cite{37}. (Later work by Fokko Du Cloux has made possible the
computation of $P_{\gamma,\delta}$ on a computer.) This paper also contains an interpretation 
of the product in the Hecke algebra and in certain modules over it in terms of convolution 
in derived categories (involving operations of inverse image, direct image and tensor 
product in derived categories). This interpretation which has become part of the folklore
has been also found around the same time by MacPherson, see [Springer, Sem. Bourbaki 589, 1982].
A proof of the results of this paper which is purely geometric (that is it does not rely
on representation theory of real groups) has been later found by [Mars and Springer, Represent.
Th., 1998].

\head\cite{52} (with P.Deligne) Duality for representations of a reductive group over a finite field, II, 1983\endhead
Let $G$ be a connected reductive group over $F_q$. In this paper it is shown that 
the "duality operator" $D$ of \cite{47} applied to the virtual representation 
$R(T,\theta)$ in \cite{22} is equal (up to sign) to $R(T,\theta)$. The proof is based 
on an inner product formula between $R(T,\theta)$ and an $R(L,r)$ (as in \cite{24})
where $L$ is a Levi subgroup over $F_q$ of a parabolic (not necessarily over $F_q$)
and $r$ is a representation of $L(F_q)$. The proof of this orthogonality
formula given in the paper contains a (not very serious) error. The
corrected proof (which I supplied to Digne and Michel at their request) 
appears in the book [Digne, Michel, Representations of finite groups of Lie 
type, 1991, 11.13].

\head\cite{53} Singularities, character formulas and a $q$-analog of weight multiplicities, 1983\endhead
This paper was written in 1981 and presented at the Luminy Conference on Analysis
and Topology on Singular Spaces (July 1981).
In this paper I find a very close connection between

-the category $A$ of finite dimensional representations of a complex simply connected group $G$ and

-the category $A'$ of $G^*[[\e]]$-equivariant perverse sheaves on the affine Grassmannian associated to the
Langlands dual $G^*$ of $G$.
\nl
In more detail, let $\L^+$ be the set of dominant weights of $G$.
For $x\in\L^+$, let $V_x$ be the finite dimensional irreducible representation of $G$ 
corresponding to $x$ and let $m_y(V_x)$ be the multiplicity of $y\in\L^+$ in $V_x$. Let $M_x$ be the element 
of the (extended) affine Weyl group $W$ of $G^*$ which has maximal length in the
double coset of $x$ with respect to the usual Weyl group $W_0$. 
Let $H$ be the affine Hecke algebra of $W$ and let $(C_w)_{w\in W}$ be the basis \cite{37} of $H$.
For $x\in\L^+$ we set $\g_x=\pi^{-1}C_{M_x}$ where $\pi=q^{-\nu/2}\sum_{w\in W_0}q^{l(w)}$;
here $\nu$ is the number of positive roots and $l(w)$ is the length of $w$.
For $x\in\L^+$ let $\bar O_x$ be the closure of the $G^*[[\e]]$-orbit coresponding to $x$ in the affine
Grassmannian and let $\Pi_x$ be the corresponding simple object of $A'$.
For $w,w'$ in $W$ let $P_{w',w}$ be the polynomial defined in \cite{37}.
Here are the main results of this paper.

(I) For $x,y$ in $\L^+$ we have $m_y(V_x)=P_{M_y,M_x}(1)$. (Thus the weight multiplicities $m_y(V_x)$ are related 
to the dimension of stalks of $\Pi_x$.)

(II) For $x,y$ in $\L^+$ we have $\g_x\g_y=\sum_{z\in\L^+}c_{x,y,z}\g_z$ where $c_{x,y,z}$ are
natural numbers (apriori they are only polynomials in $q$). An equivalent statement is that the convolution 
$\Pi_x*\Pi_y$ is a direct sum of objects $\Pi_z$ $(z\in\L^+)$ without shifts; or that the map which defines this 
convolution is semismall.

(III) For $x,y,z$ in $\L^+$, the number $c_{x,y,z}$ in (II) is equal to
the multiplicity of $V_z$ in the tensor product $V_x\otimes V_y$.

(IV) For $x$ in $\L^+$, the vector space $V_x$ is isomorphic to the total intersection cohomology
of $\bar O_x$.
\nl
(Statement (IV) appears in the last line of this paper; note that the odd intersection
cohomology of $\bar O_x$ is zero.)

Statement (II) is called the "miraculous theorem" in [Beilinson and Drinfeld,  Quantization of
Hitchin integrable system... (1991), 5.3.6].
It is equivalent to the fact that A' is a monoidal
category under convolution. Statement (III) suggests that this monoidal category is equivalent to $A$ with its
obvious monoidal structure and statement (IV) suggests the definition of a fibre functor
for $A'$ which would enter in the construction of such an equivalence.
The tensor equivalence of $A$ and $A'$ was established in [Ginzburg, arxiv:alg.geom./9511007] 
based on the results of
this paper (using (II) and the fibre functor above), except that the commutativity isomorphism for $A'$ 
given there was incorrect and was later provided by Drinfeld (whose construction is sketched in 
[Mirkovic and Vilonen, Math. Res. Lett. 2000]). 
Thus the equivalence of $A,A'$ as monoidal categories (now known as the "geometric Satake equivalence")
has been established by combining the ideas of this paper with those of Ginzburg and Drinfeld. 
A version of the geometric Satake equivalence in positive characteristic is established in 
[Mirkovic and Vilonen, Math. Res. Lett. 2000].

Now by (I) each weight multiplicity appears by setting $q=1$ in a polynomial in 
$q$ with positive coefficients; hence that polynomial can be viewed as 
a "$q$-analog of weight multiplicities", hence the title of the paper. 
Subsequently, a (partly conjectural) interpretation of these $q$-analogs was given purely in 
terms of representations of $G$ in [R.K.Gupta (later Brylinsky), Jour. Amer. Math. Soc. 1989]; 
this was later confirmed in [Joseph, Letzter and Zelikson, Jour. Amer. Math. Soc. 2000]. 
In this paper I also introduce a $q$-analogue of the Kostant partition function and prove that it is equal to the $q$-analogue of 
weight multiplicities in the stable range. The proof of (II) given in this paper relies on some dimension 
estimates in \cite{41} involving semiinfinite geometry in disguise. Another result of this paper is a description of the
affine Grassmannian as an ind-variety (as a subset of the set of selfdual orders in
a simple Lie algebra over $\CC((\e))$ given by explicit equations).

We now change notation and assume that $G$ is a simply connected algebraic group 
over an algebraically closed field of characteristic $p>0$. Now $\L^+$ still makes sense
and for $x\in\L^+$, we denote by $L_x$ the corresponding simple $G$-module. 

One of my motivations to write this paper was to produce evidence for my conjecture 
on modular representations of $G$ (Problem IV in \cite{40}). More precisely, before writing this paper, I 
understood that the conjecture in \cite{40} implies statement (I) above. Thus, a proof of (I) would be
evidence for the validity of the conjecture in \cite{40} and this was a motivation for me to try to prove (I). 
(At that time I already knew that (I) is true in type $A$, as a consequence of \cite{44}.)
I will now sketch how (I) can be proved assuming that the conjecture in \cite{40} holds. From that conjecture
we have 
$$ch(L_{px})=\sum_{y\in\L^+;y\le x}P_{M_y,M_x}(1)f_y^{-1}(\sum_{w\in W_0}sign(w)ch(V_{py+w(\rho)-\rho}))$$
where $f_y$ is the order of the stabilizer of $y$ in $W_0$ and for any weight
$\nu$ we set
$$ch(V_\nu)=\sum_{w'\in W_0}sign(w')e^{w'(\nu+\rho)}/\sum_{w'\in W_0}sign(w')e^{w'\rho}.$$
(We assume that $p$ is large compared to $x$.) We have 
$$\align&\sum_{w\in W_0}sign(w)ch(V_{py+w(\rho)-\rho})\\&
=\sum_{w',w in W_0}sign(ww')e^{pw'y+w'w(\rho)}/\sum_{w'\in W_0}sign(w')e^{w'(\rho)}=
\sum_{w'\in W_0}e^{pw'y}. \endalign$$
Thus,
$$ch(L_{px})=\sum_{y\in\L^+;y\le x}P_{M_y,M_x}(1)f_y^{-1}\sum_{w'\in W_0}e^{pw'y}.$$
By the Steinberg's tensor product theorem we have 
$$ch(L_{px})=\sum_{y\in\L^+;y\le x}m_y(V_x)f_y^{-1}\sum_{w'\in W_0}e^{pw'y}.$$
We deduce that $P_{M_y,M_x}(1)=m_y(V_x)$, as stated in (I).

\head\cite{54} Some examples of square integrable representations of semisimple $p$-adic groups, 1983\endhead
Let $G$ be the group of rational points of a simple split adjoint algebraic group
over a nonarchimedean local field whose residue field has $q$ elements.
This paper introduces the notion of unipotent representation of $G$; these are the
irreducible admissible representations of $G$ whose restriction to some parahoric subgroup of $G$
contain a unipotent cuspidal representation of the "reductive" quotient of $G$.
Let $U'$ be the set of unipotent representations of $G$ and let $U$ be the subset of $U'$ formed by the 
Iwahori-spherical representations of $G$.  

According to the Deligne-Langlands conjecture, $U$ is in finite to one correspondence with the set of pairs 
$(s,u)$ where $s,u$ are a semisimple element and a unipotent 
element (up to conjugacy) in the complex "dual" group such that $su=u^qs$.

One of the main contributions of this paper is the formulation
of a refinement for the Deligne-Langlands conjecture in which 
a third parameter is added to the Deligne-Langlands parameters, namely
an irreducible 
representation $\rho$ of the group $A(s,u)$ of connected components of the simultaneous centralizer of $s,u$
on which the centre of the dual group acts trivially.

More precisely, in this paper I state the conjecture that the triples $(s,u,\rho)$ as above are in canonical bijection
with $U'$ and that $U$ is in bijection with the set of triples $(s,u,\rho)$ such that $\rho$ appears in the 
cohomology of the variety $X$ of Borel subgroups containing $s$ and $u$. 

The idea of this paper, to enrich a Langlands parameter by adding to it an irreducible representation
of a certain finite group, has been also stated several years later (in more generality) in
[Vogan, The local Langlands conjecture, Contemp. Math. 145, 1993].

A good thing about the refined conjecture (for $U$) is that it indicates that the representations 
of the affine Hecke algebra may be constructed geometrically in terms of a space like $X$. The connection
with geometry became even stronger after the equivariant $K$-theoretic approach 
of \cite{66} was found and led to the solution of the (refined) conjecture for $U$
in \cite{67,72}. The (refined) conjecture for $U'$ was established in \cite{123}.
I have arrived at the idea of the refined conjecture by experiments performed 
in this paper: I constructed explicitly (using $W$-graphs) the reflection 
representation and some closely related representations of the affine Hecke
(and I showed that they are often square integrable by some very complicated
computation); these representations correspond conjecturally to the subregular
unipotent element and this provided evidence for the refined conjecture.
In these examples I also found that the weight structure of the representations
I construct can be interpreted in terms of the geometry of the varieties X above,
further reinforcing the idea that the geometry of $X$ should play a role in the
proof of the conjecture. 
In this paper I introduce a description of the affine Weyl group of type $A$
as a group of periodic permutations of the integers. This point of view was later
used extensively in [Shi, The Kahdan-Lusztig cells in certain affine..., Springer LNM 1179, 1986].
I also give a conjecture giving the number of left cells in each two sided cell of an 
affine Weyl group of type $A$, which was later proved by [Shi, loc.cit.] and a conjecture 
describing explicitly the two sided cells of an affine Weyl group of type $A$, which I
later proved in \cite{62}.
Another result of this paper is the construction of an imbedding of a Coxeter group
of type $H_4$ (resp. $H_3$) into the Weyl group $E_8$ (resp. $D_6$) which has the property of 
sending any simple reflection to the product of two commuting simple reflections and any
 element of length $n$ to an element of length $2n$ (this is part of a
general result about imbedding of Coxeter groups, see 3.3.)
This imbedding has been rediscovered ten years later in [Moody and Patera, J.of Physics,A, 1993].

The argument in 2.8 has been used in the later papers \cite{67,72}
to prove square integrability of certain geometrically defined representations
of an affine Hecke algebra.
The last sentence in 2.11 was later proved in \cite{78}.

\head \cite{56} Open problems in algebraic groups, 1983\endhead
In the summer of 1983 I participated in a Taniguchi conference in Katata, Japan. The participants
were asked to write up a list of open problems. Here are some of the problems on my list.

(1) Let $W$ be an affine Weyl group. Then the number of left cells contained in the two-sided cell corresponding under the bijection in [86] to the conjugacy
class of a  unipotent element $u$ in a reductive group over $\CC$ of dual type to that of $W$ is equal to the dimension of the
part ot the cohomology of the Springer fibre at $u$ invariant under the action of the centralizer of $u$.

(2) Let $W$ be as in (1). We identify $W$ with the set of (closed) alcoves in an
euclidean space in the standard way. Let $A,B$ be two alcoves in the same two
sided cell. Show that $A,B$ are in the same left cell if and only if there exists
a sequence of alcoves $A=A_0,A_1,...,A_n=B$ (all in the same two sided cell)
such that $A_i,A_{i+1}$ share a codimension $1$ face for $i=0,1,...,n-1$. Show that the union of alcoves in a left cell
is a contractible polyhedron. Show that similar results hold for a finite Weyl
group by replacing the euclidean space with the corresponding 
triangulated sphere. 

\head\cite{57} Characters of reductive groups over a finite field, 1984\endhead
Let $G$ be a connected reductive group with connected centre defined over $F_q$.
The main contribution of this book (written in 1982) is the classification of the irreducible representations of $G(F_q)$
and the computation of their multiplicities in the virtual representations $R(w,\theta)$ of \cite{22}.
(Earlier, this kind of results were known for unipotent representations with $q$ large, see \cite{42,45,46}; the classification (but 
not the multiplicities) for classical groups with any $q$ was also known \cite{29}).

Let $L$ be a $G$-equivariant line bundle over the flag manifold of $G$ and let $L-0$ 
be the complement of the zero section of $L$. Now Ch.1 contains 

($*$) the computation of the local intersection cohomology of $L-0$ with 
coefficients in certain "monodromic" local systems on some smooth subvarieties
 of $L-0$, in terms of the polynomials \cite{37} for the Weyl group of the 
centralizer of a semisimple element $s$ in the dual group. 

(This generalizes results of \cite{39} which correspond to the case $s=1$.)
Since these local intersection cohomology groups were at the time known (by 
Beilinson and Bernstein) to compute multiplicities in Verma modules with 
regular rational highest weight, $(*)$ was a new instance of a connection 
between representations of a group and geometry of the dual group. Another 
proof of the multiplicity formulas
in Verma modules with regular rational highest weight was later found in
[Soergel, Jour. Amer. Math. Soc. 1990, Theorem 11] where these multiplicities
are directly related to analogous multiplicities for integral highest weight,
thus bypassing $(*)$. Note that $(*)$ is used in [Beilinson and Bernstein, 
A proof of Jantzen's conjecture, Adv. Sov. Math. 1993]. An affine generalization 
of $(*)$ is given in \cite{117} where it is used as one of the steps in 
the proof of the character formula for quantum groups of nonsimplylaced type
at a root of $1$. Finally, $(*)$ is used in Ch.2 of this book to 
determine the local intersection cohomology 
of the closures of the varieties $X_w$ of \cite{22} with coefficients in local 
systems associated with the covering $\tilde X_w$ of $X_w$ described in \cite{22}. 
Again the result is expressed in terms of the polynomials of \cite{37} for the Weyl
group of the centralizer of a semisimple element in the dual group.

\head \cite{58} Characters of reductive groups over finite fields, 1984\endhead
This is based on my talk at the ICM-1982 held in Warsaw in 1983 (the 1982 event 
was postponed due to the martial law). This paper is an exposition of the main 
results of \cite{57} (written in 1982) which were under the assumption of connected 
centre. But in the present paper that assumption was removed. In order to remove 
two words: "connected centre" from my paper I had to do two months of intensive work 
(June/July 1983) mainly with the case of $Spin_{4n}$. These computations with spin 
groups (not included in the paper where no proofs were given) have been published 
25 years later in \cite{180} with some earlier hints given in \cite{83}.

\head \cite{59} Intersection cohomology complexes on a reductive group, 1984\endhead
This paper was written in late 1982 and early 1983. Let $G$ be a connected reductive 
group over an algebraically closed field of characteristic $p\ge0$.
Let $X$ be the (finite) set of all pairs $(C,E)$ where $C$ is a unipotent class in $G$ and
$E$ is a $G$-equivariant irreducible local system on $C$ (up to isomorphism).
In the late 1970's Springer showed that (if $p=0$ or if $p$ is large)
there is a natural bijection between a certain subset $X_0$ of $X$ and the
set $\Irr(W)$ of irreducible representations of the Weyl group $W$ of $G$.
In \cite{44} I gave a new definition of the Springer representations of $W$ using intersection
cohomology methods which is valid without restriction on $p$, but the proof that it
induces a bijection between $X_0$ (which can be defined for any $p$) and $\Irr(W)$ was first given for arbitrary $p$ 
in this paper, using a study of sheaves on the variety of semisimple classes.
In this paper I show (extending the method of \cite{44}) that a suitable enlargement
of $\Irr(W)$ is in canonical bijection ("generalized Springer correspondence") with $X$ itself.
The enlargement is a disjoint union of sets of the form $\Irr(W_i)$ where $W_i$ is a collection
of Weyl groups (one of which is $W$). Of particular
interest are the objects of $X$ for which the corresponding $W_i$ is $\{1\}$.
These are the "cuspidal local systems" (c.l.s.) which are introduced,
studied and classified in this paper.
A G-equivariant local system $E$ on a unipotent class $C$ of $G$ 
is a c.l.s. if for any proper parabolic $P$ of $G$ with unipotent radical $U_P$
and any unipotent $g$ in $P$, the $d$-th cohomology with compact support of 
$C\cap gU_P$ with coefficients in $E$ is zero (where $d$ is $\dim(C)$ minus the 
dimension of the conjugacy class of $g$ in $P/U_P$); note that if $d$ is replaced by $d'>d$ 
then the corresponding vanishing property holds for any $E$. 
A new feature of this paper is the explicit combinatorial description of the
generalized Springer correspondence in terms of some objects closely related
to the "symbols" in \cite{29}. This was new even for the ordinary 
Springer correspondence which was previously known only in the form of an
algorithm (Shoji), rather than by a closed formula. In the case where $G$ is a spin group with 
$p$ odd, the generalized Springer correspondence 
gives a combinatorial interpretation of the Jacobi triple product formula (see Section 14).
Another new result of this paper was a definition of "admissible complexes" on $G$,
a class of perverse sheaves on $G$ whose existence was conjectured in
\cite{57, 13.7,13.8} where the required class of perverse sheaves was defined for $G=GL_n$.
One of the main ingredients in the definition of admissible complexes is the notion of
c.l.s. (see above) extended from unipotent classes to "isolated classes". The 
admissible complexes on $G$ reemerged in another incarnation (as "character sheaves")
in the series \cite{63-65,68,69}. 

\head\cite{60} Cells in affine Weyl groups, 1985\endhead
Let $W$ be a Weyl group or an affine Weyl group. This paper develops some 
techniques for computing the left/two-sided cells \cite{37} of $W$. The main 
contribution of this paper is the definition of the function $a:W@>>>\NN$. 
For $w\in W$, I define $a(w)$ essentially as the order of the 
worst pole of the coefficient of $C_w$ (the Hecke algebra element of \cite{37}) in a 
product $T_xT_y$ of two (variable) standard basis elements of the Hecke algebra.
I show that $a$ is constant on the two-sided cells of $W$.
When $W$ is of affine type the fact that $a(w)$ is well defined needs a proof
(given in the paper); in fact I show that $a(w)$ is at most the number $N$ of positive roots. Therefore the set $W_*=\{w\in W;a(w)=N\}$ 
is of particular significance.
Let $W_!$ be the set of all products $abc$ where $a,b,c\in W$, the length of $abc$ is
the sum of the lengths of $a,b,c$ and $b$ has length $N$ and is contained in a
finite parabolic subgroup of $W$. In the paper it is shown that $W_!\sub W_*$; in particular, 
$W_*$ contains "almost all" elements of $W$. 
In this paper, using the function $a$, I describe explicitly the decomposition of the affine Weyl group $W$ of type 
$A_2,B_2,G_2$ into left/two-sided cells in terms of a picture in which $W$ is
viewed as the set of alcoves in a decomposition of an euclidean plane and each
alcove is colored according to the two-sided cell to which it belongs.
It turns out that, for affine $A_2,B_2,G_2$ the number of two-sided 
cells is $3,4,5$; this was one of the pieces of evidence which led to my 
conjecture \cite{40} (restated in this paper) on the relation between two-sided 
cells and unipotent classes.
From the results of this paper one can see that in rank $2$ one has $W_*=W_!$
and $W_*$ is a single two-sided cell. This was extended to arbitrary rank
in [Shi, J. Lond. Math. Soc. 1987] and [B\'edard, Commun. in Alg. 1988].

\head\cite{62} The two-sided cells of the affine Weyl group of type $A$, 1985\endhead
The results of this paper were presented at a conference at MSRI in May 1984.
In early 1983 I have learned from R.Carter about the remarkable  work of 
his Ph.D. student
J.Y.Shi (at Warwick) in which Shi determined explicitly the left cells of the
affine Weyl group $W$ of type $A_n$; it turned out that Shi's methods were not
sufficient to determine the two-sided cells of $W$ (for which I formulated a
conjecture in \cite{54}). After I introduced the function $a$ on $W$ in \cite{60}, I realized
that the results of Shi together with the use the function $a$
are sufficient to determine the two-sided cells of $W$. This is what is done
in this paper; see also [Shi, The Kazhdan-Lusztig cells in certain affine...,
Springer LNM 1179, 1986]. 

\head \cite{66} Equivariant $K$-theory and representations of Hecke algebras, 1985\endhead
The work on this paper was done at the Tata Institute, Bombay, in December 1983.
At the time when this paper was written, the parameter $q$ of a Hecke algebra 
was viewed as a number, an indeterminate, a Tate twist or a shift in a derived 
category. One of the main contributions of this paper is to formulate the idea (new at the time) to view $q$  
as the generator of the equivariant $K$-theory of a point with 
respect to the circle group and that various modules of the affine Hecke algebra $H$ 
can be realized in terms of equivariant $K$-theory with respect to a group containing 
the circle group as a factor. More specifically in this paper I show that the principal 
series representations of $H$ admits a description in terms of equivariant $K$-theory
as above and conjectures are formulated for a description in the same
spirit of other $H$-modules attached to nilpotent elements.
The idea to use equivariant K-theory to study
affine Hecke algebras was subsequently developed in the papers \cite{67, 72} 
(with Kazhdan) and in [Chriss and Ginzburg, Representation theory and complex geometry, 1997]. 
The same idea was later used 

-by Garland and Grojnowski (and by Varagnolo and Vasserot) to realize the 
Cherednik (double affine Hecke) algebra;

-by Nakajima to realize geometrically an affine quantum group.

\head\cite{73} Cells in affine Weyl groups, II, 1987\endhead
This paper is a continuation of \cite{60}. Let $W$ be a Weyl group or an affine Weyl group. One of the main contributions of this paper is a definition (in terms of
the function $a$ of \cite{60}) of a set $D$ of involutions of $W$ (which I call distinguished involutions). The definition was inspired in part by a conjecture of 
[A. Joseph, J. Algebra 1981] for finite $W$, which in fact follows from the 
results of this paper. For finite $W$, one can identify $D$ with the set of 
Duflo involutions defined in the theory of primitive ideals; but I don't know
a similar identification for affine $W$.
In this paper I show that each left cell contains exactly one element of $D$ and
that the set of left cells in $W$ is finite (hence $D$ is also finite). Note that the
set of left cells in a more general Coxeter group can be infinite,  
see [R.B\'edard, Commun. in Alg. 1986 and 1989].
The second main contribution of this paper is the definition of the
asymptotic Hecke ring $J$ of $W$. This is a $\ZZ$-module with basis $\{t_w;w\in W\}$
in which the multiplication constants are obtained from those
for the new basis \cite{37} of the Hecke algebra by making $q$ tend to $0$ (in a strange way, involving the $a$-function of \cite{60}). 
It is not immediately clear that $J$ is associative (it is so, due to \cite{43}); this ring has a rather non-obvious
unit element namely $\sum_{d\in D}t_d$. (Here the
finiteness of $D$ is used). It is also shown that the Hecke algebra admits a natural
algebra homomorphism into the algebra $J$ with scalars suitably extended (this
is again based on \cite{43}).

\head\cite{77}  Leading coefficients of character values of Hecke algebras, 1987\endhead
Let $W$ be a Weyl group and let $c$ be a two-sided cell in $W$. 
Let $G$ be the finite group attached to $c$ in \cite{57}. In this paper I show that
the "non-abelian Fourier transform" of \cite{34} can be interpreted as the 
"character table" of the equivariant complexified $K$-theory (commutative) convolution 
algebra $K_G(G)$ (where $G$ acts on itself by conjugation).

In this paper I associate to each left cell $\G$ in $c$ a subgroup $G_\G$ of $G$
defined up to conjugacy so that
the structure of $\G$ as a $W$-module has a simple description
 in terms of the permutation representation of $G$ on $G/G_\G$.
For example if $W$ is of type $E_8$ and $c$ is a two-sided cell of $W$ with
$G=S_5$ the  symmetric group in $5$ letters then there are $7$ types of left
cells in $c$; they correspond to the following $7$ subgroups of $S_5$:
$S_2,S_3,S_4,S_5, S_2\T S_2,S_3\T S_2$ and the dihedral group of order $8$.

\head\cite{78} (with C. DeConcini and C. Procesi) Homology of the zero set of a nilpotent vector field on a flag manifold, 1988\endhead
The work on this paper was done during my sabbatical leave in Rome (1985/86).
Let $\fg$ be the Lie algebra of a connected reductive group over $\CC$, let $N$ be a 
nilpotent element of $\fg$ and let $B_N$ be the variety of Borel subalgebras of $\fg$
that contain $N$. At the time this paper was written it was known that the 
rational homology of $B_N$ is zero in odd degrees. (The most difficult case, that
of type $E_8$, was done by [Beynon and Spaltenstein, J. Algebra 1984] based on computer 
calculation and then in my paper \cite{69} without computer calculation.) In this paper 
we prove a stronger result namely that the integral homology of $B_N$ is zero in
odd degrees and has no torsion in even degrees. The key case is that where $N$ is
distinguished. There are separate proofs for the case of classical groups 
(where we show the existence of a cell decomposition) and in the exceptional 
case (where we are unable to prove the existence of a cell decomposition but
instead we give an alternative argument based on blow ups and downs which in
a sense gives a more precise result than for the classical groups).
It would be interesting to complete the results of this paper by 1) extending
the method used for exceptional groups (connectedness of a certain graph) to
classical groups and 2) showing that the cell decomposition also exists for
exceptional groups. 
In this paper we also show that the Chow group of $B_N$ is the same as the
integral homology. This has the consequence that the $K$-theory of coherent 
sheaves on $B_N$ is computable, which is a necessary ingredient of \cite{140,143} and
also of [Bezrukavnikov and Mirkovic, arxiv:1001.2562].

\head\cite{79} Quantum deformations of certain simple modules over enveloping algebras, 1988\endhead
In 1986, A. Borel wrote to me a letter pointing out the interesting new work of
Jimbo in which quantized enveloping algebras (q.e.a) were introduced. As a result of 
this letter I gave a course (1986/87) at MIT on q.e.a. and this paper came out of it.
In this paper the divided powers $E_i^{(n)}, F_i^{(n)}$ are introduced for 
the first time by replacing the denominator $n!$ of the classical divided powers by a 
$q$-analogue of $n!$ (depending on $i$). The choice of denominator was such that the formulas 
for the action of $E_i^{(n)}, F_i^{(n)}$ on the standard simple modules of quantum $sl_2$ were 
as simple as possible and also the quantum Serre relations can be written in a form which 
is as simple as possible. Using these divided powers, 
in this paper I define a $\QQ[q,q^{-1}]$-form of the q.e.a. (In later papers
\cite{90,91} this was refined to a $\ZZ[q,q^{-1}]$-form which has become one of the
ingredients in the definition of the canonical basis \cite{92}.) Using this I show that
a simple integrable module of a Kac-Moody Lie algebra can be deformed to a module 
over the corresponding q.e.a. This paper also contains 
the first appearance of the braid group action on a q.e.a. at least in the simply
laced case (but the proofs appeared only in \cite{107}).

\head\cite{80} (with D.Kazhdan) Fixed points on affine flag manifolds, 1988\endhead
My motivation for this paper was as follows. Let $G$ be a semisimple adjoint
group over $\CC$ with Lie algebra $\fg$. Since \cite{40} I knew that (conjecturally) the
nilpotent classes of $\fg$ are in bijection with the two-sided cells of the 
affine Weyl group $W_{af}$ of $G^*$ (Langlands dual); moreover experiments showed that
$\dim H^*(B_x)^{A(x)}$  (where $B_x$ is the Springer fibre at a nilpotent $x$ and $A(x)$ is
the =group of components of the centralizer of $x$ in $G$) is equal to 
the number of left cells in the corresponding two-sided cell.  For example if $G$ is of type $E_8$ and 
$x$ is a subregular nilpotent element, then $B_x$ has $8$ irreducible components (all lines), 
$H^*(B_x)^{A(x)}=H^*(B_x)$ is $9$-dimensional and there are $9$ left cells in the
corresponding two-sided cell $c$. If we now take the affine analogue $x'\in \fg((t))$ of a 
subregular nilpotent $x$ in $\fg$ and if we replace $B_x$ by the set
$B'_{x'}$ of Iwahori subalgebras of $\fg((t))$ that contain $x'$ we see that
$B'_{x'}$ has exactly $9$ irreducible components (all lines). So the number of left cells in
$c$ can now be interpreted not as a dimension of a vector space but as a
number of elements in a set attached to $x'$ (the set of irreducible
components of $B'_{x'}$). This gave me some hope of finding an analogous
relation in more general cases. Although this hope remained unfulfilled
it motivated my interest in investigating sets of the form $B'_{x'}$. 
In this paper it is shown that if $N\in\fg((t))$ is, like $x'$ above, regular 
semisimple and topologically nilpotent (that is $\lim N^k=0$ as $k\to\iy$) 
then $B'_N$ (defined as for $x'$) is a finite or countable (but locally finite)
union of projective
algebraic varieties all of the same dimension; moreover if $N$ is in addition
elliptic then $B'_N$ is itself an algebraic variety. In the paper a conjectural
formula for $\dim(B'_N)$ is given and it is shown how to reduce the
proof of this formula to the case where $N$ is elliptic. (The case where $N$ is
elliptic was settled in [R.Bezrukavnikov, Math. Res. Lett. 1996].)
Also, it is shown that if $x\in\fg$ is nilpotent then for an "open dense"
subset $S(x)$ of $x+t\fg[[t]]$, all elements $N\in S(x)$ are regular semisimple
(and of course topologically nilpotent), $\dim B'_N=\dim B_x$ and the conjugacy 
class in the Weyl group which parametrizes the Cartan subalgebra of $\fg((t))$
containing $N$ depends only on $x$. 
(For example if $x\in\fg$ is subregular nilpotent then $x'$ can be taken to be
an element of $S(x)$; note that $x\to N$ is an affine analogue of the process of 
induction \cite{35}. This gives a map $\Psi$ from nilpotent orbits of $\fg$ to
the set of conjugacy classes in the Weyl group. In the paper this map is 
described explicitly in type $A$ (where it is a bijection) and in the cases 
arising from a nilpotent element of $\fg$ whose centralizer in $G$ is connected, unipotent. 
For example if $G$ has type $E_8$ and $x$ is regular/subregular/subsubregular then
$\Psi(x)$ contains an element of order $30/24/20$. The map $\Psi$ was 
later computed for $G$ of type $B,C,D$ in [N.Spaltenstein, Ast\'erisque 168(1988)], 
[N.Spaltenstein, Arch. Math. (Basel) 1990], and in many cases in 
exceptional types in [N.Spaltenstein, Adv. Math. 1990]. But there are several cases in
exceptional groups where $\Psi(x)$ remains uncomputed. 
In \cite{207} another map between the same two sets is defined using completely
different considerations (based on \cite{199} where a map in the opposite direction
is defined using properties of Bruhat decomposition). The map in \cite{207} is 
computable in all cases and I expect it to be the same as $\Psi$. 
The varieties $B'_N$ introduced in this paper play a key role in the work of 
Ngo B.C. on the fundamental lemma. 

I would like to state the following problem. Let $x\in\fg$ be a distinguished
nilpotent element and let $N\in S(x)$ (so that $B'_N$ is a well defined algebraic 
variety containing $B_x$, see Cor.2 in Sec.3 and 9.2). Let $X_N$ be the
set of irreducible components of $B'_N$. Show that $A(x)$ acts naturally on $X_N$,
that $\dim H^*(B_x)^{A(x)}=\card(X_N/A(x))$ and that the set of left 
cells in the two-sided cell of $W_{af}$ corresponding to $x$ is in natural 
bijection with $X_N/A(x)$.
For example if $G$ is of type $G_2$ (resp. $F_4$) and $x$ is subregular then $A(x)=S_3$ 
(resp. $S_2$) and $\dim H^*(B_x)^{A(x)}=3$ (resp. $5$); on the other hand, $B'_N$ is a 
Dynkin curve of type affine $E_6$ (resp. affine $E_7$) which has a natural $S_3$-action 
(resp. $S_2$-action) whose fixed point set on the set of irreducible components
has cardinal $3$ (resp. $5$).

\head\cite{81} Cuspidal local systems and graded Hecke algebras, I, 1988\endhead
The work on this paper was done in late 1987. In this paper I introduce a graded analogue 
of affine Hecke algebras (with possibly unequal parameters) associated to any root 
system. (After this paper was written I learned of the paper [Drinfeld, Funkt. Anal. Appl. 1986] 
where a similar algebra was introduced for a root system of type A and with the 
grading being disregarded.)
Another new idea of this paper is to define equivariant homology. While in Borel's
definition of equivariant cohomology with respect to an action of an algebraic
group $G$, any classifying space of $G$ can be used, the definition that I give
for equivariant homology is more subtle: it exploits the fact that the classifying 
space of $G$ can be approximated by smooth varieties. (The same idea appeared 
independently in the definition of equivariant derived category given in [Bernstein 
and Lunts, LNM, Springer Verlag 1994].)
This paper contains also a new application of the theory of character sheaves. 
Originally this theory was supposed to provide a machine to compute the values of
the irreducible characters of a reductive group over a finite field. But in this
paper character sheaves (cuspidal with unipotent support) are used to
construct geometrically representations of a graded Hecke algebra (which ultimately 
leads to representations of a $p$-adic group \cite{123,155}). In fact, in this paper I
give a geometric realization of certain graded Hecke algebras in terms of
equivariant homology of a space with group action and with a local system associated
to a cuspidal local system with unipotent support.

\head\cite{84} Modular representations and quantum groups, 1989\endhead
The work on this paper was done in the spring of 1987 and the results were presented at a
US-China conference at Tsinghua University, Beijing in the summer of 1987.
This paper introduces a new concept: that of the quantum group $U_\zeta$ ($\zeta$ is a primitive
$m$-th root of $1$ in the complex numbers) obtained from the $\QQ[q,q^{-1]}]$-form of the quantum 
group introduced in \cite{79} (which involves $q$-analogues of divided powers) by specializing $q=\zeta$.
This paper also formulates the idea (new at the time) that, in the case where $m$ is a prime number $p$,
the representation theory of $U_\zeta$ is governed by laws similar to those of the rational representation
theory of a semisimple algebraic group $G$ over a field of characteristic $p$. Most of the paper is concerned 
with providing evidence for this idea. For example, I prove an analogue for $U_\zeta$ of the Steinberg
tensor product theorem [Steinberg, Nagoya J.Math. 1963]. Thus I show that a simple module of $U_\zeta$ 
with highest weight $\lambda=\lambda_0+p\lambda_1$ (where $\lambda_0$ has coordinates strictly less than $p$) 
is the tensor product of the simple module of $U_\zeta$ with heighest weight $\lambda_0$ with a $U_\zeta$-module 
which may be viewed as the simple $U_1$-module with highest weight $\lambda_1$. (Implicit in this statement is 
the existence of a "quantum Frobenius homomorphism" from $U_\zeta$ to the classical enveloping algebra $U_1$ 
which is also one of the main new observations of this paper.) The key to this tensor product theorem is 
the following property of the Gaussian binomial coefficients specialized at $\zeta$: if $N,R$ are integers
and $N=N_0+pN_1,R=R_0+pR_1$ where $N_i,R_i$ are integers, $0\le N_0\le p-1,0\le R_0\le p-1$, then 
$[N,R]=[N_0,R_0](N_1,R_1)$ for $q=\zeta$. Here $[N,R],[N_0,R_0]$ are Gaussian binomial coefficients and $(N_1,R_1$) is an ordinary binomial 
coefficient. Similarly, the key to the classical Steinberg theorem is the following congruence
(which I learned in my student days from Steenrod's book "Cohomology operations"): 
if $N,R$ are integers and $N=N_0+pN_1+p^2N_2+...,R=R_0+pR_1+p^2R_2+...$ where $N_i,R_i$ are integers, 
$0\le N_i\le p-1,0\le R_i\le p-1$, then $(N,R)=(N_0,R_0)(N_1,R_1)(N_2,R_2)...\mod p$. In this paper I also formulate a
conjecture describing the character of an irreducible finite dimensional $U_\zeta$-module in terms of the 
polynomials \cite{37} attached to the affine Weyl group of the Langlands dual group, similar to the conjecture
that I stated in \cite{40, Problem IV}. This conjecture (which is now known to hold) is one of the steps in
the solution of Problem IV in \cite{40}.

\head\cite{86} Cells in affine Weyl groups IV, 1989\endhead
One of the main results of this paper is establishing a bijection between the set of unipotent classes in
a connected reductive group $G$ over $\CC$ and the set of two-sided cells in the (extended) affine Weyl $W$
group associated to the dual group $G^*$. (This was conjectured in \cite{40}.) The proof uses the earlier parts
of this series (especially the study of the $J$-ring associate to $W$) and the results of of \cite{72}.
Assume now that $G$ is simply connected.
One of the main contributions of this paper is the formulation (see 10.5)  of a (conjectural) basis preserving ring
isomorphism between the $J$-ring of $W$ and the direct sum over the unipotent classes of $G$ of certain
equivariant $K$-groups of certain finite sets attached to the unipotent classes.
I arrived at this conjecture after doing many explicit computations for rank $2$ and using an analogy with
finite Weyl groups \cite{77}.

A weaker form of this conjecture is proved in [Bezrukavnikov, Ostrik, in 
Adv. Studies Pure Math.40 Mat. Soc. Japan 2004]; for type $A$ the conjecture is proved in full in
[Xi, Mem. Amer. Math. Soc. 157(2002)].
A consequence of the conjecture (see 10.8) gives a conjectural bijection between the set 
of dominant weights of $G$ and the set of pairs consisting of a unipotent class of $G$ and an irreducible
rational representation of the centralizer of an element in that class.  This has now been established in 
[Bezrukavnikov, Represent. Th. 2003] with an important 
contribution by [Ostrik, Represent. Th. 2000]. 

\head\cite{88}. On quantum groups, 1990\endhead
This paper (written in early 1989) consists of two parts. In the first part it is shown that from a quantum group
associated to a positive definite symmetric Cartan matrix one can recover in a natural way the Hecke algebra 
attached to the same Cartan matrix. Namely, an explicit construction of the $q$-analog of the adjoint representation
is given (together with an explicit basis which can now be interpreted as the canonical basis \cite{92} of that
representation) and it is shown that the braid group acts naturally on this representation so that the induced
action on the $0$-weight space satisfies the relations of the Hecke algebra. In the second part two conjectures are
formulated. Conjecture 2.3 predicts an equivalence of categories between a certain category $C$ of representations
of a quantum group at a root of $1$ and a certain category $C'$ of representations of an affine Lie algebra at a 
negative central charge related to the order of the root of $1$. (This conjecture was later proved in
\cite{108,109,115,116}.) Conjecture 2.5(b) (resp.2.5(c)) predicts a character formula for the simple objects in $C'$
(resp. $C$) in terms of the polynomials \cite{37} for an affine Weyl group analogous to a conjecture I made in \cite{40} for 
modular representations of a semisimple group in characteristic $p$. Conjecture 2.5(b) has been already stated in
\cite{84} but the present paper suggested that one could prove it if one could prove Conjectures 2.3 and 2.5(c). 
Eventually that was indeed the way that Conjecture 2.5(b) was proved. (Conjecture 2.5(c) was proved by [Kashiwara
and Tanisaki, Duke Math.J. 1995]).

\head\cite{90} Finite dimensional Hopf algebras arizing from quantized universal enveloping algebras, 1990\endhead
This paper was written in the spring of 1989. Let $A=\ZZ[v,v^{-1}]$. This paper introduces a new object: the $A$-form 
${}_AU$ and ${}_AU^+$ of a quantized enveloping algebra $U$ of simplylaced type and its plus part $U^+$. While the 
definition does not need new ideas (compared to the definition of the $\QQ[v,v^{-1}]$-form of $U$ or $U^+$, already 
introduced in \cite{79}) the problem that arises is to show that one gets a well behaved object, for example that 
${}_AU^+$ is a "lattice" in $U^+$. This property is established in the present paper by constructing an $A$-basis for 
${}_AU^+$ which is also a basis of $U^+$. This basis is defined in terms of the braid group action introduced in 
\cite{79}. The proof relies on explicit calculations involving (in particular) the roots of $E_8$.
In this paper I also introduce for an integer $N>0$ a new Hopf algebra of finite dimension 
$N^{\text{number of roots}}$ which can be viewed as the Hopf algebra kernel of the quantum Frobenius
map of \cite{84}. This Hopf algebra is sometimes referred to as the "small quantum group".
This paper is a step toward the construction of the canonical basis of $U^+$ which was achieved in \cite{92}.
Indeed, the lattice ${}_AU^+$ is one of the key ingredients in the definition 
of the canonical basis of $U^+$ given in \cite{92}.

\head \cite{89} Green functions and character sheaves, 1990\endhead
I got the main idea for this paper during a visit at the College de France (May 1988) where I gave a series of 
lectures on character sheaves. The paper was completed in the fall of 1988 when I was visiting IAS, Princeton. 
This paper is a step in the program (initiated in \cite{64, p.226}) of relating (for a connected reductive group
$G$ defined over $F_q$ of characteristic $p$), the characters of representations of $G(F_q)$ and the characteristic 
functions of character sheaves on $G$ which are "defined" over $F_q$. A part of this program would be to show that the Green functions of
$G(F_q)$ (defined in \cite{22}) can be expressed in terms of character sheaves. In this paper I show that this is indeed
so assuming that $q$ is large (no restriction on $p$). The corresponding result for large $p$ was known at the time
(it could be deduced from the work of Springer and Kazhdan). The assumption that $q$ is large enough was later 
removed by [Shoji, Adv.in Math.1995]. Moreover, in this paper it is shown that the "generalized Green functions"
associated to the "induction" functor $R_{L,P}^G$ of \cite{24} can be expressed in terms of character sheaves assuming 
that $q$ is large enough and $p$ is good. This was new even for large $p$. In fact the assumption that $p$ is good can 
now be removed in view of the cleanness property \cite{204}. The methods and results in this paper were
used in [Shoji, Adv.in Math. 1995 and 1996] to study my conjecture \cite{64, p.226} on the relation of
irreducible characters of $G(F_q)$ and character sheaves. 

\head\cite{91}. Quantum groups at roots of $1$, 1990\endhead
In this paper the definition of the braid group action in \cite{79}, the results of \cite{90} about
the $\ZZ[v,v^{-1}]$-form of $U^+$ ("lattice property") and the definition of the small quantum group
in \cite{90} are extended to the nonsimplylaced case. The 
case of $G_2$ was particularly complicated since (unlike the other rank two cases) there are no simple
explicit formulas for the commutation of two divided powers of "root vectors" and for this
reason the argument becomes involved. Also the quantum Frobenius homomorphism 
(which is almost explicit in \cite{84}) is made explicit. 
In the Appendix (joint work with M. Dyer) the "Poincar\'e-Birkhoff-Witt basis" of $U^+$
corresponding to any reduced expression of the longest Weyl group element is introduced, using the
braid group action and the computations in rank $2$ from the main body of the paper. Note that the
basis introduced in \cite{89} is a special case of this PBW basis; the appendix allows one to simplify
some arguments in \cite{89}. Later, these PBW bases turned out to be another of the key ingredients in the 
definition of the canonical basis of $U^+$ (in the simplylaced case) given in \cite{92}.

\head \cite{92} Canonical bases arising from quantized enveloping algebras, 1990\endhead
The results of this paper were obtained while I was giving a course (MIT, fall 1989) on 
quantum groups and in particular on Ringel's work [Ringel, Hall algebras and quantum groups, 
Inv. Math. 1990] and were presented in that course. These results were announced at a
conference on Algebraic Groups in Hyderabad, India, in December 1989 (I did not attend that
conference, but I asked Roger Carter to present my results there).

This paper introduces a rather miraculous object: the canonical basis for $U^+$, 
the plus part of a quantized enveloping algebra $U$ of type $A,D,E$. 
In the paper this is done by two methods (which lead to the same basis):

(1) an algebraic one based on the following three ingredients: 

(i) an integer form of $U^+$ which I introduced earlier \cite{79,90}, 

(ii) a bar involution of $U^+$ and 

(iii) a basis at infinity of $U^+$ coming from any PBW basis, see \cite{91} (remarkably, 
the basis at infinity defined by a PBW basis is independent of the PBW basis);

(2) a topological method based on the local intersection cohomology of 
the orbit closures in the moduli space of representations of a quiver.
\nl
Now even in the approach (1), there is a (minimal) use of the elementary representation
theory of quivers (not intersection cohomology and not in the statements but in the 
proofs). Note that (1) (resp. (2)) 
bear some superficial similarity with things which appeared in the study of
Hecke algebras \cite{37} (resp. \cite{39}); in that study the role of PBW bases is 
played by the (single) standard basis of the Hecke algebra. 
One of the remarkable properties of the canonical basis is that it induces 
a basis in each finite dimensional irreducible module of $U$. 
This paper introduces also a natural piecewise linear structure for the 
canonical basis that is, a finite collections of bijections of the canonical 
basis with $\NN^n$ ($n$=number of positive roots) so that any two 
of these bijections differ by composition with a bijection of $\NN^n$ with itself 
given by a composition of operations which involve only the sum or difference 
of two numbers or the minimum of two numbers. Later, I found that exactly the 
same pattern appears in a rather different context: the parametrization of the 
totally positive semigroup attached to a group of type $A,D,E$, see \cite{119}. A 
similar pattern exists in the nonsimplylaced case, see \cite{193}.
Theorem 8.13 gives what is I believe the first purely combinatorial formula 
for the dimension of a finite dimensional irreducible representation (and its 
weight spaces) of a simple Lie algebra of type $A,D,E$ (it expresses the 
dimension as the result of counting the number of elements of an explicit set 
defined using the piecewise linear structure above); the previously known 
dimension formulas gave the dimension as a ratio of two integers which is 
not obviously an integer (Weyl) or as a difference of two integers which is 
not obviously positive (Kostant). Subsequently, another purely combinatorial 
formula was found by Littelman using his paths. The remarks on Fourier 
transform in Sec.13 are a precursor of \cite{97}.

The transition matrix between the canonical basis of $U^+$ and a fixed PBW basis
of $U^+$ attached to a reduced expression of the longest Weyl group element
has entries which are positive. This is proved in the paper for certain
 special reduced expressions
(adapted to an orientation of the Coxeter graph) when these entries are
interpreted as local intersection cohomology of orbit closures. The same statement
for an arbitrary reduced expression is proved in [Syu Kato, arxiv:1203.5254].
Another proof (relying on the positivity property of the comultiplication proved in \cite{97}) 
appears in [H.Oya, arxiv:1501.01416.]

\mpb

{\it Historical remark.} After this paper became available, another proof of the existence 
of the canonical basis (valid also for Kac-Moody Lie algebras) was given by
[Kashiwara, Duke Math.J. 1991] by a purely algebraic method which uses some ideas (see (i),(ii))
from my paper; and in \cite{97} which generalizes the intersection cohomology approach (2). 
For the history of this subject see the paper of N.Enomoto 
and M.Kashiwara, Symmetric crystals and affine Hecke algebras of type 
$B$, Proc.Japan Acad. 82, Ser.A, 2006 (see page 133) which contains the following statement: 
``We call it a (lower) global basis. It is first introduced by G.Lusztig [5] under 
the name canonical basis for the ADE cases.'' 

\head\cite{95} Canonical bases arising from quantized enveloping algebras, II, 1990\endhead
Some time after my paper \cite{92} became available, Kashiwara found a
different approach to the canonical basis of \cite{92} in which he preserved 
two of the ingredients in my definition ((1) it is contained in the 
``integral part'' of the quantum group and (2) is fixed by a bar operator) 
but replaced ny third ingredient with a different one which made sense in the
more general context of Kac-Moody type). One of the results of this paper
was that for ADE types, Kashiwara's definition gives the same result as my
original definition. Another contribution of this paper is the definition of
a new variety attached to an arbitrary graph. It is shown
that these varieties are equidimensional. (Each one is in fact a Lagrangian 
variety in a symplectic vector space.) Moreover the union $Z$ of the sets of 
irreducible components of these varieties is endowed with certain geometrically 
defined maps $E_k:Z\to Z$ (see 8.8); $k$ is a vertex of the graph. In this paper 
it is conjectured (10.2) that the crystal graph of the plus part of the 
quantized enveloping algebra corresponding to the graph can be geometrically 
realized as the set $Z$ together with our maps $E_k:Z\to Z$ (there are also maps 
$F_k:Z\to Z$ but they are essentially inverse to $E_k$ hence they need not be 
separately constructed). This conjecture was proved by [Kashiwara and Y.Saito, 
Duke Math.J. 1997].

\head\cite{97} Quivers, perverse sheaves and quantized enveloping algebras, 1991\endhead
Let $U^+$ be the plus part of the quantized enveloping algebra corresponding to
a symmetric Cartan matrix $C$. After writing the paper \cite{92} on the canonical
basis of $U^+$ in the case where $C$ is positive definite, I tried to consider
the similar problem for a general $C$. The main problem was to find an appropriate 
definition for the class $X$ of irreducible perverse sheaves on the space of 
representations of fixed dimension $D$ of a quiver attached to $C$ which should 
constitute the canonical basis. If $C$ is positive definite, $X$ consists of all
$G$-equivariant simple perverse sheaves ($G$=product of $GL_n$'s); but in the indefinite case
there are infinitely many $G$-equivariant simple perverse sheaves which is not 
what $X$ should be. I first tried \cite{95} to define $X$ by imposing in addition to
$G$-equivariance a condition on the singular support namely that it should
be contained in the explicit Lagrangian variety $\Lambda$ defined in \cite{95}. 
But I was not able to develop the theory from this definition. Instead I adopted a definition from
the theory of character sheaves, namely $X$ is defined as the collection of
simple perverse sheaves which appear (up to shift) as direct summands
of the direct image of the constant sheaf under the projection maps from certain 
spaces which consists of a representation of dimension $D$ of the quiver and a 
"flag" of a fixed type compatible with the representation. This makes $X$ finite
for any prescribed $D$. With this definition the collection of the various $X$
when $D$ varies can be viewed as a basis of an algebra over $\ZZ[q,q^{-1}]$ in which
multiplication is an analogue of induction of character sheaves ($q$ appears
as the shift). In this paper I prove that the resulting algebra is a
$\ZZ[q,q^{-1}]$-form of $U^+$ and that the basis of $U^+$ provided by the perverse sheaves does not
depend on the orientation of the quiver; hence it is a canonical basis of $U^+$.
I also show that this algebra has something close to a comultiplication
(it is defined as an analogue of restriction of character sheaves).
The structure constants of both the multiplication and "comultiplication" are
in $\NN[q,q^{-1}]$.
Another result of this paper is a new realization of the algebra $U^+$ (for $v=1$)
in terms of convolution of certain constructible functions 
on the Lagrangian variety $\Lambda$ (as above).
This realization actually plays a role in the proofs in this paper.

\head\cite{98}. (with J.M.Smelt) Fixed point varieties in the space of lattices, 1991\endhead
Let $V$ be a vector space of dimension $n$ over $\CC[[\e]]$ with a basis
$e_1,...,e_n$. Let $I$ be the space of Iwahori subalgebras of $SL(V)$ (an affine flag
manifold). Let $N$ be the linear map from $V@>>>V$ such that $N(e_i)=e_{i+1}$ for
$i=1,...,n-1$, $N(e_n)=\e e_1$. Let $t>0$ be an integer relatively
prime to $n$. In this paper we study the space $X_t=\{B\in I;N^t\in B\}$ (by \cite{80}, $X_t$ is 
a projective algebraic variety over $\CC$). It is shown that the Euler characteristic 
of $X_t$ is $\c(X_t)=t^{n-1}$ and that $X_t$ can be paved with affine spaces. After this
paper appeared, I defined a generalization of $N^t$ for any simple Lie algebra
$\fg$ over $\CC$; namely for an integer $t\ge1$ prime to the Coxeter number $h$ we write $t=ah+b$, 
$1\le b\le h-1$ and let $N_t=\e^a\sum_{\a:\text{ root}}c_\a e_\a$ where $e_\a$ are the root vectors
and $c_\a=1$ if the height of $\a$ is $b$, $c_\a=\e$ if the height of $\a$ is $h-b$,
$c_\a=0$ if the height of $\a$ is not $b$ or $h-b$; then $N_t$ is a topologically
nilpotent regular semisimple elliptic element of Coxeter type. Let $X_t$ be the variety
of Iwahori subalgebras of $\fg[[\e]]$ that contain $N_t$ (a projective variety); I conjectured that
the Euler characteristic of $X_t$ is $\c(X_t)=t^{\text{rank}(\fg)}$, which in type $A$ reduces 
to the formula in this paper. This conjecture was proved in [Fan, Transfor. Groups, 1996].
The result on paving was generalized in [Goresky, Kottwitz and MacPherson, Represent. Th.,2006]. In this paper there is also an
explicit formula for the Euler characteristic in the case where the space of Iwahori
subalgebras is replaced by that of maximal parahoric algebras (type $A$); this was
generalized to arbitrary $\fg$ in [Sommers, Nilpotent orbits and ...(Ph.D.Thesis at MIT), 1997].
The formula for $\c(X_t)$ in this paper plays a role in [Berest, Etingof and Ginzburg, IMRN, 2003]. 
The variety $X_t$ (type $A$) and its paving in this paper also plays a role in [Laumon, Fibres de 
Springer et jacobiennes compactifi\'ees, Springer 2006]. 

\head\cite{100} A unipotent support for irreducible representations, 1992\endhead
Let $G$ be a connected reductive group defined over a finite field $F_q$ of sufficiently large characteristic. 
For any unipotent element $u\in G(F_q)$ let $\G_u$ be the generalized Gelfand-Graev representation (GGGR)
associated by Kawanaka to $u$; this is a representation of $G(F_q)$ whose character is zero outside the unipotent set. Let 
$\r$ be an irreducible complex representation of $G(F_q)$; let $\r'$ be the representation of $G(F_q)$ which is dual to 
$\r$ in the sense of \cite{47}. In \cite{57, 13.4} a unipotent conjugacy class $C$ of $G$ was attached to $\r$. In this paper the 
following properties of $C$ are proved (see Theorem 11.2).

(i) The average value of the character of $\r$ on $C(F_q)$ is nonzero and $C$ is characterized by having
maximum dimension among unipotent classes with this property.

(ii) If $g\in G$ is such that $\tr(g,\r)\ne0$ then the unipotent part of $g$ lies in $C$ or in a conjugacy
class of dimension $<\dim C$.

(iii) For some $u\in C(F_q)$, $\r'$ appears with non-zero multiplicity in $\G_u$; for any $u\in C(F_q)$, $\r'$ appears 
with small multiplicity in $\G_u$; if $C'$ is a unipotent class in $G$ such that $\dim(C')>\dim(C)$ or $\dim(C')=\dim(C)$,
$C'\ne C$, then $\r'$ does not appear in $\G_u$ for $u\in C'(F_q)$.
\nl
Note that something close to (i) has been conjectured in \cite{40}; (ii) has been hinted at in \cite{76,p.177,line 13};
(iii) has been conjectured by Kawanaka. It is natural to call $C$ the unipotent support of $\r$. 
One of the keys to the proof of (i)-(iii) is Theorem 7.3 of this paper which gives an explicit decomposition of a
GGGR in terms of intersection cohomology complexes of closures of unipotent classes with coefficients in various 
local systems. A step in the proof of this theorem is a formula for the Fourier transform of a GGGR viewed as a 
function on $Lie(G(F_q))$, involving a Slodowy slice. The connection between GGGR and Slodowy slices found in this
paper is perhaps related to the observation made several years later by [Premet, Special transversal slices ...,
Adv.in Math.2002] that a $W$-algebra (a characteristic zero analogue of the endomorphism algebra of a GGGR) is a 
quantized version of the coordinate ring of a Slodowy slice. 
In this paper we also give a (provisional) definition (see Theorem 10.7) of the unipotent support of a character 
sheaf on $G$. The actual definition (partly conjectural) is given in \cite{212}.

\head\cite{104} Affine quivers and canonical bases, 1992\endhead
In this paper I fix an affine quiver of type $A,D$ or $E$ (but not $A_{2n}$)
with one of the two orientations in which every vertex is a sink or a source.
In this case I construct explicitly the perverse sheaves on the
space of representations of fixed dimension of the quiver which constitute the
canonical basis introduced in \cite{97}. Unlike in the finite type case, these
perverse sheaves can be higher dimensional local systems on an open subset
of their support (the dimension is that of an irreducible representation of a symmetric group).
Also, I describe explicitly (enumerate) the irreducible components of the
Lagrangian variety $\Lambda$ attached in \cite{95} to the affine quiver and show 
that they are
in natural bijection with the perverse sheaves in the canonical basis.
In this paper, the affine quivers are studied in terms of a finite subgroup
of $SL_2(\CC)$ (MacKay correspondence) and I reprove from this point of view the
classification of the indecomposable representations of this quiver,
which goes back in various degrees of generality to Weierstrass and Kronecker
(affine $A_1$), Gelfand and Ponomarev (affine $D_4$), Donovan and Freislich, Nazarova
and [Dlab and Ringel, Memoirs AMS, 1976].
Another result of this paper is the construction of a new basis of the
algebra $U^+$ (with $v=1$) attached to our quiver (later called the semicanonical basis \cite{147})
in which the basis elements appear as constructible functions on the Lagrangian variety $\Lambda$.

\head\cite{110} Coxeter groups and unipotent representations, 1993\endhead
This paper contains things that I did in 1982. One of the results of the classification 
\cite{57} of unipotent representations of a Chevalley group over $F_q$ was that the set of unipotent 
representations depends only on the Weyl group $W$, not on the underlying root system or Chevalley
group. Therefore one can asks whether the set of unipotent representations makes sense if $W$ is
replaced by a finite Coxeter group when the root system and the Chevalley group are not defined.
(One indication that this may be true was provided by the results of \cite{49} which computed
what should be the degrees of the principal series of unipotent representations in type $H_4$
and these degrees turned out to be polynomials in $q$.)
This question is answered in this paper: the set of unipotent representations is
attached to the finite Coxeter group $W$ by heuristic considerations by postulating certain properties
that this set should have which are known in the crystallographic case and showing that these postulates
have a unique solution in the general case. The degrees of the unipotent representations are computed
(extending the results of \cite{49}) and 
the classification of representations in families, the classification of unipotent cuspidal representations
are given in each noncrystallographic case.
For example if $W$ is of type $H_4$ there are $104$ unipotent
representations of which $50$ are cuspidal; the largest family contains $74$ representations of degree
$cq^6+$higher powers of $q$ where $c$ is an algebraic integer (independent of $q$) divided by $120$. 
This result has been found independently by Brou\'e and Malle (unpublished).
It has become a part of a heuristic theory (Brou\'e, Malle, Michel) of 
unipotent representations associated to complex reflection groups.
In the case of finite Coxeter group this theory is no longer heuristic:
it now has a concrete meaning described in \cite{226} in terms of J-rings
\cite{73}.

\head\cite{111} (with I.Grojnowski) A comparison of bases of quantized enveloping algebras, 1993\endhead
At the end of 1991 there were two definitions of a canonical basis of the plus part $U^+$
of the quantized enveloping algebra of a Kac-Moody Lie algebra with symmetric Cartan  
matrix: the algebraic one in [Kashiwara, Duke Math.J. 1991] and a topological one in \cite{[97}. 
(But it was already known that, for finite types, both these definitions agree
with the original definition \cite{92}, see \cite{95},\cite{97}.) In this paper it is shown that these
two bases agree in the general case. The new idea of this paper is a geometric interpretation
of the symmetric bilinear form $(,)$ on $U^+$. Namely for $b,b'$ in the basis \cite{97}, it is shown 
that the rational function $(b,b')$ expanded in a power series in $v^{-1}$ has coefficients 
given by the dimensions of the equivariant $Ext$ groups between the equivariant 
simple perverse sheaves which represent $b,b'$. (These $Ext$ groups can be defined along the same
lines as the equivariant homology spaces in \cite{81}.) In particular these coefficients are
natural numbers. 
The direct sum of the equivariant $Ext$ groups above (for various degrees and
various $b,b'$) is naturally an algebra which, by [Varagnolo and Vasserot, arxiv:0901.3992],
coincides with the KLR-algebra introduced combinatorially in [Khovanov and Lauda, arxiv:0803.4121] and 
[Rouquier, arxiv:0812.5023]. 

\head\cite{112} Tight monomials in quantized enveloping algebras, 1993\endhead
In this paper I show that the construction \cite{97} in terms of quivers of the 
canonical basis of the plus part $U^+$ of a quantized enveloping algebra can be 
generalized to the case where the quiver is allowed to have loops (this was
not allowed in \cite{97}). The resulting class of algebras includes the usual $U^+$
but also the classical Hall algebra with their canonical bases. Moreover, the 
plus part of a quantized (Borcherds) generalized Kac-Moody Lie algebra 
as described in [Kang and Schiffmann, Adv. Math. 2006]) is in fact a subalgebra of 
one of our $U^+$, and the canonical basis described in [loc.cit.] is
closely connected with the canonical basis of $U^+$ introduced in this paper,
see [Kang and Schiffmann, arxiv:0711.1948].
In the $U^+$ of this paper there are elements $F_i^{(a)}$ of the canonical basis
indexed by a vertex $i$ of the quiver and a natural number $a$. 
In the case without loops these elements are divided powers
of a single element $F_i$ but in the general case this is not so.
In the paper I conjectured that the elements $F_i^{(a)}$ generate the algebra $U^+$. 
(This was known from \cite{87} in the case without loops and was proved in the paper in the
case where there is only one vertex and any number of loops.) The conjecture
is now proved by [T. Bozec, 2014]. 
Consider now a monomial $m=F_{i_1}^{(a_1)}F_{i_2}^{(a_2)}...F_{i_n}^{(a_n)}$
in the $F_i^{(a)}$. We say that $m$ is tight if it belongs to the canonical basis.
In this paper I give a criterion to determine whether $m$ is tight.
The criterion is in terms of a certain positivity property of a quadratic form.
Using this criterion I show that $m$ is always tight if there is exactly one
vertex and at least two loops. I also investigate the existence of tight
monomials in the loop free case of small rank. It was already known from
\cite{92} that in type $A_2$ all elements of the canonical basis are tight monomials.
In the paper I show that in type $A_3$ there is an abundance of tight
monomials. In some sense (explained in the paper), $80/100$ of the canonical basis
consists of tight monomials; they fall into 8 families indexed by the
various reduced expressions of the longest Weyl group element). Later, in
[N.Xi, Commun. in Alg. 1999], the remaining elements of the canonical basis were
described explicitly in this case; they are not tight monomials.
The tight monomials in type $A_4$ are described in [Y.Hu, J.Ye and X.Yue, J.Alg. 2003].
But in higher rank there are fewer and fewer tight monomials.

\head\cite{122} Quantum groups at $v=\iy$, 1995\endhead
The main contribution of this paper is that the idea of the $J$-ring (an asymptotic version of the Hecke 
algebra) introduced in \cite{73} makes sense in other contexts. In this paper we try to develop this idea
in the case where the Hecke algebra with its canonical basis and its $a$-function is replaced by the modified 
quantum group $\dot U$ with its canonical basis $\dot B$ introduced in \cite{101} and an appropriate
$a$-function on it. This leads to a ring version at infinity $\dot U^{\iy}$ of $\dot U$. In the
paper this is made explicit for quantum groups of finite type and is stated as a conjecture
for affine type. The  conjecture has now been proved for type $A$ in [K. McGerty, 
Int. Math. Res. Not. 2003] and in general in [J. Beck and H. Nakajima, Duke Math.J. 
2004].

\head\cite{126} Braid group actions and canonical bases, 1996\endhead
Let $U$ be the quantized enveloping algebra corresponding to a given root datum.
Let $U^+$ be the plus part of $U$. Let $E_i$ be the standard generators of $U^+$.
Let $T_i$ be the symmetries of $U$ defined in \cite{107, Part VI} and let $B$ be the canonical basis of $U^+$ defined in \cite{107, 14.4}.
In this paper I show that $T_i$ respects $B$ as much as possible. 
More precisely, we have $U^+=(U^+\cap T_i^{-1}U^+)\oplus U^+E_i$ and I show
that the associated projection $U^+\m(U^+\cap T_i^{-1}U^+)$ applies $B$ to a
basis of $U^+\cap T_i^{-1}U^+$ union with $0$. Similarly we have
$U^+=(T_iU^+\cap U^+)\oplus E_iU^+$ and I show that the associated projection $U^+\m(T_iU^+\cap U^+)$ applies $B$ to a
basis of $T_iU^+\cap U^+$ union with $0$. I then show that these bases of
$U^+\cap T_i^{-1}U^+$, $T_iU^+\cap U^+$ correspond to each other under $T_i$.
According to [Baumann, arxiv:1104.0907], an analogus result holds when the 
canonical basis B is replaced by the semicanonical basis \cite{147} assuming that
the root datum is simply laced and v=1.
The results of this paper have been used in
[Beck, Chari and Pressley, Duke Math.J. 1999] to give a characterization of
the canonical basis $B$ of $U^+$ (in the affine case) in terms of a basis $B'$ of 
$U^+$ of PBW type, constructed using (in part) iterations of
symmetries $T_i$; the results of this paper are used to show that any element
of $B'$ is congruent to a unique element of $B$ modulo $v^{-1}$ times the
$\ZZ[v^{-1}]$-lattice generated by $B$. (This extends the results of \cite{92} in the
finite type case.)

\head\cite{131} Notes on unipotent classes, 1997\endhead
Let $G$ be a semisimple almost simple algebraic group over an algebraically
closed field $k$ whose characteristic is $0$ or a good prime.
In this paper I study a partition of the unipotent variety of G into loccally
closed strata, called special pieces. Each special piece is a union of
unipotent classes of which exactly one is special in the sense of \cite{36} (that is
the corresponding Springer representation of the Weyl group is special in the
sense of \cite{36}; the other unipotent classes in the piece are in the closure
of the special class in the piece but not in the closure of any smaller
special piece. The fact that this partition of G is well defined was shown
by Spaltenstein in his book.

In \cite{44} it was conjectured that any special piece is a rational homology manifold. This was later shown to be
true by [Beynon and Spaltenstein, J.Alg. 1984] and [Kraft and Procesi, Asterisque 1989].
In this paper we state a refinement of this conjecture: a special piece
is the quotient of a smooth variety by a finite group. 

In \cite{44} it was also conjectured that the polynomials in $q$ which give
the number of $F_q$-points of a special piece (when $k=\bar F_q$) depends only on
the Weyl group (not on the root system). This conjecture is proved in this
paper by a complicated computation. Another proof of the conjecture based on
Kato's exotic nilpotent cone was given in [Achar, Henderson and Sommers, Repres. Th.,
2011]. A conjectural explanation for why the conjecture should hold was given
in [Geck and Malle, Experimental Math., 1999].

In this paper I give the following characterization
of special pieces which doesn't use the notion of closure: two unipotent
classes belong to the same
special piece if and only if the corresponding Springer representations belong to the same two-sided cell of the Weyl group. 

\head\cite{132} Cells in affine Weyl groups and tensor categories, 1997\endhead
The main conjecture of this paper is proved in [Bezrukavnikov, Adv.Studies 
Pure Math.40, Mat. Soc. Japan 2004].

\head \cite{138} On quiver varieties, 1998\endhead
Theorem 5.5 has been strengthened in [Malkin, Ostrik, Vybornov, Adv.in Math. 2006]
where it is shown that the morphism in that Theorem is in fact an isomorphism of 
algebraic varieties.

\head\cite{148} Fermionic form and Betti numbers, 2000\endhead
This paper contains a conjecture which expresses the Betti numbers of the Nakajima
quiver varieties in terms of a certain complicated but in principle computable
fermionic form. This conjecture has now been proved in [Kodera, Naoi, arxiv:1103.4207].

\head \cite{157} Rationality properties of unipotent representations, 2002\endhead
Let $G$ be a split connected reductive group over $F_q$. For each $w$ in the Weyl
group $W$ of $G$ let $R_w$ be the virtual representation of $G(F_q)$ associated to
$w$ in \cite{22}. Let $r$ be a unipotent representation of $G(F_q)$ that is, an 
irreducible representation appearing in $R_w$ for some $w$. Let $A(r)$ be the set
of $w$ in $W$ such that $r$ appears in $R_w$. Let $A'(r)$ be the set of elements of minimal
length of $A(r)$. One of the main observations of this paper is that if $r$ is
cuspidal then  $A'(r)$ is contained in a single conjugacy class $C(r)$ of $W$ and that 
for $w$ in $C(r)$, the multiplicity of $r$ in $R_w$ is equal to $(-1)^{\text{semisimple rank of }G}$. From 
this it is deduced that a unipotent representations of $G(F_q)$
whose character has values in rational numbers is actually defined over the
rational numbers; in particular if $G$ is of classical type any unipotent 
representation is defined over the rational numbers. (This is not true for
unitary groups over $F_q$). The proof is not
constructive since it uses the Hasse principle for division algebras. It is
also observed that an analogue $A\mapsto C(A)$ of the correspondence 
$r\mapsto C(r)$ holds when $r$ is replaced by a unipotent cuspidal character 
sheaf $A$. For example if $G$ is of type $E_8/F_4/G_2$ and $A$ is the unique 
unipotent cuspidal character sheaf with unipotent support (the closure of the 
conjugacy class $\gamma$ of a unipotent element whose centralizer has group of 
components $S_5/S_4/S_3$) then $C(A)$ contains an element which is "regular" of 
order $6/4/3$ ($=$largest order of an element of $S_5/S_4/S_3$). 
In the paper it is noted that in these three cases $C(A)$ consists of elements
of a single length $40/12/4$. It is
interesting that $C(A)$ also corresponds to $\gamma$ under a quite different
correspondence described in \cite{197}. The rationality property of unipotent
representations described in this paper was known to me (with a different
proof, also explained in the paper) since 1982 
when it was the object of a lecture that I gave at a US-France Conference on Representation Theory in Paris.
The results of this paper were presented at a conference in Rome (June 2001)
and one in Isle de Berder (Bretagne) in September 2001. 

\head\cite{167} An induction theorem for Springer's representations, 2004\endhead
The theorem in the title was stated without proof in \cite{48} for reductive groups
in characteristic zero and it was one of the main tools in the computation in 
\cite{48} of the Springer correspondence for groups of type $E_n$. This paper (written
in 2001) contains a proof of that theorem, valid in arbitrary characteristic. 
It uses the connection between Green functions of a reductive group over a 
finite field and character sheaves \cite{89} and also some arithmetic considerations.

\head\cite{227} Exceptional representations of Weyl groups, 2017\endhead
Springer has discovered that certain representations of the Hecke algebra with parameter $p$
(a prime number) of an irreducible Weyl group $W$ cannot be defined over the rational numbers,
contradicting an assertion of Benson and Curtis. (This occurs only if $W$ is of type $E_7$ or $E_8$.)
In this paper I study various properties of this kind of ``exceptional'' representations. One of the
themes of this paper is that the phenomenon of exceptional representations also occurs in finite 
non-crystallographic Coxeter groups (namely those of type $H_3$ and $H_4$). For example I observe that for 
$W$ of type $E_7,E_8,H_3$ or $H_4$, the number of exceptional representations of $W$ times the dimension of 
any exceptional representation of $W$ is equal to the largest power of $2$ dividing the order of $W$.
I also show that the fake degree of a nonexceptional representation of $W$ is
given by a palindromic polynomial (this was earlier known by a case by case argument.)   

\head\cite{228} Action of longest element on a Hecke algebra cell module, 2015\endhead
In this paper I studied the action of the basis element of the Hecke algebra
corresponding to the longest element of the Weyl group on a left cell
module. My goal was to show that this action has a simple form in
terms of the canonical basis of the left cell module, namely that it
is given by a permutation (with square one) of the basis elements 
times a sign and times a power of q. My proof applies in the case in which 
the Hecke algebra is allowed to have unequal parameters. 
After a first version of this paper was written, M. Douglass pointed out to me that 
for Hecke algebras with equal parameters the result was proved earlier by A. Mathas.
(However, for the Hecke algebra of type A, the result was proved even earlier; it is a
consequence of Cor.5.9 in my 1990 paper \cite{95} on quantum groups, applied
to the intersection of $B[d]$ in loc.cit. with the zero weight of a suitable $L_d$.)
The result of this paper has been used by Bonnafe to extend Losev's work on cacti to unequal parameters.

\head\cite{229} On the character of certain irreducible modular representations, 2015\endhead
Let $G$ be a connected simply connected almost simple algebraic group over an algebraically closed 
field of characteristic $p>0$ with a fixed maximal torus $T$ and a fixed Borel subgroup containing it.
Let $X^+$ be the set of dominant characters of $T$. We consider the category C of finite dimensional
rational representations of $G$. The simple objects of C have been classified by Chevalley; up to isomorphism
they are indexed by $X^+$; let $L_\l$ be the simple object indexed by $\l\in X^+$. Let $E^0_\l\in C$ be the
Weyl module indexed by $\l\in X^+$. The Weyl modules form another basis of the Grothendieck group $\cg(C)$ 
of $C$. Hence for any $\l\in X^+$ we have $L_\l=\sum_{\l'}c_{\l',\l}E^0_{\l'}$ where $c_{\l',\l}$ are integers, 
zero for all but finitely many $\l'$. It is of considerable interest to understand the character of each $L_\l$ 
or, equivalenty, to understand the integers $c_{\l',\l}$. (The character of $E^0_\l$ is known; it is
given by Weyl's formula.) In \cite{40}, Problem IV, I stated a conjecture 
which expresses the integers $c_{\l',\l}$ (with $\l$ in a certain finite subset of $X^+$ containg the 
``restricted weights'') in terms of the polynomials $P_{y,w}$ \cite{37} (evaluated at $1$) associated to the 
affine Weyl group of the Langlands dual of $G$, assuming that $p$ is sufficiently large relative to the type of 
$G$. (Then the case of general $\l$ can be deduced by appealing to Steinberg's tensor product theorem.)
This conjecture has been proved for $p$ very large (see the comments to \cite{40}). In this paper I reformulate
the conjecture in such a way that the tensor product theorem is not used. Namely, for any $k\ge1$ I define 
a basis $\{E^k_\l;\l\in X^+\}$ of $\cg(C)$ with the following properties: if $\k\ge0$ the matrix expressing the 
elements of the basis $\{E^{k+1}_\l;\l\in X^+\}$ in terms of the basis $\{E^k_\l;\l\in X^+\}$ has entries 
expressed in terms of $P_{y,w}(1)$ (suitably rescaled) with no restriction on $\l$; for any $\l\in X^+$ we have 
$E^k_\l=E^{k+1}_\l=E^{k+2}_\l=\do$ large $k$ (we denote the common value by $E^\iy_\l$); if $p$ is very large 
relative to the type of $G$ then $E^\iy_\l=L_\l$ for any $\l\in X^+$. But even without an assumption on $p$, 
the basis $\{E^\iy_\l;\l\in X^+\}$ of $\cg(C)$ is well defined. For arbitrary $p$ we can write
$L_\l=\sum_{\l'}\tc_{\l',\l}E^\iy_{\l'}$ where $\tc_{\l',\l}$ are integers, zero for all but finitely many 
$\l'$. Since $E^\iy_{\l'}$ are in principle computable, to understand the character of $L_\l$ is the same as to 
understand the unknowns $\tc_{\l',\l}$. In this formulation the Weyl modules $E^0_\l$ are replaced by 
$E^\iy_\l$ which can be viewed as the ``new Weyl modules''. (At least if $p$ is not very small, $E^\iy_\l$ is 
indeed a module. If $\l=\l_0+\l_1p+\l_2p^2+\do$ with $\l_i\in X^+$ restricted, we have
$E^\iy_\l=E^1_{\l_0}\ot(E^1_{\l_1})^{Fr}\ot\do$ where $Fr$ denotes twisting by Frobenius and
$E^1_{\l_0},(E^1)_{\l_1},\do$ represent reductions modulo $p$ of simple modules over a quantum group at
 a root of $1$.) Thus the problem of computing the character of the $L_\l$'s is decomposed into two
steps: the first one is to compute the character of each $E^\iy_\l$ (this can be regarded as understood);
the second one is to find the transition matrix from $E^\iy_\l$ to $L_\l$ (this is not understood, except for
very large $p$ when it is the identity matrix).

\mpb

We consider the example where $G=SL_2$, $p=2$.
In this case we can identify $X^+=\NN$. We write the elements $E^0_\l$ as a sequence of dimensions:
$$1,\un2,3,\un4,5,\un6,7,\un8,9,\un{10},11,\un{12},13,\un{14},15,\un{16},\do$$
We underline the numbers divisible by $2$. They form mirrors which are used to get the list of elements 
$E^1_\l$. In this new list $1,\un2$ remain the same; $3$ is reflected in the mirror $2$ giving 
$1$; we change $3$ to $3-1=2$. We keep $\un4$ the same; $5$ is reflected in the mirrors $4$ and $2$ giving 
$3,1$; we change $5$ to $5-3+1=3$. We keep $\un6$ the same; $7$ is reflected in the mirrors $6$,$4$ and $2$ 
giving $5,3,1$; we change $7$ to $7-5+3-1=4$. We continue in this way and we find the sequence of dimensions
of elements $E^1_\l$:
$$[1,2],\un{[2,4]},[3,6],\un{[4,8]},[5,10],\un{[6,12]},[7,14],\un{[8,10]},\do$$
We have arranged these elements in groups of two as shown and we have underlined the groups in which each
element is divisible by $2$; these will play the role of mirrors which are used to get the list of elements 
$E^2_\l$.
In this new list, $[1,2],[2,4]$ remain the same; $[3,6]$ is reflected in the mirror $[2,4]$ giving $[1,2]$; we
change $[3,6]$ to $[3,6]-[1,2]=[2,4]$. We keep $\un{[4,8]}$ the same;
$[5,10]$ is reflected in the mirrors $[4,8]$ and $[2,4]$ giving $[3,6],[1,2]$; we change
$[5,10]$ to $[5,10]-[3,6]+[1,2]=[3,6]$. We continue in this way 
and we find the sequence of dimensions of elements $E^2_\l$:
$$[1,2,2,4],\un{[2,4,4,8]},[3,6,6,12],\un{[4,8,8,16]},\do$$
We have arranged these elements in groups of four as shown and we have underlined the groups in which each
element is divisible by $2$; these will play the role of mirrors which are used to get the list of elements 
$E^3_\l$:
$$[1,2,2,4,2,4,4,8],\un{[2,4,4,8,4,8,8,16]},\do$$
We arrange these in groups of eight, etc. 

\mpb

We consider the example where $G=SL_2$, $p=3$.
In this case we can identify $X^+=\NN$. We write the elements $E^0_\l$ as a sequence of dimensions:
$$1,2,\un3,4,5,\un6,7,8,\un9,10,11,\un{12},13,14,\un{15},16,17,\un{18},19,20,\un{21},\do$$
We underline the numbers divisible by $3$. They form mirrors which are used to get the list of elements 
$E^1_\l$. In this new list $1,2,\un3$ remain the same; $4$ and $5$ are reflected in the mirror $3$ giving 
$2$ and $1$; we change $4$ to $4-2=2$ and $5$ to $5-1=4$. 
We keep the underlined $\un6$ the same; $7$ and $8$ are reflected
in the mirrors $6$ and $3$ giving $5,1$ and $4,2$; we change $7$ to $7-5+1$ and $8$ to $8-4+2=6$.
We keep the underlined $\un9$ the same; $10$ and $11$ are reflected 
in the mirrors $9$, $6$ and $3$ giving $8,4,2$ and $7,5,1$; we change $10$ to $10-8+4-2=4$ and $11$ to
$11-7+5-1=8$. We keep the underlined $\un{12}$ the same; $13$ and $14$ are reflected in 
the mirrors $12$, $9$, $6$ and $3$ giving $11,7,5,1$ and $10,8,4,2$; we change $13$ to 
$13-11+7-5+1=5$ and $14$ to $14-10+8-4+2=10$. We continue in this way and we find the sequence of dimensions
of elements $E^1_\l$:
$$[1,2,3],[2,4,6],\un{[3,6,9]},[4,8,12],[5,10,15],\un{[6,12,18]},[7,14,21],\do$$
We have arranged these elements in groups of three as shown and we have underlined the groups in which each
element is divisible by $3$; these will play the role of mirrors which are used to get the list of elements $E^2_\l$.
In this new list, $[1,2,3],[2,4,6],\un{[3,6,9]}$ remain the same; $[4,8,12]$ and $[5,10,15]$
are reflected in the mirror $[3,6,9]$ giving $[2,4,6]$ and $[1,2,3]$; we change
$[4,8,12]$ to $[4,8,12]-[2,4,6]=[2,4,6]$ and $[5,10,15]$ to $[5,10,15]-[1,2,3]=[4,8,12]$. We keep
the underlined $\un{[6,12,18]}$ the same; $[7,14,21]$ is reflected in the mirrors $[6,12,18]$ and $[3,6,9]$
giving $[5,10,15],[1,2,3]$; we change $[7,14,21]$ to $[7,14,21]-[5,10,15]+[1,2,3]=[3,6,9]$. We continue in 
this way and we find the sequence of dimensions of elements $E^2_\l$:
$$[1,2,3,2,4,6,3,6,9],[2,4,6,4,8,12,6,12,18],\un{[3,6,9,\do},\do$$
We have arranged these elements in groups of nine as shown and we have underlined the groups in which each
element is divisible by $3$; these will play the role of mirrors which are used to get the list of elements 
$E^3_\l$:
$$1,2,3,2,4,6,3,6,9,2,4,6,4,8,12,6,12,18,3,6,9,\do$$
We arrange these in groups of $27$, etc.

\head\cite{231} Some power series involving involutions in Coxeter groups, 2015\endhead
In this paper I study an analogue of the Poincare series of a Coxeter group in the case where the summation 
is taken not over the entire Coxeter group but only over the
involutions in the Coxeter group with an additional weight function. 
The main result is that this new Poncare series can be expressed as a quotient of two classical
Poincare series, one with parameter $q^2$ and one with parameter $q$.
The proof depends on properties of the Hecke algebra module in \cite{208},\cite{209}.
A key step in the proof is in common with a proof in \cite{238}.

\head \cite{232} Nonsplit Hecke algebras and perverse sheaves, 2016\endhead
It is well known that the Hecke algebra with equal parameters associated to a Weyl
group can be interpreted geometrically in terms of convolution of perverse sheaves on a
flag manifold. This has important consequences such as positivity for various structure
constants of the Hecke algebra. An interesting question is whether this can be generalized to Hecke algebra
with unequal parameters arising from the induction of unipotent cuspidal representations of
reductive groups over a finite field.
In my 2003 book I conjectured such an interpretation which involved the theory of
parabolic character sheaves which I have introduced in the early 2000's and which 
generalizes the theory of character sheaves. 
In this paper I verify this conjecture in the first non-trivial case, the Hecke
algebra of type $B_2$ with parameters $q,q^3$. The proof involves the explicit
knowledge of the Kazhdan-Lusztig polynomials for type $B_4$. This provides a very compelling
support for the conjecture.

\head\cite{233} (with G. Williamson) On the character of certain tilting modules, 2017\endhead
Let $G$ be a connected, simply connected reductive group over an algebraically closed field 
$\kk$ of characteristic $p>0$. Assume that $p$ is large. In a paper that I wrote in 2014, I 
found that one can define for any $n\ge0$ a class of $C_n$ of
 characters of rational representations of $G$ which for $n=0$ are the Weyl modules, for
$n=1$ they look like the irreducible characters of the corresponding quantum
group at a root of $1$ and for large $n$ they look more and more like the 
irreducible characters of $G$; moreover the transition matrix from $C_n$ to $C_{n-1}$ is
the same as that from $C_1$ to $C_0$ (suitably rescaled). In this paper we try to extend this point of
view from the case of simple modules to the case of tilting modules. 
We succeed in doing so for most but not all indecomposable tilting modules.  (We are using Donkin's tensor
product theorem for tilting modules is valid for most but not all dominant weights).

\head \cite{234} Non-unipotent character sheaves as a categorical centre, 2016; 
\cite{245} Non-unipotent representations and categorical centers, 2017\endhead
Let $G$ be a connected reductive group over an algebraically closed field $\kk$. In the case where
$\kk$ is an algebraic closure of the field $F_q$ with $q$ elements and $G$ has a fixed $Fq$-rational
structure we denote by $G(F_q)$ the group of rational points of $G$ over $F_q$.
The irreducible representations of $G(F_q)$ were classified in my 1984 book. Soon after that I have
introduced the character sheaves on $G$ and classified them. The two classifications are very similar.
One key ingredient in my proof of those classifications was the use of leading coefficients of
characters of the Hecke algebra. Later this approach was formalized in my definition of the J-ring
or asymptotic Hecke algebra and even later in my definition of the categorified version of the J-ring
(a monoidal category). It is therefore not unreasonable that the study of this categorified J-ring can 
provide insight about the classification problems about. Now Joyal-Street, Majid and Drinfeld have
associated to a monoidal category a new category called the center of the monoidal category. (The
relation between these two categories is similar to the relation between a ring and its centre.)
In particular the centre of the categorified J-ring is defined. 
In 2012, Bezrukavnikov, Finkelberg and Ostrik (BFO) showed that in the case where
$\kk$ is the complex numbers, the classification of unipotent character sheaves on $G$ is equivalent to
the classification of simple objects of the centre of the categorified J-ring. They proved a similar
result for the classification of character sheaves on $G$ with fixed central character assuming that $G$ has
connected centre (and again $\kk$ is the complex numbers). This was an important contribution to the
subject, but its deficiency was the assumption of characteristic zero and the assumption of connected
centre. In fact BFO used techniques that are not available in positive characteristic. Since from my
point of view character sheaves are particularly interesting in positive characteristic (due to their
usefulness in computing character tables of finite reductive groups) I was very interested in finding
a proof of a characteristic p analogue of the BFO result. In 2015 I published two papers where I
gave a characteristic $p$ analogue of the BFO result in the case of unipotent character sheaves
and also a result in the same spirit for the unipotent representations of $G(F_q)$ (which has no
counterpart in the BFO approach). The new ingredient was a definition of the truncated convolution
of two character sheaves which involves taking usual convolution but taking a certain fixed perverse cohomology
sheaf of it and the part of a certain fixed weight in it. (Here we taking advantage that we are in positive
charscteristic so that the notion of weight is defined.) One biproduct was an explanation of the known
fact that the unipotent character sheaves and the unipotent representations have the same classification;
they are both indexed by the simple objects in the same categorical centre. There remained the
problem of extending these results to not necessarily unipotent character sheaves or representations (again
in positive characteristic). This problem is resolved in my two papers  \cite{234},\cite{245} in which I
show that the character sheaves on $G$ (with $\kk$ of characteristic $p>0$) and the irreducible 
representations of $G(F_q)$ can be classified in terms of twisted categorical centers of certain
monoidal categories (with an automorphism) defined in terms of monodromic perverse sheaves on $G/U$
($U$ is the unipotent radical of a Borel subgroup of $G$). 
One ingredient in the solution of this problem  was the fact that in my 1984 book \cite{57} I computed the local 
intersection cohomology of the inverse image of a Schubert variety in $G/U$ with coefficients in a 
monodromic local system. Using this computation I was able to 
define an analogue of the $J$-ring and its categorified version when sheaves on the flag manifold are 
replaced by (monodromic) sheaves on $G/U$. The $J$-ring was defined using a new kind of Hecke algebra (which
I call extended Hecke algebra; it was introduced in my 2005 paper \cite{172}.) This algebra has as a specialization the algebra considered by
Yokonuma in the 1960's which describes the endomorphisms of the representation of $G(F_q)$ induced by the
unit representation of $U(F_q)$. This algebra is similar to but not the same as an algebra considered 
in the 1990's by Mars and Springer. For this extended Hecke algebra one can define the notion of canonical
basis, the notion of two-sided cell and the notion of $J$-ring (generalizing the analogous notions for 
the usual Hecke algebra.) One can also define a categorified $J$-ring which decomposes
according to the two-sided cells in the extended Hecke algebra. Moreover the problem of classifying character
sheaves or representations can be decomposed into a separate problem for each two-sided cell. This is the
problem which I solved in \cite{234},\cite{245}.

\head \cite{235} An involution based left ideal in the Hecke algebra, 2017\endhead
In the 2012 paper \cite{208} I and Vogan have introduced a module for the Hecke algebra of a Weyl group
with basis indexed by the involutions in that Weyl group. The results of that paper have been
used by subsequent authors to study the problem of unitarizability of representations of real
reductive groups. In this paper I show that the module above can be realized as a left ideal
of the Hecke algebra generated by a very remarkable element of the Hecke algebra. 
The basis of the module becomes a basis of this ideal indexed by the involutions. By studying this basis
one obtains a canonical surjective map from the Weyl group to the
set of involutions of the Weyl group. Another interesting consequence is a characterization of
the special representations of a Weyl group (of classical type) by a positivity property.

\head \cite{236} Generic character sheaves on groups over $\kk[\e]/(\e^r)$, 2017\endhead
Let $G$ be a connected reductive group over an algebraically closed field $\kk$. Let
$\e$ be an indeterminate and let $r$ be an integer $\ge1$. We consider the group
$G_r$ of points of $G$ over the truncated polynomial ring $\kk[\e]/(\e^r)$. We view $G_r$
as a not necessarily reductive, connected algebraic group of dimension equal to $r\dim G$.
In the case where $r=1$ we have $G=G_1$; in this case I have defined in 1985-1986 a class
of simple perverse sheaves on $G$ which constitute a geometric theory of characters for $G$;
this theory is very useful for the computation of the characters of irreducible representations
of the points of $G$ over a finite field. In a paper in 2006 I observed that 
a theory of character sheaves does not exist for $r\ge2$ and conjectured that at least
a theory of ``generic character sheaves'' should exist; I have proved that this is so for 
$G=SL_2$, $r=2$. This would again provide insight into the characters of generic irreducible representations
of the analogue of $G_r$ over a finite field. In this paper I prove the conjecture in the
case where $r=2$ or $4$ (no restriction on $G$) and give a partial proof in the case where $r=3$.

\head \cite{237} Generalized Springer theory and weight functions, 2017\endhead
Let $G$ be a connected reductive group over an algebraically closed field.
Let $X$ be the set of pairs $(\co,\cl)$ where $\co$ is a unipotent class in $G$ and $\cl$ is an irreducible
$l$-adic local system on $\co$, equivariant for the conjugation action of $G$. In my 1984 paper \cite{59} 
I showed 
that $X$ has a canonical partition $X=\sqc_{i\in I}X_i$ where each $X_i$ is in natural bijection with the
$\Irr W_i$, set of irreducible representations of a Weyl group $W_i$ (which appears as the group of 
components of the normalizer of a Levi subgroup of a 
parabolic subgroup of $G$). For some $i=i_0$ in $I$, $W_i$ is the full
Weyl group of $G$ and the bijection $X_{i_0}\lra \Irr W_{i_0}$ was originally defined by Springer in 1976.
The main contribution of this paper is that each $W_i$ has an additional structure, namely there is a
natural weight function $L_i:S_i@>>>\ZZ_{>0}$ on the set $S_i$ of simple reflections of $W_i$ that is,
 a function which takes the same value on any two simple reflections which are conjugate in $W_i$.

\head\cite{238} On the definition of almost characters, 2017\endhead
Let $G$ be a connected reductive group over an algebraic closure of the finite field $F_q$ with
 a fixed $F_q$-structure. Let $F:G@>>>G$ be the Frobenius map. Let $V$ be the vector space of 
class functions on $G(F_q)$ with values in an algebraic closure of the $l$-adic numbers. Now $V$ has
two bases $B_1,B_2$; here $B_1$ consists of the characters of the irreducible representations of $G(F_q)$
and $B_2$ consists of the characteristic functions of the character sheaves $A$ on $G$ such that $F^*A\cong A$.
(Note that the functions in $B_2$ are defined only up to multiplication by roots of $1$.)
In the case where the centre of $G$ is connected  a third basis $B_3$ (again defined
only up to multiplication by roots of $1$) was described in my 1984 book \cite{57}. The functions in $B_3$ were
called {\it almost characters}. They were explicit linear
combinations of elements in $B_1$ in which the coefficients were given essentially by a non-abelian
Fourier transform. At that time I conjectured that $B_2$ and $B_3$ coincide. (This is known in many cases.)
In this paper I give a definition of $B_3$ without the assumption that the centre of $G$ is connected.
(It answers a question of M. Geck.) The definition again involves something similar to a non-abelian
Fourier transform. It is again expected that $B_2=B_3$.

\head \cite{239} Special representation of Weyl groups: a positivity property, 2017\endhead
Let $W$ be a Weyl group. In my 1979 paper \cite{36} 
I introduced a class $\cs_W$ of irreducible representations of $W$.
Later these were called special representations. They play a key role in the classification of irreducible
representations of finite reductive groups, in the classification of character sheaves, in the classification
of primitive ideals and in other questions of representation theory. In this paper I give a new
characterization of special representations: an irreducible representation of $W$ is special if and only
if the corresponding irreducible representation of the $J$-ring associated to $W$ admits a basis in which
each canonical basis element of the $J$-ring acts by a matrix with all entries in $\NN$. (In the case where
$W$ is of classical type, I have proved this result in an earlier paper \cite{235}.)

\head \cite{240},\cite{242},\cite{243},\cite{254} (with Z. Yun) 
$\ZZ/m$-graded Lie algebras and perverse sheaves I-IV, 2017-2018\endhead 
Let $G$ be a connected, simply connected, semisimple algebraic group over an algebraically closed field $\kk$.
We assume that the characteristic of $\kk$ is either $0$ or a large prime number. Let $\fg$ be the Lie
algebra of $G$. We assume that we are given an integer $m\ge1$ and a $\ZZ/m$-grading $\fg=\op_{j\in\ZZ/m}\fg_j$
compatible with the Lie algebra structure. Let $G_0$ be a closed connected subgroup of $G$ with Lie algebra
$\fg_0$. Now $G_0$ acts by the adjoint action on $\fg_1$ and this induces an action of $G_0$ on $\fg_1^{nil}$,
the variety of elements of $\fg_1$ which are nilpotent in $\fg$; this action has finitely many orbits.
Let $\ci$ be the set of pairs $(\co,\cl)$ where $\co$ is a $G_0$-orbit in $\fg_1^{nil}$ and $\cl$ is an
irreducible $G_0$-equivariant $l$-adic local system on $\co$ (up to isomorphism). This is a finite set.
For $(\co,\cl)\in\ci$ let $\co^\sha$ be the intersection cohomology complex of the closure $\bar\co$ of $\co$ 
with coefficients in $\cl$; if we are also given $(\co',\cl')\in\ci$ such that $\co'\sub\bar\co$ and an integer
$i$, then we can consider the multiplicity $\mu^i_{\co',\cl';\co,\cl}$ of $\cl'$ in the restriction of the 
$i$-th cohomology sheaf of $\cl^\sha$ to $\co'$. 
The numbers $\mu^i_{\co',\cl';\co,\cl}$ contain important representation theoretic information.
When $m=1$ they are closely connected with the generalized Green functions which enter in the character
formulas for Chevalley groups over a finite field. In the limiting case $m=\iy$ these numbers enter in the
multiplicity formulas in standard modules for affine Hecke algebras with possibly unequal parameters.
In the case $1<m<\iy$ these numbers enter (at least conjecturally) in the multiplicity formulas in standard 
modules for half-degenerate double affine Hecke algebras with possibly unequal parameters. For these reasons
and also for purely geometrical reasons it is of interest to find ways to compute the numbers
$\mu^i_{\co',\cl';\co,\cl}$. When $m=1$ these numbers are computable by an algorithm which appears in my
papers on character sheaves. In the limiting case $m=\iy$ an algorithm to compute these numbers was essentially
given in two of my papers \cite{120}, \cite{191} in 1995 and 2010. 
In this paper we essentially compute the numbers
$\mu^i_{\co',\cl';\co,\cl}$ for any $m<\iy$. Our approach is quite different from the old approach for $m=1$;
it has on the other hand much in common with the approach for $m=\iy$ ($\ZZ$-graded case); it also has some
similarities to the theory of canonical bases for quantum groups.
The main ingredient in our proof is the process of ``spiral induction''; this involves a family of proper maps 
into $\fg^{nil}_1$ which allows us to ``induce'' perverse sheaves from the $ZZ$-graded case to the 
$\ZZ/m$ graded case. 
These maps are indexed by the facets of an affine hyperplane arrangement (or rather several such
arrangements). These maps appear to be new even in the case $m=1$ and give enough information to 
compute all the numbers $\mu^i_{\co',\cl';\co,\cl}$. For example we prove that these numbers are $0$ when
$i$ is odd. (This was known earlier when $m=1$ or $m=\iy$.) We also define a partition of $\ci$ into blocks
so that $\mu^i_{\co',\cl';\co,\cl}=0$ unless $(\co',\cl'),(\co,\cl)$ are in the same block. Again this was 
known earlier when $m=1$ or $m=\iy$. In part III we discuss some applications of our methods to the study
of half-degenerate double affine graded algebras. In part IV we focus on one particular block, called the
principal block. (For $m=1$ this is the block of $(\co,\cl)$ where $\co=\{0\}$.) In this case we show that
the classical Springer correspondence is determined by purely combinatorial means in terms of the root
system.

\head\cite{241} The canonical basis of the quantum adjoint representation, 2017\endhead
Let $\fg$ be a simple Lie algebra over $\CC$ 
with a fixed pinning and let $V$ be the adjoint representation of $\fg$.
Let $V_v$ be the deformation of $V$ to a representation of the quantized universal enveloping algebra $U_v$
 of $\fg$. Here $v$ is an indeterminate. In my 1990 paper \cite{90} I showed that (in type $ADE$)
$V_q$ admits a basis in which the standard generators $E_i,F_i$ of $U_q$ act by very simple formulas in which 
all coeffocients are in $\NN[v,v\i]$; in fact they are either $1$ or $0$ or $v+v\i$. (Then in \cite{91} I 
showed that a similar result holds for the nonsimply laced types.) In this paper I show that this basis of 
$V_v$ is the canonical basis (which I defined in \cite{92} and which was not yet defined at the 
time I wrote this paper.) This gives a construction
of simple Lie algebras which is simpler than the standard one (which involves complicated signs). It also
gives a simple construction of Chevalley groups. 

\head \cite{244} On the generalized Springer correspondence, 2016\endhead
Let $G$ be a connected almost simple simply connected algebraic group of type $E_6$ over an algebraically 
closed field of characteristic $\ne3$. In a 1985 paper, Spaltenstein has determined the generalized
Springer correspondence for $G$ with one indeterminacy; there was a block governed by a Weyl group $W$ of
type $G_2$ for which one could not tell how to match the two irreducible representations of $W$ of 
dimension $2$ could be matched with corresponding local systems on a unipotent class. This was the only gap
in the determination of the generalized Springer correspondence if one restricts to characteristic $\ne2$.
(Another such gap remains for type $E_8$ in characteristic $2$.) In my 1986 paper on character sheaves I
proposed a remedy for this gap (for $E_6$), but later it turned out that there was an error in my calculation.
In this paper I redo the calculation and I find a matching which is the opposite of the one claimed in
my 1986 paper. The arguments involve a detailed study of various irreducible components of Springer 
fibres in type $E_6$. This 
completes the determination of the generalized Springer correspondence in characteristic $\ne2$.

\head\cite{246} (with G. Williamson) Billiards and tilting characters of $SL_3$, 2018\endhead
Let $G$ be a connected, simply connected reductive group over an algebraically closed field 
$\kk$ of characteristic $p>0$. The isomorphism classes of indecomposable tilting rational
representations of $G$ are indexed by the dominant weights of $G$. It is of considerable
interest to compute the character of these modules. This is known for $G=SL_2$ but not
known already for $G=SL_3$ even when $p$ is large. In this case only sporadic examples
were found by Andersen and his students around 2000 from which no pattern has emerged.
For any given dominant weight and given prime $p$, the desired character can in principle be 
calculated by a computer, but it is not clear that this fits any pattern. In this paper we have 
analyzed the data provided by the computer calculation in the case of $SL_3$ with $p=3,5,7$
and we found that these data obey some remarkable regularity which allowed us to make a conjectural 
statement for any $p\ge3$ and any dominant weights which in some sense are less than $p^3$. 
The conjectural statement uses some unexpected dynamical systems one of which involves the
 movement of a ball on a billiard table shaped as an equilateral triangle and the other is similar 
to a brownian motion on an infinite one-dimensional simplicial complex. A consequence of the conjecture 
is that the dimension of indecomposable tilting modules grows exponentially in the weight. Another 
consequence is a prediction of the pattern of decomposition of reduction modulo $p$ of the symmetric 
group applied to representations indexed by partitions with at most three parts.
It is likely that the complex behaviour observed in this paper for $SL_3$ reflects an interaction
between the subregular unipotent class $\co$ in $SL_3$ and the class of $1$ (the closure of $\co$ has a
rather serious singularity at $1$).

\head\cite{247} Conjugacy classes in reductive groups and two-sided cells, 2018\endhead
Let $G$ be a connected semisimple group over $\CC$. In my 1989 paper \cite{86} I proved that the
unipotent conjugacy classes of $G$ are in bijection with the two-sided cells of the affine Hecke algebra
attached to the Langlands dual of $G$. In this paper I prove an extension of this result in which I 
replace the set of unipotent classes of $G$ by the set of all conjugacy classes of $G$ of elements with 
semisimple part of finite order and I replace the affine Hecke algebra by an extended affine Hecke algebra for 
which the notion of two-sided cell can be still defined. The extended affine Hecke algebra has the
same relation to the unextended one as the extended Hecke algebra in the comments to \cite{234},\cite{245} 
to the
ordinary Hecke algebra.

\head\cite{248} Comments on my papers, 2017\endhead
This is a continuing project in which I write comments on some of my papers, 
point out connections with subsequent developments and also add some new points of view.

\head\cite{249} Lifting involutions in a Weyl group to the torus normalizer, 2018\endhead
Let $G$ be a connected reductive group with a given split rational structure over a finite field $F_q$
with Frobenius map $F:G@>>>G$ and with a fixed pinning defined over $F_q$.
Let $T$ be a maximal torus of $G$ and let $B$ be a Borel subgroup containing $T$. We assume that
$T$ and $B$ are defined over $F_q$. Let $NT$ be the normalizer of $T$ in $G$ and let $W=NT/T$ be the
Weyl group. Note that $F:NT@>>>NT$ induces the identity map on $W$. Let $\r:NT@>>>W$ be the obvious map. 
Let $I$ be the set of involutions in $W$. In this paper we define a canonical subset $\ti I$ of $N$ such that
$\r$ defines a bijection $\ti I@>>>I$ and such that $F(w)=w\i$ for any $w\in\ti I$. The definition of $\ti I$ 
is surprisingly complicated. It is based in part on
Kostant's cascades (which I rediscovered while working on this paper without knowing about them) and its 
proof in the exceptional cases relies in part on a computer calculation. The motivation for this paper 
came from the needs of the paper \cite{250} where the results of this paper are used.

\head\cite{250} Hecke modules based on involutions in extended Weyl groups, 2018\endhead
Let $G$ be a connected reductive group over an algebraically closed field of characteristic $p>0$. 
Let $T$ be a maximal torus of $G$ and let $W$ be the Weyl group. In my 2012 paper \cite{208}, I and 
Vogan have introduced a module for the Hecke algebra of $W$ with basis indexed by the 
involutions in $W$. This is extended in the present paper as follows. We consider the dual torus $T^*$
and we consider the group $\ti W$ consisting of all pairs $(w,t)\in W\T T^*$ where $t$ has finite order (a
semidirect product). We say that $\ti W$ is the extended Weyl group. Let $\ti I$ be the set of involutions in
$\ti W$. In this paper we show that the vector space with basis $\ti I$ is naturally a module over the
``extended Hecke algebra'' (see the comments to \cite{234},\cite{245}). 
The proof uses the connection of the extended Hecke
algebra with the Hecke algebra considered by Yokonuma in the 1960's. The proof also uses the canonical
lifting of involutions in $W$ to the normalizer of $T$ which I studied in \cite{249}.

\head\cite{251} Discretization of Springer fibres, 2017\endhead
Let $G$ be a connected almost simple simply connected algebraic group of adjoint type
over $\CC$ with Lie algebra $\fg$, let $e\in\fg$ be a nilpotent element 
and let $\cb_e$ be the variety of Borel subalgebras of $\fg$ that contain $e$. It is known that $\cb_e$ has 
a natural $\CC^*$-action coming from the Morozov-Jacobson theorem. We consider the equivariant $K$-group 
$K_{\CC^*}(\cb_e)$ of $\cb_e$. This is naturally a module over $\ZZ[v,v\i]$ where $v$ is an indeterminate.
From the analysis of $\cb_e$ given in my 1988 paper \cite{78} with De Concini and Procesi, it follows that 
$K_{\CC^*}(\cb_e)$ is in fact a free $\ZZ[v,v\i]$-module of finite rank. In my 1999 paper \cite{143} 
I defined a subset 
$B_e^\pm$ of $K_{\CC^*}(\cb_e)$ stable under multiplication by $-1$ and I conjectured that this is a signed 
basis of the $\ZZ[v,v\i]$-module $K_{\CC^*}(\cb_e)$. (This conjecture has been later proved by
Bezrukavnikov and Mirkovic.) The set $B_e^\pm/\pm1$ has a natural action of $A(e)$, the centralizer of $e$
in $G$ which factors through a finite quotient. Now $B_e^\pm/\pm1$ can be viewed as a discrete analogue of the 
Springer fibre. For example the cardinal of $B_e^\pm/\pm1$ is equal to the sum of Betti numbers of $\cb_e$.
In this paper we propose a conjectural description of the $A(e)$-set $B_e^\pm/\pm1$. Namely we give
a conjectural description of the isotropy groups of the various points of
$B_e^\pm/\pm1$ and of the number of times each isotropy group occurs.
This is useful since the $A(e)$-set $B_e^\pm/\pm1$ is conjecturally the main ingredient in a description
(as an equivariant $K$-group) 
of the $J$-ring of an (extended) affine Weyl group which I proposed in my 1989 paper \cite{86}.

\head \cite{252} A new basis for the representation ring of a Weyl group; \cite{255} II, 2018\endhead
Let $W$ be a Weyl group.  In my 1979 paper \cite{36} 
I introduced a class $\cs_W$ of irreducible representations of $W$ (later called
special representations). In my 1982 paper \cite{50} I introduced a class $C_W$ of not necessarily irreducible
representations of $W$ (later called constructible). In my 1986 paper \cite{70} 
I showed that $C_W$ consists exactly of the
representations of $W$ carried by the various left cells of $W$. In this paper I introduce a class $\fC_W$
of representations of $W$ which interpolates between $\cs_W$ and $C_W$ (it contains both). One of the
main results is that $\fC_W$ provides a $\ZZ$-basis for the Grothendieck group of representations of $W$.
Moreover I show that this new basis and the standard basis are in natural bijection and that 
the transition matrix between the new basis and the standard basis is triangular with $1$ on diagonal.
I also show that any representation $\ce$ in the new basis can be described by a pair $\ch,\ch'$ of subgroups
of a certain finite group associated to the two-sided cell containing the irreducible components of $\ce$.
This gives a  new way to index the irreducible representations in a fixed two-sided cell. In this paper
I conjecture that the representations in $\fC_W$ have a positivity property analogous to that
characterizing the special representations (see \cite{239}).

\head \cite{253} Positive conjugacy classes in Weyl groups, 2018\endhead
Let $W$ be a Weyl group. Let $H$ be the corresponding Hecke algebra over $\CC(q)$ (with $q$ an indeterminate)
and let $\{T_w;w\in W\}$ be the standard basis of $H$. For any irreducible representation $E$ of $W$
let $E_q$ be the corresponding $H$-module. Let $C$ be a conjugacy class of $W$ and let $C_{min}$ be the
set of minimal length elements of $C$. For any $w\in C_{min}$ we can form the polynomial
$N_w=\sum_E\tr(T_w,E_q)^2$ where $E$ runs over the irreducible representsations of $W$. From Geck-Pfeiffer one
can deduce that $N_w$ depends only on $C$, not on $w$. We say that $C$ is positive if $N_w\in\NN[q]$ for 
some/any $w\in C_{min}$. The main observation of this paper is that many elliptic conjugacy classes in $W$ are
positive. For example all elliptic conjugacy classes which are regular in the sense of Springer are positive.
We give also a number of examples of positive conjugacy classes which are not regular. In exceptional types
we give a complete list, but for classical types we do not know a complete list.


\widestnumber\key {ABC}
\Refs
\ref\key{1} \by G. Lusztig\paper Model de geometrie afina plana peste un corp finit\jour Studii Cerc. Mat.\vol17\yr1965)\pages1337-1340\endref

\ref\key{2} \by G.Lusztig \paper 
Constructia fibrarilor universale peste poliedre arbitrare\jour  Studii Cerc. Mat.\vol18\yr1965\pages1215-1219\endref

\ref\key{3}\by G.Lusztig and H.Moscovici\paper Demonstration du th\'eoreme sur la suite spectrale d'un fibr\'e au sens de Kan\jour Proc. Camb.   
Phil. Soc.\vol64\yr1968 \pages293-297\endref  

\ref\key{4} \by G.Lusztig \paper  Sur les complexes elliptiques fibr\'es\jour C.R. Acad. Sci. Paris(A)\vol266\yr1968
\pages 914-917\endref

\ref\key{5} \by G.Lusztig \paper  Sur les actions libres des groupes finis\jour Bull. Acad. Polon. Sci.\vol16
\yr1968\pages461-463\endref

\ref\key {6} \by G.Lusztig \paper  Coomologia complexelor eliptice\jour Studii Cerc. Mat.\vol21\yr1969\pages 38-83
\endref

\ref\key {7} \by G.Lusztig \paper  A property of certain non-degenerate holomorphic vector fields\jour
 An. Univ. Timisoara\vol 7\yr1969\pages 73-76\endref

\ref\key {8} \by G.Lusztig, J.Milnor and F.P.Peterson\paper Semicharacteristics and cobordism\jour Topology\vol 8
\yr1969\pages 357-359\endref

\ref\key {9} \by G.Lusztig \paper  Remarks on the holomorphic Lefschetz formula\inbook Analyse globale
\publ Presses de l'Univ.de Montr\'eal\yr 1971\pages193-204\endref

\ref\key {10} \by G.Lusztig and J.Dupont\paper On manifolds satisfying $w_1^2 =0$\jour Topology\vol 10\yr1971
\pages 81-92\endref

\ref\key {11} \by G.Lusztig \paper  Novikov's higher signature and families of elliptic operators\jour J. Diff. Geom.
\vol7\yr1972 \pages229-256\endref

\ref\key {12} \by G.Lusztig \paper  On the discrete series  representations of the general 
linear groups over a finite field\jour Bull. Amer. Math. Soc.\vol79\yr1973\pages 550-554\endref

\ref\key {13} \by G.Lusztig \book The discrete  series of $GL_n$ over a finite field\bookinfo Ann. Math. Studies 81
\publ Princeton U.Press \yr1974\endref

\ref\key {14} \by G.Lusztig \paper  Introduction to elliptic operators\inbook Global Analysis and applications
\bookinfo Internat.Atomic Energy Agency, Vienna \yr1974\pages 187-193\endref

\ref\key {15} \by G.Lusztig and R.W.Carter \paper  On the modular representations of the general linear and symmetric groups\jour Math. Z.\vol136\yr1974\pages193-242\endref

\ref\key {16} \by G.Lusztig and R.W.Carter\paper Modular representations of the general linear and symmetric groups
\inbook Proc.2nd Int.Conf. Th.Groups 1973, LNM 372\publ Springer Verlag \yr1974\pages218-220\endref

\ref\key {17} \by G.Lusztig \paper On the discrete series representations of the classical groups over a finite field
 \inbook Proc.Int.Congr.Math.,Vancouver 1974\pages 465-470\endref

\ref\key {18} \by G.Lusztig \paper Sur la conjecture de Macdonald\jour C.R. Acad. Sci. Paris(A)\vol280\yr1975
\pages 371-320\endref

\ref\key {19} \by G.Lusztig \paper  A note on counting nilpotent matrices of fixed rank\jour Bull. Lond. Math. Soc.
\vol8\yr1976\pages 77-80\endref

\ref\key {20} \by G.Lusztig \paper  Divisibility of projective modules of finite Chevalley groups by the 
Steinberg module\jour Bull. Lond. Math. Soc.\vol8\yr1976\pages 130-134\endref

\ref\key {21} \by G.Lusztig and R.W.Carter \paper Modular representations of finite groups of Lie type
\jour Proc. Lond. Math. Soc.\vol32\yr1976\pages 347-384\endref

\ref\key {22} \by G.Lusztig and P.Deligne\paper Representations of reductive groups over finite fields\jour
 Ann. Math.\vol103\yr1976\pages 103-161\endref

\ref\key {23} \by G.Lusztig \paper  On the Green polynomials of classical groups\jour Proc. Lond. Math. Soc.\vol33
\yr1976\pages443-475\endref

\ref\key {24} \by G.Lusztig \paper  Coxeter orbits and eigenspaces of Frobenius\jour Inv. Math.\vol28\yr1976\pages 
101-159\endref

\ref\key {25} \by G.Lusztig \paper  On the finiteness of the number of unipotent classes\jour Inv. Math.\vol 34
\yr1976\pages 201-213\endref

\ref\key{26} \by G.Lusztig, J.A.Green and G.I.Lehrer\paper On the degrees of certain group characters\jour Quart.
J. Math.\vol 27\yr1976\pages1-4\endref

\ref\key {27} \by G.Lusztig and B.Srinivasan\paper The characters of the finite unitary groups\jour J. Alg.\vol49
\yr1977\pages 167-171\endref

\ref\key {28} \by G.Lusztig \paper  Classification des repr\'esentations irr\'eductibles des groupes classiques 
finis\jour C.R. Acad. Sci. Paris(A)\vol284\yr1977\pages 473-476\endref

\ref\key {29} \by G.Lusztig \paper  Irreducible representations of finite classical groups\jour Inv. Math.\vol43
\yr1977\pages125-175\endref

\ref\key {30} \by G.Lusztig \book  Representations of finite Chevalley groups \bookinfo Regional Conf. Series in 
Math. 39\publ Amer. Math. Soc.\yr 1978\endref

\ref\key {31} \by W.M.Beynon and G.Lusztig \paper Some numerical results on the characters of exceptional Weyl 
groups\jour Math. Proc. Camb. Phil. Soc.\vol84\yr1978\pages 417-426\endref

\ref\key{32} \by G.Lusztig \paper  Some remarks on the supercuspidal representations of p-adic semisimple \lb groups
\inbook Proc. Symp. Pure Math.33(1)\publ Amer. Math. Soc.\yr 1979\pages 171-175\endref 

\ref\key {33} \by G.Lusztig \paper  On the reflection representation of a finite Chevalley group, \inbook
Representation theory of Lie groups \bookinfo LMS 
Lect.Notes Ser.34\publ Cambridge U.Press \yr1979\pages 325-337\endref

\ref\key {34} \by G.Lusztig \paper  Unipotent representations of a finite Chevalley group of type $E_8$\jour 
Quart. J. Math.\vol30\yr1979\pages 315-338\endref

\ref\key {35} \by G.Lusztig and N.Spaltenstein\paper Induced unipotent classes\jour J. Lond. Math. Soc.\vol19\yr1979
\pages 41-52\endref

\ref\key {36} \by G.Lusztig \paper  A class of irreducible representations of a Weyl group\jour Proc. Kon. Nederl.
Akad.(A)\vol82\yr1979\pages 323-335\endref

\ref\key {37} \by D.Kazhdan and  G.Lusztig \paper  Representations of Coxeter groups and Hecke algebras\jour 
Inv. Math.\vol53\yr1979\pages 165-184\endref

\ref\key {38} \by D.Kazhdan and  G.Lusztig \paper  A topological approach to Springer's representations\jour 
Adv. Math.\vol38\yr1980\pages 222-228\endref

\ref\key {39} \by D.Kazhdan and G.Lusztig \paper  Schubert varieties and Poincar\'e duality\inbook Proc. Symp.
Pure Math.36\publ Amer. Math. Soc.\yr 1980\pages185-203\endref

\ref\key {40} \by G.Lusztig \paper  Some problems in the representation theory of finite Chevalley groups
\inbook Proc. Symp. Pure Math.37\publ Amer. Math. Soc.\yr1980\pages 313-317\endref 

\ref\key {41} \by G.Lusztig \paper  Hecke algebras and Jantzen's generic decomposition patterns\jour Adv.Math.
\vol37\yr1980\pages 121-164\endref

\ref\key{42} \by G.Lusztig \paper On the unipotent characters of the exceptional groups over finite fields\jour 
Inv. Math.\vol60\yr1980\pages173-192\endref

\ref\key {43} \by G.Lusztig \paper  On a theorem of Benson and Curtis\jour J. Alg.\vol71\yr1981\pages 490-498\endref

\ref\key {44} \by G.Lusztig \paper  Green polynomials and singularities of unipotent classes \jour Adv. Math.\vol42
\yr1981\pages 169-178\endref

\ref\key {45} \by G.Lusztig \paper  Unipotent characters of the symplectic and odd orthogonal groups over 
a finite field\jour Inv. Math.\vol64\yr1981\pages263-296\endref

\ref\key {46} \by G.Lusztig \paper  Unipotent characters of the even orthogonal groups over a finite field
\jour Trans. Amer. Math. Soc.\vol272\yr1982 \pages733-751\endref

\ref\key {47} \by P.Deligne and G.Lusztig \paper  Duality for representations of a reductive group over a 
finite field\jour J. Alg.\vol74\yr1982\pages284-291\endref

\ref\key {48} \by  D.Alvis and G.Lusztig \paper  On Springer's 
correspondence for simple groups of type $E_n(n=6,7,8)$\jour Math. Proc. Camb. Phil. Soc.
\vol92\yr1982\pages 65-72\endref

\ref\key {49} \by  D.Alvis and G.Lusztig \paper  The representations and generic degrees of the Hecke algebras of type $H_4$, J. reine und angew. math.\vol336\yr1982\pages 201-212\moreref Erratum\vol 449\yr1994 \pages217-281\endref

\ref\key {50} \by G.Lusztig \paper  A class of irreducible representations of a Weyl group II\jour
 Proc. Kon. Nederl. Akad.(A)\vol85\yr1982\pages 219-226\endref

\ref\key {51} \by G.Lusztig and D.Vogan\paper Singularities of closures of K-orbits on a flag manifold\jour Inv.
Math.\vol71\yr1983\pages 365-379\endref

\ref\key {52} \by P.Deligne and G.Lusztig \paper Duality for representations of a reductive group over a finite 
field II\jour J. Alg.\vol81\yr1983\pages 540-549\endref

\ref\key {53} \by G.Lusztig \paper  Singularities, character formulas and a $q$-analog of weight multiplicities
\lb\jour Ast\'erisque\vol101-102\yr1983\pages 208-229\endref

\ref\key {54} \by G.Lusztig \paper  Some examples of square integrable representations of semisimple p-adic \lb groups
\jour Trans. Amer. Math. Soc.\vol227\yr1983\pages623-653\endref

\ref\key {55} \by G.Lusztig \paper  Left cells in Weyl groups\inbook  Lie groups representations\bookinfo
 LNM 1024\publ Springer Verlag\yr 1983\pages 99-111\endref

\ref\key {56} \by G.Lusztig \paper Open problems in algebraic groups\inbook  Proc.12th Int.Symp., Taniguchi 
Foundation, Katata\yr 1983\pages 14-14\endref

\ref\key {57} \by G.Lusztig \book Characters of 
reductive groups over a finite field\bookinfo  Ann.Math.Studies 107\publ Princeton U.Press \yr1984\endref

\ref\key {58} \by G.Lusztig \paper  Characters of reductive groups over finite fields\inbook Proc.Int.Congr.Math. 
Warsaw 1983\publ North Holland \yr1984
\pages877-880\endref

\ref\key {59} \by G.Lusztig \paper  Intersection cohomology complexes on a reductive group\jour Inv. Math.\vol75
\yr1984\pages205-272\endref

\ref\key {60} \by G.Lusztig \paper  Cells in affine Weyl groups\inbook Algebraic groups and related topics
\bookinfo Adv. Stud. Pure Math. 6\publ  North-Holland and Kinokuniya\yr 1985\pages 255-287\endref

\ref\key {61} \by G.Lusztig and N.Spaltenstein\paper On the generalized Springer correspondence for classical groups
\inbook Algebraic groups and related topics
\bookinfo Adv. Stud. Pure Math. 6\publ  North-Holland and Kinokuniya\yr 1985\pages 289-316\endref

\ref\key {62} \by G.Lusztig \paper  The two sided cells of the affine Weyl group of type A\inbook
Infinite dimensional groups with applications\bookinfo MSRI Publ.4\publ Springer Verlag\yr 1985\pages 275-283\endref

\ref\key {63} \by G.Lusztig \paper  Character sheaves I\jour Adv. Math.\vol56\yr1985\pages 193-237\endref

\ref\key {64} \by G.Lusztig \paper  Character sheaves II\jour Adv.Math.\vol57\yr1985\pages 226-265\endref

\ref\key {65} \by G.Lusztig \paper  Character sheaves III\jour Adv.Math.\vol57\yr1985\pages 266-315\endref

\ref\key {66} \by G.Lusztig \paper Equivariant K-theory and representations of Hecke algebras\jour 
Proc. Amer. Math. Soc.\vol94\yr1985\pages 337-342\endref

\ref\key {67} \by D.Kazhdan  and G.Lusztig\paper Equivariant K-theory and representations of Hecke algebras II
\jour Inv. Math.\vol80\yr1985\pages 209-231\endref

\ref\key {68} \by G.Lusztig \paper  Character sheaves IV\jour Adv. Math.\vol59\yr1986\pages 1-63\endref

\ref\key {69} \by G.Lusztig \paper Character sheaves V\jour Adv. Math.\vol61\yr1986\pages103-155\endref

\ref\key {70} \by G.Lusztig \paper  Sur les cellules gauches des groupes de Weyl\jour C.R. Acad. Sci. Paris(A)
\vol302\yr1986\pages5-8\endref

\ref\key {71} \by G.Lusztig \paper  On the character values of finite Chevalley groups at unipotent elements\jour
 J. Alg.\vol104\yr1986\pages 146-194\endref

\ref\key {72} \by  D.Kazhdan and  G.Lusztig \paper  Proof of the Deligne-Langlands conjecture for Hecke algebras
\jour Inv. Math.\vol87\yr1987\pages 153-215\endref

\ref\key {73} \by G.Lusztig \paper  Cells in affine Weyl groups II\jour J. Alg.\vol109\yr1987\pages 536-548\endref

\ref\key {74} \by G.Lusztig \paper  Fourier transforms on a semisimple Lie algebra over $F_q$, \inbook
Algebraic Groups Utrecht 1986\bookinfo LNM 1271\publ Springer Verlag\yr 1987\pages177-188\endref

\ref\key {75} \by G.Lusztig \paper  Cells in affine Weyl groups III \jour J. Fac. Sci. Tokyo U.(IA)\vol34\yr1987
\pages 223-243\endref

\ref\key {76} \by G.Lusztig \paper  Introduction to character sheaves\inbook Proc. Symp. Pure  Math. 47(1)
\publ Amer. Math. Soc.\yr 1987\pages 165-180\endref

\ref\key{77} \by G.Lusztig \paper Leading coefficients of character values of Hecke algebras\inbook
 Proc. Symp. Pure Math. 47(2)\publ Amer. Math. Soc.\yr 1987\pages235-262\endref

\ref\key {78} \by C.De Concini, G.Lusztig and C.Procesi\paper Homology of the zero set of a nilpotent 
vector field on a flag manifold \jour J. Amer. Math. Soc.\vol1\yr1988\pages15-34\endref

\ref\key {79} \by G.Lusztig \paper  Quantum deformations of certain simple modules over enveloping algebras
\jour Adv.Math.\vol70\yr1988\pages 237-249\endref

\ref\key{80}\by D.Kazhdan and G.Lusztig \paper  Fixed point varieties on affine flag  manifolds\jour Isr. J. Math.
\vol62\yr1988 \pages 129-168\endref

\ref\key {81} \by G.Lusztig \paper Cuspidal local systems and graded Hecke algebras I\jour Publ. Math. I.H.E.S.
\vol67\yr1988\pages 145-202\endref

\ref\key {82} \by G.Lusztig and  N.Xi\paper Canonical left cells in affine Weyl groups\jour Adv.Math.\vol72\yr1988
\pages284-288\endref

\ref\key {83} \by G.Lusztig \paper On representations of reductive groups with disconnected center\jour 
Ast\'erisque\vol168\yr1988\pages 157-166\endref

\ref\key {84} \by G.Lusztig \paper  Modular representations and quantum groups\jour Contemp. Math.\vol82\yr1989
\pages 59-77\endref

\ref\key {85} \by G.Lusztig \paper  Affine Hecke algebras and their graded version\jour J. Amer. Math. Soc.\vol2
\yr1989\pages 599-635\endref

\ref\key {86} \by G.Lusztig \paper  Cells in affine Weyl groups IV\jour J. Fac. Sci. Tokyo U.(IA)\vol36\yr1989\pages
 297-328\endref

\ref\key {87} \by G.Lusztig \paper  Representations of affine Hecke algebras\jour Ast\'erisque\vol 171-172\yr1989
\pages 73-84\endref

\ref\key {88} \by G.Lusztig \paper  On quantum groups\jour J. Alg.\vol131\yr1990\pages 466-475\endref

\ref\key {89} \by G.Lusztig \paper  Green functions and character sheaves\jour Ann. Math.\vol131\yr1990
\pages 355-408\endref

\ref\key {90} \by G.Lusztig \paper  Finite dimensional Hopf algebras arising from quantized universal enveloping 
algebras\jour J. Amer. Math. Soc.\vol3\yr1990\pages 257-296\endref

\ref\key {91} \by G.Lusztig \paper  Quantum groups at roots of 1\jour Geom.Ded.\vol35\yr1990\pages 89-114\endref

\ref\key {92} \by G.Lusztig \paper  Canonical bases arising from quantized enveloping algebras\jour J. Amer. Math. 
Soc.\vol3\yr1990\pages 447-498\endref

\ref\key {93} \by A.A.Beilinson, G.Lusztig and R.MacPherson\paper A geometric setting for the quantum deformation of 
$GL_n$ \jour Duke Math. J.\vol61\yr1990\pages 655-677\endref

\ref\key {94} \by G.Lusztig \paper  Symmetric spaces over a finite field\inbook The Grothendieck Festschrift III
\bookinfo Progr. in Math. 88\publ Birkh\"auser Boston \yr1990\pages 57-81\endref

\ref\key{95} \by G.Lusztig \paper  Canonical bases arising from quantized enveloping algebras II\jour
Progr.of Theor. Phys. Suppl.\vol102\yr1990\pages 175-201\endref

\ref\key {96}\by  D.Kazhdan and G.Lusztig \paper  Affine Lie algebras and quantum groups\jour Int. Math. Res. Notices 
\yr1991\pages 21-29\endref

\ref\key {97} \by G.Lusztig \paper  Quivers, perverse sheaves and enveloping algebras
 J. Amer. Math. Soc.\vol4\yr1991\pages 365-421\endref

\ref\key {98} \by G.Lusztig and J.M.Smelt\paper Fixed point varieties in the space of lattices
 Bull. Lond. Math. Soc.\vol23\yr1991\pages 213-218\endref

\ref\key {99} \by G.Lusztig \paper  Intersection cohomology methods in representation theory\inbook Proc. Int.
Congr. Math. Kyoto 1990\publ Springer Verlag\yr1991\pages 155-174\endref

\ref\key {100} \by G.Lusztig \paper  A unipotent support for irreducible representations\jour Adv. Math.\vol94\yr1992
\pages 139-179\endref

\ref\key {101} \by G.Lusztig \paper  Canonical bases in tensor products 
\jour Proc. Nat. Acad. Sci.\vol89\yr1992\pages 8177-8179\endref
     
\ref\key{102} \by G.Lusztig \paper  Remarks on computing irreducible characters\jour
 J. Amer. Math. Soc.\vol5\yr1992\pages 971-986\endref

\ref\key {103} \by G.Lusztig \paper  Introduction to quantized enveloping algebras \inbook Progr.in Math.105
\publ Birkh\"auser Boston\yr 1992\pages 49-65\endref

\ref\key {104} \by G.Lusztig \paper  Affine quivers and canonical bases\jour Publ. Math. I.H.E.S.\vol76\yr1992\pages 
111-163\endref

\ref\key {105} \by G.Lusztig and J.Tits\paper The inverse of a Cartan matrix\jour An.Univ.Timisoara\vol 30\yr1992
\pages17-23\endref

\ref\key {106} \by I.Grojnowski and  G.Lusztig \paper  On bases of irreducible representations of quantum $GL_n$
\inbook Kazhdan-Lusztig theory and
related topics\bookinfo Contemp.Math.139\yr1992\pages 167-174\endref

\ref\key {107} \by G.Lusztig \book Introduction to quantum groups\bookinfo Progr.in Math.110\publ Birkh\"auser 
Boston\yr 1993\endref

\ref\key {108} \by D.Kazhdan and G.Lusztig\paper Tensor structures arising from affine Lie algebras I\jour
 J. Amer. Math. Soc.\vol6\yr1993\pages 905-947\endref

\ref\key {109}\by D.Kazhdan and  G.Lusztig \paper  Tensor structures arising from affine Lie algebras II\jour 
J. Amer. Math. Soc.\vol6\yr1993\pages 949-1011\endref

\ref\key {110} \by G.Lusztig \paper  Coxeter groups and unipotent representations \jour Ast\'erisque\vol 212\yr1993
\pages 191-203\endref

\ref\key {111}\by I.Grojnowski and  G.Lusztig \paper  A comparison of bases of quantized enveloping algebras\inbook
Linear algebraic groups and their representations\bookinfo Contemp.Math.153\yr1993\pages 11-19\endref

\ref\key {112} \by G.Lusztig \paper  Tight monomials in quantized enveloping algebras\inbook 
Quantum deformations of algebras and their representations\bookinfo  ed. A.Joseph et al., Isr. Math. Conf. Proc. 7
\publ Amer. Math. Soc.\yr 1993\pages 117-132\endref

\ref\key {113} \by G.Lusztig \paper  Exotic Fourier transform\jour Duke Math.J.\vol73\yr1994\pages 227-241\endref

\ref\key {114} \by G.Lusztig \paper  Vanishing properties of cuspidal local systems\jour Proc. Nat. Acad. Sci.\vol91
\yr1994\pages 1438-1439\endref

\ref\key {115}\by D.Kazhdan and G.Lusztig \paper Tensor structures arising from affine Lie algebras III\jour 
J. Amer. Math. Soc.\vol7\yr1994\pages 335-381\endref

\ref\key {116}\by  D.Kazhdan and G.Lusztig \paper Tensor structures arising from affine Lie algebras IV\jour
 J. Amer. Math. Soc.\vol7\yr1994\pages 383-453\endref

\ref\key {117} \by G.Lusztig \paper  Monodromic systems on affine flag manifolds\jour Proc. Roy. Soc. Lond.(A) 
\vol445\yr1994\pages231-246\moreref Errata\vol 450\yr1995\pages731-732\endref

\ref\key {118} \by G.Lusztig \paper  Problems on canonical bases \inbook Algebraic groups and their 
generalizations: quantum and infinite dimensionalmethods\bookinfo Proc. Symp. Pure Math. 56(2)\publ Amer. Math. Soc. 
\yr1994\pages 169-176\endref

\ref\key {119} \by G.Lusztig \paper  Total positivity in reductive groups\inbook Lie theory and geometry\bookinfo
 Progr.in Math. 123\publ Birkh\"auser Boston \yr1994\pages 531-568\endref

\ref\key {120} \by G.Lusztig \paper  Study of perverse sheaves arising from graded Lie algebras\jour Adv.Math.\vol112
\yr1995\pages 147-217\endref

\ref\key {121} \by G.Lusztig \paper  Cuspidal local systems and graded Hecke algebras II\inbook Representations of 
groups\bookinfo ed. B.Allison et al., Canad. Math. Soc. Conf. Proc.16\publ Amer. Math. Soc.\yr 1995\pages217-275\endref

\ref\key {122} \by G.Lusztig \paper  Quantum groups at $v=\iy$\inbook Functional analysis on the eve of the 21st 
century, vol.I\bookinfo Progr.in Math. 131, Birkh\"auser Boston\yr 1995\pages 199-221\endref

\ref\key {123} \by G.Lusztig \paper  Classification of unipotent representations of simple $p$-adic groups\jour
 Int. Math. Res. Notices\yr1995\pages 517-589\endref

\ref\key {124} \by G.Lusztig \paper  An algebraic-geometric parametrization of the canonical basis\jour Adv. Math.
\vol120\yr1996\pages 173-190\endref

\ref\key {125} \by G.Lusztig \paper Affine Weyl groups and conjugacy classes in Weyl groups\jour Transform. Groups 
\yr1996\pages83-97\endref

\ref\key {126} \by G.Lusztig \paper  Braid group actions and canonical bases\jour Adv. Math.\vol122\yr1996\pages 237-261\endref

\ref\key {127} \by G.Lusztig \paper  Non local finiteness of a $W$-graph\jour Represent.Th\vol1\yr1997\pages25-30\endref

\ref\key{128} \by G.Lusztig \paper  Cohomology of classifying spaces and hermitian representations\jour Represent.Th.
\vol1\yr1997\pages 31-36\endref

\ref\key {129}\by C.K.Fan and  G.Lusztig \paper Factorization of certain exponentials in Lie groups\inbook
Algebraic groups and Lie groups\bookinfo ed. G.I.Lehrer\publ Cambridge U.Press\yr 1997\pages 215-218\endref

\ref\key {130} \by G.Lusztig \paper  Total positivity and canonical bases\inbook Algebraic groups and Lie groups 
\bookinfo ed. G.I.Lehrer\publ Cambridge U.Press\yr1997\pages 281-295\endref

\ref\key {131} \by G.Lusztig \paper  Notes on unipotent classes\jour Asian J.Math.\vol1\yr1997\pages 194-207\endref

\ref\key {132} \by G.Lusztig \paper  Cells in affine Weyl groups and tensor categories\jour Adv. Math.\vol129
\yr1997\pages 85-98\endref

\ref\key {133} \by G.Lusztig \paper  Periodic $W$-graphs\jour Represent.Th.\vol1\yr1997\pages 207-279\endref
   
\ref\key {134} \by G.Lusztig \paper  A comparison of two graphs\jour Int. Math. Res. Notices\yr1997\pages 639-640\endref

\ref\key {135} \by G.Lusztig \paper  Constructible functions on the Steinberg variety\jour Adv. Math.\vol130\yr1997
\pages287-310\endref

\ref\key {136} \by G.Lusztig \paper  Total positivity in partial flag manifolds\jour Represent.Th.\vol2\yr1998
\pages 70-78\endref

\ref\key {137} \by G.Lusztig \paper  Introduction to total positivity\inbook Positivity in Lie theory: open problems
\bookinfo  ed. J.Hilgert et al.\publ de Gruyter\yr1998\pages 133-145\endref

\ref\key {138} \by G.Lusztig \paper  On quiver varieties\jour Adv.Math.\vol136\yr1998\pages 141-182\endref

\ref\key {139} \by G.Lusztig \paper  Canonical bases and Hall algebras \inbook Representation Theories and 
Algebraic Geometry\bookinfo ed. A.Broer et al.\publ Kluwer Acad.Publ.\yr 1998\pages 365-399\endref

\ref\key {140} \by G.Lusztig \paper  Bases in equivariant $K$-theory\jour Represent.Th.\vol2\yr1998\pages 298-369\endref

\ref\key {141} \by G.Lusztig \paper  Homology bases arising from reductive groups over a finite field\inbook
Algebraic groups and their 
representations\bookinfo ed. R.W.Carter et al.\publ Kluwer Acad. Publ.\yr 1998\pages 53-72\endref

\ref\key {142} \by G.Lusztig \paper  Aperiodicity in quantum affine $\fg\fl_n$\jour Asian J. Math.\vol3\yr1999
\pages 147-178\endref

\ref\key {143} \by G.Lusztig \paper  Bases in equivariant $K$-theory II\jour Represent. Th.\vol3\yr1999\pages 281-353\endref

\ref\key {144} \by G.Lusztig \paper  A survey of group representations\jour Nieuw Archief voor Wiskunde \vol17\yr1999
\pages 483-489\endref
                                          
\ref\key {145} \by G.Lusztig \paper  Subregular nilpotent elements and bases in $K$-theory\jour Canad. J. Math.\vol51
\yr1999\pages1194-1225\endref

\ref\key {146} \by G.Lusztig \paper  Recollections about my teacher, Michael Atiyah\jour Asian J. Math.\vol3\yr1999
\pages iv-v\endref

\ref\key {147} \by G.Lusztig \paper  Semicanonical bases arising from enveloping algebras\jour Adv. Math.\vol151
\yr2000\pages129-139\endref

\ref\key {148} \by G.Lusztig \paper  Fermionic form and Betti numbers, arxiv:QA/0005010\endref

\ref\key {149} \by G.Lusztig \paper  Quiver varieties and Weyl group actions\jour Ann. Inst. Fourier\vol 50\yr2000
\pages 461-489\endref

\ref\key {150} \by G.Lusztig \paper  $G(F_q)$-invariants in irreducible $G(F_{q^2})$-modules\jour Represent. Th.\vol4
\yr2000 \pages446-465\endref

\ref\key {151} \by G.Lusztig \paper  Remarks on quiver varieties\jour Duke Math. J.\vol105\yr2000\pages 239-265\endref

\ref\key {152} \by G.Lusztig \paper  Transfer maps for quantum affine $\fs\fl_n$ \inbook 
Representations and quantizations\bookinfo ed. J.Wang et al.\publ China Higher Ed.Press and Springer Verlag\yr 2000
\endref

\ref\key {153} \by G.Lusztig \paper  Representation theory in characteristic $p$\inbook Taniguchi Conf. on Math.
Nara'98\bookinfo Adv. Stud. Pure Math. 31\publ Math. Soc. Japan \yr2001\pages 167-178\endref

\ref\key {154} \by G.Lusztig \paper  Cuspidal local systems and graded Hecke algebras III\jour Represent.Th.\vol6
\yr2002\pages 202-242\endref

\ref\key {155} \by G.Lusztig \paper  Classification of unipotent representations of simple $p$-adic groups II
\jour Represent.Th.\vol6\yr2002\pages243-289\endref

\ref\key {156} \by G.Lusztig \paper  Constructible functions on varieties attached to quivers\inbook Studies in 
memory of I. Schur\bookinfo Progr. in Math. 210\publ Birkh\"auser Boston\yr 2002\pages177-223\endref

\ref\key {157} \by G.Lusztig \paper  Rationality properties of unipotent representations\jour J. Alg.\vol258
\yr2002\pages 1-22\endref

\ref\key {158} \by G.Lusztig \paper  Notes on affine Hecke algebras\inbook Iwahori-Hecke algebras and their 
representation theory\bookinfo ed. M.W.Baldoni et al., LNM 1804\publ Springer Verlag\yr 2002\pages 71-103\endref

\ref\key {159} \by G.Lusztig \book Hecke algebras with unequal parameters\bookinfo  CRM Monograph Ser.18\publ
 Amer. Math. Soc. 2003,136p\finalinfo additional material in version 2 (2014), arxiv:math/0208154 \endref 

\ref\key {160} \by G.Lusztig \paper  Homomorphisms of the alternating group $A_5$ into reductive groups\jour J. Alg.
\vol260\yr2003\pages 298-322\endref

\ref\key {161} \by G.Lusztig \paper  Character sheaves on disconnected groups I \jour Represent. Th.\vol7\yr2003
\pages 374-403 \moreref Errata \vol8\yr2004\pages179-179\endref

\ref\key {162} \by G.Lusztig \paper  Representations of reductive groups over finite rings\jour Represent. Th.\vol8
\yr2004\pages 1-14\endref

\ref\key {163} \by G.Lusztig \paper  Character sheaves on disconnected groups II\jour Represent. Th.\vol8\yr2004
\pages 72-124\endref

\ref\key {164} \by G.Lusztig \paper  Character sheaves on disconnected groups III\jour Represent. Th.\vol8\yr2004
\pages 125-144
\endref

\ref\key {165} \by G.Lusztig \paper  Character sheaves on disconnected groups IV\jour Represent. Th.\vol8\yr2004
\pages 145-178
\endref

\ref\key {166} \by G.Lusztig \paper  Parabolic character sheaves I\jour Moscow Math.J.\vol4\yr2004\pages 153-179\endref

\ref\key {167} \by G.Lusztig \paper  An induction theorem for Springer's representations\inbook Adv.Stud.Pure Math.40
\publ Math. Soc. Japan, Kinokuniya\yr 2004\pages 253-259\endref

\ref\key {168} \by G.Lusztig \paper  Character sheaves on disconnected groups V\jour Represent. Th.
\vol8\yr2004\pages 346-376\endref

\ref\key {169} \by G.Lusztig \paper  Character sheaves on disconnected groups VI\jour Represent.Th.\vol8\yr2004\pages 377-413\endref

\ref\key {170} \by G.Lusztig \paper  Parabolic character sheaves II\jour Moscow Math.J.\vol4\yr2004\pages 869-896\endref

\ref\key {171} \by G.Lusztig \paper  Convolution of almost characters\jour Asian J. Math.\vol8\yr2004\pages 769-772
\endref   

\ref\key {172} \by G.Lusztig \paper  Character sheaves on disconnected groups VII\jour Represent. Th.\vol9\yr2005
\pages 209-266\endref

\ref\key {173} \by G.Lusztig \paper  Unipotent elements in small characteristic\jour Transform. Groups\vol 10\yr2005
\pages 449-487\endref

\ref\key {174} \by G.Lusztig \paper  Character sheaves and generalizations\inbook The Unity of Mathematics\lb\bookinfo
 ed. P.Etingof et al., Progress in Math.244\publ Birkh\"auser Boston\yr 2006\pages 443-455\endref

\ref\key {175} \by G.Lusztig \paper  A $q$-analogue of an identity of N.Wallach\inbook Studies in Lie theory\lb
\bookinfo ed. J.Bernstein et al., Progress in Math. 243\publ Birkh\"auser Boston\yr 2006\pages 405-410\endref

\ref\key {176} \by G.Lusztig \paper  Character sheaves on disconnected groups VIII\jour Represent. Th.\vol10\yr2006
\pages 314-352\endref

\ref\key {177} \by G.Lusztig \paper  Character sheaves on disconnected groups IX\jour Represent. Th.\vol10\yr2006
\pages 353-379\endref

\ref\key {178} \by G.Lusztig \paper  A class of perverse sheaves on a partial flag manifold\jour Represent. Th.
\vol11\yr2007\pages 122-171\endref

\ref\key {179}\by  X.He and  G.Lusztig \paper  Singular supports for character sheaves on a group compactification\jour
 Geom. and Funct.Analysis 
\vol17\yr2007\pages 1915-1923\endref

\ref\key {180} \by G.Lusztig \paper  Irreducible representations of finite spin groups\jour Represent. Th.\vol12\yr
2008\pages 1-36\endref

\ref\key {181} \by G.Lusztig \paper  A survey of total positivity\jour Milan J. Math.\vol76\yr2007\pages 1-10\endref

\ref\key {182} \by G.Lusztig \paper  Generic character sheaves on disconnected groups and character values
\jour Represent. Th.\vol12\yr2008\pages 225-235
\endref

\ref\key {183} \by G.Lusztig \paper  Unipotent elements in small characteristic II\jour Transform. Groups\vol 13
\yr2008\pages 773-797\endref

\ref\key {184} \by G.Lusztig \paper  Study of a $\ZZ$-form of the coordinate ring of a reductive group\jour 
J. Amer. Math. Soc.\vol22\yr2009\pages 739-769\endref

\ref\key {185} \by S.Kumar, G.Lusztig and D.Prasad\paper Characters of simplylaced nonconnected groups versus 
characters of nonsimplylaced connected groups\inbook Representation theory\bookinfo ed. Z.Lin, Contemp. Math. 478
\yr2009\pages 99-101\endref

\ref\key {186} \by G.Lusztig \paper  Twelve bridges from a reductive group to its Langlands dual\inbook 
Representation theory\bookinfo ed. Z.Lin, Contemp. Math.478\yr2009\pages 125-143\endref

\ref\key {187} \by G.Lusztig \paper  Character sheaves on disconnected groups X\jour Represent. Th.\vol13\yr2009
\pages 82-140\endref

\ref\key {188} \by G.Lusztig \paper  Unipotent classes and special Weyl group representations\jour J. Alg.\vol321
\yr2009\pages 3418-3449\endref

\ref\key {189} \by G.Lusztig \paper  Remarks on Springer's representations\jour Represent. Th.\vol13\yr2009\pages 
391-400\endref

\ref\key {190} \by G.Lusztig \paper  Notes on character sheaves\jour Moscow Math.J.\vol9\yr2009\pages 91-109\endref

\ref\key{191} \by G.Lusztig \paper  Graded Lie algebras and intersection cohomology\inbook Representation theory of 
algebraic groups and quantum groups\bookinfo ed. A.Gyoja et al., Progress in Math.284\publ Birkh\"auser\yr 2010
 \pages191-224\endref

\ref\key {192} \by G.Lusztig \paper  Unipotent elements in small characteristic IV\jour Transform. Groups\vol 14
\yr2010\endref

\ref\key {193} \by G.Lusztig \paper  Parabolic character sheaves III\jour Moscow Math.J.\vol 10\yr2010\pages 603-609
\endref

\ref\key{194} \by G.Lusztig \paper  Unipotent elements in small characteristic III\jour J. Alg.\vol329\yr2011\pages 
163-189\endref

\ref\key {195} \by G.Lusztig \paper  Piecewise linear parametrization of canonical bases\jour Pure Appl. Math. Quart.
\vol7\yr2011\pages 783-796\endref

\ref\key {196} \by G.Lusztig \paper  On some partitions of a flag manifold\jour Asian J. Math.\vol 15\yr2011\pages 
1-8\endref

\ref\key{197} \by G.Lusztig \paper  From conjugacy classes in the Weyl group to unipotent classes\jour Represent.Th. 
\vol15\yr2011\pages 494-530\endref

\ref\key {198} \by G.Lusztig \paper  From groups to symmetric spaces\jour Contemp. Math.\vol557\yr2011
\pages 245-258\endref

\ref\key {199} \by G.Lusztig \paper  Study of antiorbital complexes\jour Contemp. Math.\vol557\yr2011\pages 259-287\endref

\ref\key{ 200} \by G.Lusztig \paper  On C-small conjugacy classes in a reductive group\jour Transfor. Groups\vol 16
\yr2011\pages 807-825\endref

\ref\key {201} \by G.Lusztig\paper Bruhat decomposition and applications, arxiv:1006.5004\endref

\ref\key {202} \by G.Lusztig\paper On certain varieties attached to a Weyl group element\jour Bull. Inst. Math. Acad.
Sinica (N.S.)\vol 6\yr2011 
\pages377-414\endref

\ref\key {203} \by X.He and G.Lusztig \paper A generalization of Steinberg's cross-section\jour J. Amer. Math. Soc.
\vol 25\yr2012\pages 739-757\endref

\ref\key {204} \by G.Lusztig \paper  Elliptic elements in a Weyl group: a homogeneity property\jour Represent. Th. 
\vol16\yr2012\pages 127-151\endref

\ref\key {205} \by G.Lusztig \paper  From conjugacy classes in the Weyl group to unipotent classes II\jour Represent.
Th.\vol 16\yr2012\pages 189-211\endref

\ref\key {206} \by G.Lusztig \paper  On the cleanness of cuspidal character sheaves\jour Moscow Math.J.\vol 12\yr2012
\pages621-631\endref

\ref\key {207 } \by G.Lusztig and T.Xue\paper Elliptic Weyl group elements and unipotent isometries with $p=2$\jour
 Represent. Th.\vol 16\yr2012\pages 
270-275\endref

\ref\key {208} \by G.Lusztig and D.Vogan\paper Hecke algebras and involutions in Weyl groups\jour Bull. Inst. Math. 
Acad. Sinica(N.S.) \vol7\yr2012
\pages323-354\endref

\ref\key {209} \by G.Lusztig \paper  A bar operator for involutions in a Coxeter group\jour Bull. Inst. Math. Acad.
Sinica (N.S.)\vol 7\yr2012\pages 
355-404\endref

\ref\key {210} \by G.Lusztig \paper  From conjugacy classes in the Weyl group to unipotent classes III\jour Represent.
Th.\vol16\yr2012\pages450-488\endref

\ref\key {211} \by G.Lusztig \paper  On the representations of disconnected reductive groups over $F_q$\inbook
"Recent developments in Lie Algebras, Groups and Representation theory\bookinfo ed. K.Misra, Proc. Symp. Pure Math.
\publ Amer. Math. Soc.\vol86\yr2012\endref

\ref\key {212} \by G.Lusztig and Z.Yun\paper A (-q)-analogue of weight multiplicities\jour Jour. Ramanujan Math. Soc.
\vol29A\yr2013\pages311-340\endref

\ref\key {213} \by J.-L.Kim and G.Lusztig\paper On the characters of unipotent representations of a semisimple p-adic 
group\jour Represent. Th.\vol17\yr2013\pages426-441\endref

\ref\key {214} \by G.Lusztig \paper  Asymptotic Hecke algebras and involutions\inbook Perspectives in Representation
Theory\bookinfo ed. P.Etingof et.al., Contemp.Math.610 \yr2014\pages  267-278\endref

\ref\key {215} \by G.Lusztig \paper  Families and Springer's correspondence\jour Pacific J.Math.\vol267\yr2014
\pages 431-450\endref

\ref\key {216} \by G.Lusztig \paper  Restriction of a character sheaf to conjugacy classes\jour Bulletin Math\'em. 
\vol58\yr2015\pages 297-309\endref

\ref\key {217}\by J.-L.Kim and G.Lusztig\paper On the Steinberg character of a semisimple p-adic group\jour Pacific 
J.Math.\vol 265\yr2013
\pages499-509\endref

\ref\key {218} \by G.Lusztig and D.Vogan\paper Quasisplit Hecke algebras and symmetric spaces\jour Duke Math. J. 
\vol163\yr2014\pages983-1070\endref

\ref\key {219} \by G.Lusztig \paper  Unipotent almost characters of simple $p$-adic groups\jour Ast\'erisque
\vol 369-370\yr2015 \pages243-267\endref

\ref\key {220} \by G.Lusztig \paper  Unipotent almost characters of simple $p$-adic groups, II\jour Transfor. Gr. 
\vol19\yr2014 \pages527-547\endref

\ref\key {221} \by G.Lusztig \paper  Distinguished conjugacy classes and elliptic Weyl group elements\jour
 Represent.Th.\vol 18\yr2014\pages223-277\endref

\ref\key {222} \by G.Lusztig \paper   On conjugacy classes in a reductive group\inbook Representations of Reductive 
Groups\bookinfo Progr.in Math. 312\publ Birkh\"auser\yr 2015\pages 333-363\endref

\ref\key {223} \by G.Lusztig \paper Truncated convolution of character sheaves\jour
 Bull. Inst. Math. Acad. Sinica (N.S.)\vol10\yr2015\pages1-72\endref

\ref\key{224} \by G.Lusztig \paper  On conjugacy classes in the Lie group $E_8$, arxiv:1309.1382\endref

\ref\key {225} \by G.Lusztig and D.Vogan\paper Hecke algebras and involutions in Coxeter groups\inbook 
Representations of Reductive Groups\bookinfo Progr.in Math. 312\publ Birkh\"auser\yr 2015\pages 365-398\endref

\ref\key {226} \by G.Lusztig\paper Unipotent representations as a categorical centre\jour Represent.Th.\vol19
\yr2015\pages211-235\endref

\ref\key {227} \by G.Lusztig \paper  Exceptional representations of Weyl groups \jour J. Alg.\vol475\yr2017\pages14-20
\endref

\ref\key {228} \by G.Lusztig \paper  Action of longest element on a Hecke algebra cell module\jour
Pacific.J.Math. \vol279\yr2015\pages383-396\endref

\ref\key {229} \by G.Lusztig \paper  On the character of certain irreducible modular representations\jour
 Represent. Th.\vol 19\yr2015\pages 3-8\endref

\ref\key {230} \by G.Lusztig \paper  Algebraic and geometric methods in representation theory, arxiv:1409.8003\endref

\ref\key {231} \by G.Lusztig \paper  Some power series involving involutions in Coxeter groups\jour  Repres.Th.\vol19
\yr2015\pages281-289\endref

\ref\key {232} \by G.Lusztig \paper  Nonsplit Hecke algebras and perverse sheaves\jour  Selecta Math.\vol22\yr 2016
\pages1953-1986\endref

\ref\key {233} \by G.Lusztig and G.Williamson\paper On the character of certain tilting modules
\jour Science China\endref

\ref\key {234} \by G.Lusztig \paper  Non-unipotent character sheaves as a categorical centre\jour  Bull. Inst. Math.
Acad. Sinica (N.S.)\vol 11\yr 2016\pages603-731 \endref

\ref\key {235} \by G.Lusztig \paper  An involution based left ideal in the Hecke algebra\jour  Represent .Th.\vol20
\yr2016\pages172-186 \endref

\ref\key {236} \by G.Lusztig \paper  Generic character sheaves on groups over $\kk[\e]/(\e^r)$ \jour Contemp. Math.
\vol683\yr2017\pages227-246 \endref

\ref\key {237} \by G.Lusztig \paper  Generalized Springer theory and weight functions\jour Ann. Univ. Ferrara Sez.VII
Sci. Mat. \vol63\yr2017\pages159-167 \endref

\ref\key {238} \by G.Lusztig \paper  On the definition of almost characters \jour
Transfor.Groups\toappear \endref

\ref\key {239} \by G.Lusztig \paper Special representation of Weyl groups: a positivity property
\jour Adv.in Math.\toappear  \endref

\ref\key {240} \by G.Lusztig and Z.Yun\paper $\ZZ/m$-graded Lie algebras and perverse sheaves I\jour
Represent.Th.\vol21\yr2017\pages277-321\endref

\ref\key {241} \by G.Lusztig \paper The canonical basis of the quantum adjoint representation\jour J. Comb. Alg.
\vol1\yr2017\pages45-57\endref

\ref\key {242} \by G.Lusztig and  Z.Yun\paper $\ZZ/m$-graded Lie algebras and perverse sheaves II\jour
Represent.Th.\vol21\yr2017\pages322-353\endref

\ref\key {243} \by G.Lusztig and Z.Yun\paper $\ZZ/m$-graded Lie algebras and perverse sheaves III: graded double 
affine Hecke algebra\jour Represent.Th.\vol22\yr2018\pages87-118\endref

\ref\key {244} \by G.Lusztig\paper On the generalized Springer correspondence\jour arxiv:1608:02222\toappear
\endref

\ref\key {245}\by G.Lusztig\paper Non-unipotent representations and categorical centers\jour 
Bull. Inst. Math. Acad. Sinica (N.S.)\vol12\yr2017\pages205-296\endref

\ref\key{246}\by G.Lusztig and G.Williamson \paper Billiards and tilting characters of $SL_3$
\jour SIGMA Symmetry, Integrability, Geometric Methods Appl.\yr2018\endref

\ref\key{247}\by G.Lusztig \paper Conjugacy classes in reductive groups and two-sided cells\jour
Bull. Inst. Math. Acad. Sinica\toappear\endref

\ref\key{248}\by G.Lusztig \paper Comments on my papers\jour arxiv:1707.09368\endref

\ref\key{249}\by G.Lusztig\paper Lifting involutions in a Weyl group to the
torus normalizer\jour Represent.Th.\yr2018\endref

\ref\key{250}\by G.Lusztig\paper Hecke modules based on involutions in extended Weyl groups\jour 
arxiv:1710.03670\finalinfo submitted\endref

\ref\key{251}\by G.Lusztig\paper Discretization of Springer fibres\jour arxiv:1712.07530\endref

\ref\key{252}\by G.Lusztig\paper  A new basis for the representation ring of a Weyl 
group\jour arxiv:1805.03770\finalinfo submitted\endref

\ref\key{253}\by G.Lusztig\paper Positive conjugacy classes in Weyl groups\jour
arxiv:1805.03772\finalinfo to appear\endref

\ref\key{254} \by G.Lusztig and Z.Yun\paper $\ZZ/m$-graded Lie algebras and
perverse sheaves IV\jour arxiv:1805.10550\endref

\ref\key{255}\by G.Lusztig\paper A new basis for the representation ring of a Weyl 
group, II\jour arxiv:1808.06896\endref
\endRefs
\enddocument